\documentclass{amsart}
\usepackage{amssymb}
\usepackage[all]{xy}
\xyoption{dvips}
\newcommand{\pair}[1]{\langle #1 \rangle}
\newcommand{\bbox}{\openbox}
\newcommand{\Cole}{{\textup{Col}}}
\newcommand{\Cone}{\operatorname{Cone}}
\newcommand{\dlog}{\operatorname{\dd\!\log}}
\newcommand{\dotcup}{\cup}
\newcommand{\spf}{\operatorname{Spf}}

\newcommand{\dr}{{\textup{dR}}}
\newcommand{\hdr}{H_{\dr }}
\newcommand{\hr}{H_{\rig }}
\newcommand{\htms}{\thmsyn}
\newcommand{\id}{\operatorname{id}}
\newcommand{\im}{\operatorname{im}}
\newcommand{\inject}{\hookrightarrow}
\newcommand{\isom}{\cong}
\newcommand{\Li}{\textup{Li}}
\newcommand{\syn}{{\textup{syn}}}
\newcommand{\ms}{{\textup{ms}}}
\newcommand{\rig}{{\textup{rig}}}
\newcommand{\thmsyn}{\tilde{H}_{\ms}}
\newcommand{\hms}{H_{\ms}}
\newcommand{\hsyn}{H_{\syn}}
\newcommand{\loc}{{\text{loc}}}

\newcommand{\Xbar}{\overline{X}}
\newcommand{\Bbar}{\overline{B}}
\newcommand{\Cbar}{\overline{C}}
\newcommand{\Fbar}{\overline{\F}}
\newcommand{\Ybar}{\overline{Y}}
\newcommand{\Phat}{\widehat{P}}

\newcommand{\Xbarhat}{\widehat{\Xbar}}

\renewcommand{\O}{\mathcal{O}}
\newcommand{\bcat}{\mathcal{S}}
\newcommand{\smap}{r}
\newcommand{\calC}{\mathcal{C}}
\newcommand{\PF}{{\mathcal{P}}}
\newcommand{\vv}{{R}}
\newcommand{\PP}{{\mathbb{P}}}
\def\dirlim_#1{\mathchoice{\underset{#1}{\varinjlim}}{\varinjlim_{#1}}{}{}}
\newcommand{\rgrig}{{\mathbb{R}\Gamma_\rig}}

\newcommand{\rgtms}{\widetilde{\mathbb{R}\Gamma}_{\ms}}
\newcommand{\rg}{\mathbb{R}\Gamma}
\newcommand{\tu}[2]{{]#1[_{#2}}}
\newcommand{\tub}[1]{{]#1[}}
\newcommand{\jdag}{j^\dagger}
\newcommand{\rreg}{\operatorname{\mathcal{R}}}
\newcommand{\Ulocb}[1]{U_{B,\loc}^{#1}}
\newcommand{\Uloc}[1]{U_{\loc}^{#1}}
\def\a{\alpha}
\def\ab #1 {\square(#1)}
\def\alt{{\fam\fcrm alt}}
\def\b{\beta}
\def\calC{{\mathcal C}}
\def\calfinalC{{\mathcal C}_{\square}}
\def\Cr #1 #2 #3 #4 {{C\left(#1, \{ #2_{#3},\dots, #2_{#4}\}\right) }}
\def\Crshort #1 #2 #3 #4 {{C\left(#3,\dots,#4\}\right) }}
\def\cSymb #1 {{\mathcal Symb}_{(#1)}^\bullet}
\def\Coprod{\displaystyle \coprod}
\def\dd{\textup{d}}

\def\FQ{{F^*_\Q}}

\def\I{(1+I)^*}
\def\J #1 {{\mathcal J}_{(#1)}^\bullet}
\def\iso{\cong}
\def\isoot{{\overset \sim \leftarrow}}
\def\K #1 #2 #3 {K_{#1}^{(#2)}(#3)}
\def\KQ{K_\Q^*}
\def\Krloc #1 #2 #3 #4 {\K #1 #2 {\xlocb #3 #4 } }
\def\Krel #1 #2 #3 #4 {\K #1 #2 {\xb #3 #4 } }
\def\Ks #1 #2 #3 #4 {{\scriptstyle K_{#1, #3}^{(#2), #4}}}
\def\Ksnot #1 #2 #3 #4 {{ K_{#1, #3}^{(#2), #4}}}
\def\L{L}
\def\Lmod #1 {L_{\mdfd,#1}}
\def\mdfd{{\fam\fcrm mod}}
\def\N #1 {{\mathcal N}_{(#1)}^\bullet}
\def\OQ{\O_\Q^*}
\def\Osp{\O^\flat} % the special units
\def\reg{{\fam\fcrm reg}}
\def\s{\sigma}
\def\scp{\textstyle\coprod}
\def\scriptM #1 {{\mathcal M}_{(#1)}^\bullet}
\def\Span #1 {{<}#1{>}}
\def\Spec{\textup{Spec}}
\def\Sym{{\fam\fcrm Sym}}
\def\symb #1 #2 {[#1]_{#2}}
\def\Symb #1 #2 {{\fam\fcrm Symb}_{#1}(#2)}
\def\tcup{\tilde\cup}
\def\tensor{\otimes}
\def\tildescriptM #1 {\widetilde{\mathcal M}_{(#1)}^\bullet}
\def\tM{\widetilde M}
\def\w{\wedge}
\def\W{\bigwedge}
\def\xb #1 #2 {{X_{#1}^{#2};\bbox^{\,#2}}}
\def\xl #1 #2 {X_{#1,\loc}^{#2}}
\def\xloc #1 {{X_{\loc}^{#1};\bbox^{\,#1}}}
\def\xlocb #1 #2 {{X_{#1,\loc}^{#2};\bbox^{\,#2}}}
\def\loc{{\fam\fcrm loc}}
\def\A{{\mathbb A}}
\def\C{{\mathbb C}}
\def\F{{\mathbb F}}

\def\Q{{\mathbb Q}}

\def\Z{{\mathbb Z}}
\def\soule{Soul\'e}

\begin{document}

\newfam\fcrm
\font\tenrm=cmr10
\font\sevenrm=cmr7
\font\fiverm=cmr5
\textfont\fcrm\tenrm
\scriptfont\fcrm\sevenrm
\scriptscriptfont\fcrm\fiverm

\newtheorem{philos}[equation]{Philosophy}
\newtheorem{thm}[equation]{Theorem}
\newtheorem{conjecture}[equation]{Conjecture}
\newtheorem{notation}[equation]{Notation}
\newtheorem{assumption}[equation]{Assumption}
\newtheorem{proposition}[equation]{Proposition}
\newtheorem{prop}[equation]{Proposition} %%Delete [thm] to re-start numbering
\newtheorem{lemma}[equation]{Lemma}
\newtheorem{corollary}[equation]{Corollary}
\theoremstyle{definition}
\newtheorem{definition}[equation]{Definition}
\newtheorem{remark}[equation]{Remark}
\newtheorem{example}[equation]{Example}
\theoremstyle{plain}
\numberwithin{equation}{section} %% Comment out for sequentially-numbered
\numberwithin{figure}{section} %% Comment out for sequentially-numbered

\title{The syntomic regulator for $K$--theory of fields}
\author{Amnon Besser}
\address{
Department of Mathematics\\
Ben-Gurion University of the Negev\\
P.O.B. 653\\
Be'er-Sheva 84105\\
Israel
}
\author{Rob de Jeu}
\address{Department of Mathematical Sciences\\
University of Durham\\
Science Laboratories, South Road\\
Durham DH1 3LE\\
United Kingdom}
\begin{abstract}
We define complexes analogous to Goncharov's complexes for the
$ K $--theory of discrete valuation rings of characteristic zero.
Under suitable assumptions in $ K $--theory, there is a map from
the cohomology of those complexes to the $ K $--theory of the
ring.  In case the ring is the localization of the ring of integers
in a number field, there are no assumptions necessary.  We compute
the composition of our map to the $ K $--theory with the syntomic
regulator.  The result can be described in terms of a $ p $--adic
polylogarithm.  Finally, we apply our theory in order to compute
the regulator to syntomic cohomology on Beilinson's cyclotomic elements.
The result is again given by the $ p $--adic polylogarithm.
This last result is related to one by Somekawa and generalizes
work by Gros.
\end{abstract}

\maketitle

\section{Introduction}

Let $ K $ be a complete discrete valuation field of
characteristic zero, $ R $ its valuation ring,
and $ \kappa $ its residue field. Assume $ \kappa $ has positive
characteristic $ p $ and is algebraic over $\F_p$.
If $ X/R $ is smooth, separated
and of finite type, there is a regulator
map from K-theory to syntomic cohomology
\[
\K n j X \to \hsyn^{2j-n}(X, j)\;,
\]
see \cite{Bes98a}. In many interesting cases
the target group of the regulator is isomorphic to the
rigid cohomology group, in the sense of Berthelot,
$\hr^{2j-n-1}(X_\kappa/K)$, where $X_\kappa$ is the special fiber of $X$.
We will be most interested in the situation where $ X = \textup{Spec}(R) $,
and the $ K $--group is $ \K 2n-1 n R $ for $ n \geq 2 $.  The target group
for the regulator in this case is $\hr^0(\Spec(\kappa)/K)\isom K$
(see Definition~\ref{normalization} for the precise identification).
Because $ \kappa $ is algebraic over $ \F_p $,  $ K_n(\kappa) $ is torsion
for all $ n \geq 1 $, so from the long exact localization sequence
\[
\dots \to \K n j-1 {\kappa} \to \K n j R \to \K n j K \to \K n-1 j-1 {\kappa} \to \cdots
\]
we get an isomorphism $ \K 2n-1 n R \iso \K 2n-1 n K $ for $ n\geq 2 $.
Hence we get a regulator map (for $ n \geq 2 $)
\[
\reg: \K 2n-1 n K \iso \K 2n-1 n R \to K\;.
\]
In this paper we try to explicitly compute this regulator map.
We note that if $ F $ is a number field
with an embedding $ F \to K $, we can combine the natural map
$ \K 2n-1 n F \to \K 2n-1 n K $
with this regulator map to obtain a regulator map on $ \K 2n-1 n F $.
Also,  for a number field $ F $, all $ K_n(F) $ are torsion
if $ n $ is even and positive.  For the odd ones, all
$ K_{2n-1}(F) \tensor_\Z \Q $
are $ \K 2n-1 n F $, so the computation for $\K 2n-1 n K $ is the most
interesting from the point of view of number fields.
Our principal tool of study will be the complexes $ \tildescriptM n (K) $
constructed in \cite{Jeu95}.
Write $ \KQ$ for $ K^*\tensor_\Z \Q $.
The complex $ \tildescriptM n (K) $
for $ n\geq 2 $ is of the form
\[
\tM_n \to \tM_{n-1} \tensor \KQ \to \tM_{n-2} \tensor \W^2 \KQ \to
\dots \to
\tM_2 \tensor \W^{n-2} \KQ \to \tensor \W^n \KQ
\]
where $ \tM_k = \tM_k(K) $ is a $ \Q $--vector space generated
by symbols $ \symb x k $, with $ x $ in $ K $, $ x \neq 0 $ or
1, and the differential is given by
\[
\dd (\symb x k \tensor y_1\w\dots\w y_{n-k}) =
     \symb x k-1 \tensor x\w y_1\w\dots\w y_{n-k}
\]
if $ k \geq 3 $, and 
\[
\dd (\symb x 2 \tensor y_1\w\dots\w y_{n-k}) =
        (1-x)\w x\w y_1\w\dots\w y_{n-k}\;.
\]
We give this complex a cohomological grading in degrees 1 through
$ n $.
Under suitable assumptions about weights in $K$--theory 
(as formulated in the Beilinson--\soule\ conjecture, see
Definition~\ref{no-low-weight-def}), there is a map
\begin{equation}
  \label{eq:kbasicmapa}
  H^r(\tildescriptM n (K)) \to \K 2n-r n K \;.
\end{equation}

In Section~\ref{sec:kstuff}, we construct analogous complexes
$ \tildescriptM n (R) $ for the
ring $ R $,
whose cohomology (again under suitable assumptions) maps directly to
the $ K $--theory
of $ R $, and in Section~\ref{sec:end} we compute the regulator map on
its image. 
In the cases we are interested in $ \tildescriptM n (R) $
can be identified with the subcomplex of the
complex for $ K $ spanned in degree $ k+1 = 1,\dots,n-1 $ by all
$ \symb u n-k \tensor v_1 \wedge \dots \wedge v_k $
with $ v_1,\dots,v_k $ in $ R^* $, $ u $ 
in $ R^* $ such that $ 1-u $ is also in $ R^* $, and
in degree $ n $ by all $ v_1 \w \dots \w v_n $ with all $ v_i $
in $ R^* $.  (A $ u $ in $ R^* $ such that $ 1-u $ is also in
$ R^* $ is called a \emph{special unit}.)  But redoing the construction has the advantage that
we can work over $R$ all the time, which is required for the computation
of the regulator.

The case that the field is a number field deserves special mentioning.
First of all, no assumptions about weights are necessary in
this case.  Furthermore,
it is known that if $ F $ is a number field, the map
$ H^1(\tildescriptM n (F) ) \to \K 2n-1 n F $ is an isomorphism
for $ n = 2 $ and $ n = 3 $, as well as when $ F $ is a cyclotomic
field for all $ n \geq 2 $.  (There is also substantial numerical
evidence that it should be an isomorphism for all $ n $ for number
fields, which is part of a conjecture by Zagier, as well as a
corresponding conjecture for infinite fields by Goncharov.)
Therefore one would have a complete description of the
syntomic regulator for our discrete valuation ring $R$
if we knew that the image of $ H^1(\tildescriptM n (R)) $
in $ \K 2n-1 n R $ would be everything.  This may not be the
case, as perhaps $ R $ does not have enough special units. One can try
to overcome this difficulty by rewriting elements
in the image of $ H^1(\tildescriptM n (F)) $ as being part of
the image of $ H^1(\tildescriptM n (R')) $ where $ F'/F $ is
a finite field extension, $ R' $ the corresponding ring in $ F' $.
We do this in the case of cyclotomic fields, so that we obtain
a full description of the syntomic regulator in this case.
We also state a Conjecture that the formulas found for the regulator
on the complex for $ R $ generalize to be the regulator on the
complex for $ F $.

In order to present our results,
we shall need the following functions.  Let $ \log : \C_p^* \to \C_p $
be a branch of the $p$--adic logarithm.  This means we define
$ \log $ on the elements $ z $ with
$ |1-z| <1 $ by the usual power series, and we extend this to
$ \C_p^* $ by choosing $\pi\in \C_p$ with $|\pi|<1$, declaring
$\log(\pi)=0$,  and extending to a homomorphism
from $ \C_p^* $ to $ \C_p $ (see Definition~\ref{branchlog}).  Note
that the values on the elements
in $ \C_p^* $ with $ |z|=1 $ is independent of the choice of
$\pi$, but $ \log $ and the functions $ \Li_n(z) $ about
to be described depend on this choice.  For the relation, see
Proposition~\ref{branch-prop}.

Let $ \Li_1(z) = -\log(1-z) $ for $ z \neq 0 $ or 1. We follow Coleman
to recursively define, using his integration theory, functions
$ \Li_n(z) $ for $ n \geq 2 $. The defining relations are
$ \dd \Li_n(z) = \Li_{n-1}(z) \dlog z $
and $ \lim_{z\to 0} \Li_n(z) = 0 $, and they have a unique solution in
the class of functions defined by Coleman. It is shown in \cite{Col82}
that those functions are locally analytic in the naive topology on $\C_p$,
and are given by a
convergent power series $ \sum_{k=1}^\infty z^k/k^n $
on the open unit disc in $ \C_p $.  The function $ \Li_n(z) $ extends
to a locally analytic function on $ \C_p \setminus \{1\} $ with
$ \Li_n(z) = 0 $ for $ n \geq 1 $.
These functions satisfy the functional equation
\begin{equation}
\label{functionaleq}
\Li_n(z) + (-1)^n \Li_n(1/z) = \frac{-1}{n!} \log^n(z)\;,
\end{equation}
see Proposition~6.4 of \cite{Col82}.
We also introduce the function $ \L_n $, defined as
\begin{equation}\label{Ldef}
\L_n(z) = \sum_{m=0}^{n-1}\frac{(-1)^m}{m!}\Li_{n-m}(z) \log^m(z)\;.
\end{equation}

In order to state the Theorems below easily, we shall need
linear combinations of these functions.  Namely,
we want a suitable combination that satisfies a clean functional
equation for $ z $ and $ 1/z $.
It follows from~\eqref{functionaleq} that
$ L_k(z) + (-1)^k L_k(1/z) = \frac{(-1)^k}{k!} \log^k(z) $.  Therefore
the function
\begin{equation*}
  \Lmod n (z) = \sum_{m=1}^n a_m L_m(z) \log^{n-m}(z)
\end{equation*}
with $ a_n=1 $ satisfies 
\begin{equation}
  \label{zooz}
  \Lmod n (z) + (-1)^n \Lmod n (1/z) = 0
\end{equation}
if $ \sum_{m=1}^n a_m \frac{(-1)^m}{m!} = 0 $.
Below, $ \Lmod n $ will mean any of those choices.
For $ n=2 $, there is a unique such function, namely
\[
L_2(z)+\frac12 \log(z)L_1(z) =
\Li_2(z)- \frac12 \log(z)\Li_1(z)\;,
\]
which is studied in Section~VI and beyond in \cite{Col82}, where
it is called  $ D(z) $.

It is easily deduced from
Coleman's theory (see Remark~\ref{galequiv}) that $\Li_n$ is Galois
equivariant. In particular, if $K\subset \C_p$ is a complete subfield, then
$\Li_n$, and as a result also $L_n$ and $\Lmod n $,
send $K$ to $K$ provided $\log$ was defined such that $\log(\pi)=0$
with $\pi\in K$.

\begin{remark}
By considering the coefficients of the terms $ \Li_m(z)\log^{n-m}(z) $
in the functions above, one sees that the functions
$ \L_m(z)\log^{n-m}(z) $ for $ m=1,\dots,n $ and
$ \Lmod m (z)\log^{n-m}(z) $ as above for $ m=1,\dots,n $ span
the same $ \C_p $--vector space (or even $ \Q $--vector space
in case all $ a_i $'s are in $ \Q $), namely the space spanned
by $ \Li_m(z)\log^{n-m}(z) $ for $ m=1,\dots,n $.  
Therefore one can consider any function of the form
$ \sum_{j=0}^{n-1} b_j \Li_{n-j}(z)\log^j(z) $
with all $ b_j $ in $ \C_p $ as a candidate for $ \Lmod n $, provided
$ b_0=1 $ and $ \sum_{j=0}^{n-1} \frac{b_j}{(n-j)!} =0 $.
Let $ B_i $ for $ i=0,1,\dots $ be the Bernoulli numbers, defined
by the identity of formal power series
\begin{equation*}
  \sum_{i=0}^\infty \frac{B_i}{i!} t^i = \frac{t}{e^t-1}\; .
\end{equation*}
Then the functions $ \Lmod n (z) $
defined by  $ \sum_{j=0}^{n-1} \frac{B_j}{j!} \Li_{n-j}(z) \log^j(z) $
satisfy the above requirements as $ B_0=1 $, and the other identity
holds by definition of the $ B_j $ if $ n\geq 2 $.
Note that this formula is different from the classical case,
where one uses the real or imaginary part of the functions
$ \sum_{j=0}^{n-1} \frac{B_j 2^j}{j!} \Li_{n-j}(z) \log^j|z| $,
see \cite{Zag91} and \cite[Remark~5.2]{Jeu95}. Another possible
natural candidate for the function $\Lmod m (z)$ is $\Lmod m (z)=
L_m(z) + L_{m-1}(z) \log(z)/m$. This function is distinguished by the
following fact proved in~\cite[Theorem~\ref{thetheorem}]{Bes00a}: It
is the unique combination of
type $\Lmod m $ with coefficients independent of $p$ such that the
function $-m p^{1-m}z (1-z) \dd \Lmod m (z)/\dd z$ has a reduction modulo $p$,
for sufficiently large $p$,
which is the so called $m-1$-polylogarithm function introduced in the
$m=2$ case by Kontsevich~\cite{Kon99} and by Elbaz-Vincent and
Gangl~\cite{Elb-Gan99} in general.
\end{remark}

If $ R $ is a ring with 1, let $ R^\flat $ be the set of elements  $ u $ in $ R $
such that both $ u $  and $ 1-u $ are units.  We shall refer to those elements
as \emph{special units}.

We are now ready to state our main results.

\begin{thm}
\label{general-field}
Let $ F $ be a field of characteristic zero.  Let $ \O \subset F $
be a discrete valuation ring, and let $ \kappa $ be the residue field.
Assume that the Beilinson--\soule\
conjecture holds for fields of characteristic 0 and for $ \kappa $.
For $ n\geq 2 $ let $ \tildescriptM n (\O) $ be
the subcomplex of the complex $ \tildescriptM n (F) $ constructed
in \cite{Jeu95} (see also Section~\textup{\ref{sec:kstuff}})
generated by symbols of
the form $ \symb x k \tensor y_1 \w\dots\w y_{n-k} $, where all
$ y_i $ are elements in $ \O^* $, and $ x $ is in $ \O^\flat $.
Then
\begin{enumerate}
\item  There is a map
$ H^r(\tildescriptM n (\O)) \to \K 2n-r r {\O} $ such that the diagram
\[
\xymatrix{
H^r(\tildescriptM n (\O)) \ar[r]\ar[d] & \K 2n-r n {\O} \ar[d]
\\
H^r(\tildescriptM n (F)) \ar[r]        & \K 2n-r n F
}
\]
commutes, where the lower horizontal map is the map in \eqref{eq:kbasicmapa}.
\item If in addition $ \sigma : F \to K $ is an embedding with $ \sigma(\O) \subset R $,
then for $ r = 1 $, the regulator map
\[
H^1( \tildescriptM n (\O) )\to \K 2n-1 n {\O} \overset{\sigma}{\to} \K 2n-1 n R \to K
\]
is given by mapping $ \symb x n $ to $  \pm  (n-1)! \, \Lmod n (\sigma(x)) $.
\end{enumerate}
Moreover, if $ n=2 $, those results hold without
any assumptions on the Beilinson-\soule\ conjecture.
\end{thm}

\begin{remark}\label{sign-indeterminacy}
The indeterminacy of the sign in Theorem~\ref{general-field}
comes from the fact that the maps
$ H^r(\tildescriptM n (\O)) \to \K 2n-r n {\O} $ and
$ H^r(\tildescriptM n (F)) \to \K 2n-r n F $
in Theorem~\ref{general-field} are
natural only up to a choice of sign (see the introduction to Section~\ref{sec:kstuff}).
This indeterminacy of the
sign will occur every time those maps are referred to in the paper.
\end{remark}

\begin{remark}
Our computations will show that there is a 
map $ \tM_n(R) \to K $
given by mapping $ \symb x n $ to $ \pm  (n-1)!\, \Lmod n (x) $,
and that this map is compatible with the regulator map
$\K 2n-1 n R \to K $ if the assumptions in
Theorem~\ref{general-field} are fulfilled.
\end{remark}

\begin{remark}
In the exact localization sequence
\[
\dots\to
\K 2n-r n-1 {\kappa} \to
\K 2n-r n {\O} \to
\K 2n-r n F \to
\K 2n-r-1 n-1 {\kappa} \to
\cdots
\]
we have that $ \K 2n-1 n-1 {\kappa} $ and $ \K 2n-r-1 n-1 {\kappa} $
are both zero if $ r = 1 $ and $ n \geq 2 $
because we are assuming that the Beilinson--\soule\
conjecture holds for $ \kappa $.  Hence for $r=1$ the map
$ \K 2n-1 n {\O} \to \K 2n-1 n F $
in Theorem~\ref{general-field} above is an isomorphism.
Note that if an embedding $ \sigma : F \to K $ exists as in the
Theorem,  this implies that $ \kappa $ is algebraic over $ \F_p $,
so all $ K_m(\kappa) $ are torsion for $ m \geq 1 $.
In particular, in that case we have an isomorphism
$ \K 2n-r n {\O} \iso \K 2n-r n F $ for all $ n \geq 2 $ and $ r=1,\dots,2n-2 $.
\end{remark}

If $ F $ is a number field, the Beilinson-\soule\ conjecture
is known for $ F $, and one can get the map
$ H^r( \tildescriptM n (F) ) \to \K 2n-r n F $ without making
assumptions.  In fact, for $ n=2 $ and $ n=3 $, as well as in
the case $ F $ is a cyclotomic field, for all $ n \geq 2 $, one
gets an isomorphism $ H^1(\tildescriptM n (F) ) \iso \K 2n-1 n F $
this way, see \cite[Theorem~5.3]{Jeu95}.  We formulate our
results for number fields separately, as there are no assumptions
involved about weights in this case.

Note that because $ \kappa $ will be a finite field in this case,
as before we get an isomorphism $ \K 2n-r n {\O} \to \K 2n-r n F $
for all $ n \geq 2 $ and $ r=1,\dots,n $.

\begin{thm}
\label{number-field}
Let $ F $ be a number field.  Let $ \O $ be a localization of the ring
of integers of $ F $ at a nonzero prime ideal.
Then
\begin{enumerate}
\item
There is a map
$ H^r(\tildescriptM n (\O)) \to \K 2n-r r {\O} $ such that the diagram
\[
\xymatrix{
H^r(\tildescriptM n (\O)) \ar[r]\ar[d] & \K 2n-r n {\O} \ar[d]^*[@]{\sim}
\\
H^r(\tildescriptM n (F)) \ar[r]       & \K 2n-r n F
}
\]
commutes, where the lower horizontal map is the map in \eqref{eq:kbasicmapa}.
\item
If in addition $ \sigma : F \to K $ is an embedding with $ \sigma(\O) \subset R $,
then for $ r=1 $, the regulator map
\[
\reg_\sigma :
H^1( \tildescriptM n (\O) )\to \K 2n-1 n {\O} \overset{\sigma}{\to} \K 2n-1 n R \to K
\]
is given by mapping $ \symb x n $ to $  \pm  (n-1)! \, \Lmod n (\sigma(x)) $.
\end{enumerate}
\end{thm}

\begin{remark}
For any fixed element in $ H^1(\tildescriptM n (F)) $, all elements
involved are special units for almost all $ p $, so that this
Theorem applies to any given element in $ H^1(\tildescriptM n (F)) $
for almost all $ p $.
\end{remark}

If we try to apply Theorem~\ref{number-field} to cyclotomic
elements, i.e., to the elements $ \symb {\zeta} n $ corresponding to a
$m$-th root of unity $\zeta$ we notice that it applies directly only
if $m$ is prime to $p$ since otherwise $\zeta$ may not be
special. However, using relations between symbols (the distribution
relation, see Proposition~\ref{distribution}) we are able to prove the
following Theorem.

\begin{thm}
  \label{cyclotomic}
  Under the assumptions of Theorem~\ref{number-field}, the regulator map
  \[
  \reg_\sigma :
  H^1(\tildescriptM n (F)) \to \K 2n-1 n F \iso \K 2n-1 n {\O}
  \overset{\sigma}{\to} \K 2n-1 n R \to K 
  \]
  maps $ \symb {\zeta} n $ to
  $  \pm  (n-1)! \, \Lmod n (\sigma(\zeta)) $
  if $ \zeta $ is any root of unity in $ F^* $.
\end{thm}

\begin{remark}
Because it is known that the elements $ \symb {\zeta} n $ for $ n\geq 2 $
where $ \zeta $ runs through the primitive $ m $--th roots of unity
generate $ \K 2n-1 n \Q(\mu_m) $, this gives a complete description
of the syntomic regulator for cyclotomic fields.  Note also that
this particular result
extends the results of \cite{Gro94}, where the corresponding
result was proved only for roots of unity of order $ m $ with
$ (m,p) = 1 $ (see Th\'eor\`eme~2.22),
and is equivalent to the results of \cite{Som92}. That paper has a different
formulation, with another version of a syntomic regulator and also using a
specialized version of the polylogarithm at roots of unity. The relation
with Coleman's polylogarithm was proved by Barsky (unpublished).
The result of Gros is that the element $[\zeta]_n$ is
mapped under the syntomic regulator to $\Li_n^{(p)}(\zeta)$, where
$\Li_n^{(p)}$ is defined by
\[
  \Li_n^{(p)}(z)= \Li_n(z)- \frac{1}{p^n} \Li_n(z^p)\;.
\]
Note that the expansion of $\Li_n^{(p)}$ at $0$ is $\sum_{(k,p)=1} z^k/k^n$.
The difference between the results is caused by the different
normalizations of the regulators. One has the relation
\[
  \reg_{\textup{Gros}} = 
  \left(1-\frac{\operatorname{Frob}}{p^n}\right) \reg\;,
\]
where $\operatorname{Frob}$ is the Frobenius automorphism (The Gros
regulator is only defined for unramified fields.) From Galois
equivariance it follows that
\[
  \operatorname{Frob}(\Li_n(\zeta)) = \Li_n(\operatorname{Frob}(\zeta)) =
  \Li_n (\zeta^p)\;.
\]
The relation between the two results is therefore clear.
\end{remark}

We now state the following Conjecture.

\begin{conjecture}\label{main-conjecture}
Let $ K \subset \C_p $ be a discrete valuation subfield (i.e.,
the valuation is induced from the one on $ \C_p $).
Let $ R $ be the valuation ring of $ K $.
Assume the Beilinson--\soule\ conjecture holds in characteristic
zero if $ n \geq 3 $.
Then, for all $ n\geq 2 $, the regulator map
\[
 H^1(\tildescriptM n (K) ) \to \K 2n-1 n K \iso \K 2n-1 n R \to
 K
\]
is given by the same formula as before,
$ \symb x n $ being mapped to $ \pm  (n-1)! \, \Lmod n (x) $.
\end{conjecture}

\begin{remark}
\label{independence-of-branch}
Of course, Conjecture~\ref{main-conjecture} would imply that the map
\[
\tildescriptM n (K) \to K
\]
given by mapping $ \symb x n $ to $ \Lmod n (x) $
is well defined, as any relation among the $ \symb x n $ would give
rise to the zero element in $ H^1(\tildescriptM n (K) ) $.
Moreover, it implies that the map
\[
H^1( \tildescriptM n (K) ) \to K
\]
given by mapping $ \symb x n $ to $ \Lmod n (x) $ is independent
of the choice of the branch of the logarithm.
We shall verify this last statement under the assumption that
the map $ \tildescriptM k (K) \to K $ given by mapping $ \symb x k $
to $ \Lmod k (x) $ is well defined for all $ k \leq n $
in Proposition~\ref{independent} below,
after determining the dependence of $ \Li_n(z) $ on the choice
of the logarithm.
\end{remark}

\begin{remark}
As will be described in Section~\ref{sec:kstuff}, the complexes
$ \tildescriptM n (K) $ and $ \tildescriptM n (R) $
are quotient complexes of corresponding complexes $ \scriptM n (K) $
and $ \scriptM n (R) $, obtained by imposing the relations $ \symb x k + (-1)^k \symb 1/x k $
for all $ k\geq 2 $.  The general assumptions about the Beilinson--\soule\
conjecture are necessary in order to prove this quotient map
to be a quasi--isomorphism.  
There is a map $ H^r(\scriptM n (R) ) \to \K 2n-r n R $
assuming only the Beilinson--\soule\ conjecture for $ K $ and $ \kappa $.
Therefore, assuming only the Beilinson--\soule\
conjecture for $ K $ and $ \kappa $, we get a commutative diagram
\[
\xymatrix{
H^1(\scriptM n (R)) \ar[r]\ar[d] &
   H^1(\tildescriptM n (R))\ar[d]^-{\symb z n \mapsto \pm (n-1)! \, \Lmod n (z)}
\\
\K 2n-1 n R   \ar[r]^-{\reg}     & K.
}
\]
\end{remark}

As each of the steps in the proofs of
Theorems~\ref{general-field} and~\ref{number-field}
is fairly technical, we give a
brief outline of the main steps and where in the
paper they occur.

Using multi relative $K$--theory 
and localization (both discussed in Section~\ref{sec:kstuff}), we get a diagram
$$
\xymatrix{
K_{2n-1}^{(n)}(\O) \ar[r]^-{\iso}\ar[d] & K_n^{(n)} (X_\O^{n-1} ;\bbox^{n-1}) \ar[r]\ar[d]
& K_n^{(n)} (X_{\O,\loc}^{n-1} ;\bbox^{n-1}) \ar[r]\ar[d] & \dots
\cr
K_{2n-1}^{(n)}(\vv) \ar[r]^-{\iso} & K_n^{(n)} (X_\vv^{n-1} ;\bbox^{n-1}) \ar[r]
& K_n^{(n)} (X_{\vv,\loc}^{n-1};\bbox^{n-1}) \ar[r] & \dots
}
$$
The rows (except the first term) will be used to construct the
complexes $\tildescriptM n (\O) $ and $\tildescriptM n (R) $ in
Section~\ref{sec:kstuff}.  If the Beilinson--\soule\ conjecture
holds generally enough, there is a map $H^1(\tildescriptM n (\O) )
\to \K n n {X_\O^{n-1};\bbox} \iso \K 2n-1 n {\O} $.  On the
other hand, there is a syntomic regulator
$ \K 2n-1 n {\vv} \to K $.  Using the embedding $\O \to \vv$ as
in Theorem~\ref{general-field} gives us the map 
$$
H^1(\tildescriptM n (\O) ) \to \K 2n-1 n {\O} \to \K 2n-1 n {\vv}  \to K.
$$
Similar results will be proved if $ K $ is replaced with a number
$ F $, but without assumptions on weights in algebraic $ K $--theory.
We will also compare those complexes with the complexes
$\tildescriptM n (K) $ constructed in \cite[Section~3]{Jeu95}.

After reviewing syntomic regulators in Section~\ref{sec:syntomic},
we analyze this in Sections~\ref{sec:down} through~\ref{sec:end}
by extending the regulator map over the maps
$$
K_{2n-1}^{(n)}(\vv) \iso  K_n^{(n)} {(X_\vv^{n-1};\bbox^{n-1})}
\to
\K n n { X_{\vv,\loc}^{n-1} ;\bbox^{n-1} }
.
$$
The computation involves Coleman integration, and we start the
paper by reviewing it, and tying up some loose ends, in Section~\ref{sec:prelim}.

Finally, we would like to thank the support of the European
community through the TMR network \emph{Arithmetic Algebraic
Geometry}.

{\bf Notation:}
Throughout the paper, if $ A $ is an Abelian group, we shall write $ A_\Q $
for $ A \tensor _\Z \Q $.
Throughout the paper, $ R $ will be a complete discrete valuation ring of
characteristic zero, with field of fractions $K$, and residue
$\kappa $ of characteristic $ p > 0 $ and algebraic over $ \F_p $.
(Please note that in the Appendix~\ref{sec:appendix} $K$ will
have a different meaning, whereas $R$ will not.)
\section{Some preliminary material}
\label{sec:prelim}

We begin with recalling Coleman's integration theory in the form and
to the extent
that it will be needed for this work. References for the theory are
\cite{Col82} and \cite{Col-de88}. There is also a short summary
in~\cite{Bes98b}.

Our basic data is a ``basic wide open'' in the sense of Coleman. The data
defining such an object consist of a complete curve $C/\C_p$, which is
defined over some discretely valued subfield and which has good reduction
$\Cbar$ (the reader may take $\PP^1$ for $C$ since this is the only case
that will be used in this paper), together with a finite nonempty set of points
$S\subset \Cbar(\Fbar_p)$ where $\Fbar_p$ is the algebraic closure of
$\F_p$.
To every point $y\in \Cbar(\Fbar_p)$ corresponds
a ``residue disc'' $U_y$, a subspace of the rigid analytic space associated
with $C$, consisting of all points in $C$ whose reduction is $y$. The basic
wide open $U=U_\lambda$ associated with the data above is a rigid
analytic space obtained
from $C$ by ``removing discs of radius $\lambda<1$ from the insides of the
residue discs $U_y$ for $y\in S$''. Technically this means that if the
point $y$ is locally defined by the equation $\bar{z}=0$, with $z=z_y$ some
local parameter near $y$, then one removes the points in $U_y$ where $|z|\le
\lambda$. This procedure depends on the choice of $z$ but becomes
independent of this choice as $\lambda$ approaches $1$. We will not fix
$\lambda$ but think of it as approaching $1$ and will take it as large as
needed. From now on we will use $U_y$ to denote the residue disc of $y$
\emph{in $U$}, which is the intersection of the residue disc with $U$. This
is the same as before unless $y\in S$ in which case $U_y$ is an annulus
given by the equation $\lambda < |z_y|<1$. Our final basic datum is a
Frobenius endomorphism. This is a rigid analytic map
$\phi:U_{\lambda_1}\to U_{\lambda_2}$, for some $\lambda_1$ and
$\lambda_2$, whose
reduction $\bar{\phi}$ is some power of the Frobenius endomorphism of some
model of $\Cbar$ over a finite field, extended $\Fbar_p$ linearly. A good
example of such a morphism is $\phi(z)=z^q$ on $\PP^1$ for some power $q$
of $p$.

The goal of Coleman's theory is to integrate certain differential forms on
$U$. This is first done locally, on each residue disc $U_y$. If $y\notin S$
this residue disc is isomorphic to the open disc $\{ |z|<1\}$. a rigid
differential form on such a disc has a convergent power series expansion
$\sum_{n\ge 0}z^n \dd z$ and integration is done term by
term. When $y\in S$ the form $\dd z_y/z_y$ is also analytic on $U_y$ and so
there is no choice but to introduce a logarithm
\begin{definition}\label{branchlog}
  Let $\pi\in \C_p^*$ be such that $|\pi|<1$. The branch of the $p$-adic
  logarithm determined by $\pi$ is the unique function
  $\log=\log_\pi:\C_p^\times \to \C_p$ which is multiplicative,
  defined by the usual power series when $|z-1|<1$ and satisfies
  $\log(\pi)=0$.
\end{definition}
We fix once and for all a branch of the logarithm. Then the integral of
$\dd z_y/z_y$ can be taken to be $\log(z_y)$ and that allows integration of an
arbitrary form on any $U_y$.

Let $A(U)$ (respectively $\Omega^1(U)$) be the ring of rigid analytic functions
(respectively one forms) on $U$ and let
$A_{\textup{loc}}(U)$ (resp.\ $\Omega_{\textup{loc}}^1(U)$) be the ring of $\C_p$-valued
functions (resp.\ module of one forms) on $U$ which are in
$A(U_y)$ (resp.\ $\Omega^1(U_y)$) for each $y\notin S$ and are in the
polynomial ring
$A(U_y)[\log z_y]$ (resp.\ in $A(U_y)[\log z_y]dz_y$) when $y\in S$. It is
implicit in this definition that it is independent of the choice of the
local parameter $z_y$, a fact which follows because for any two choices of
$z_y$ the difference between the $\log(z_y)$ is in $A(U_y)$.

Each $\omega\in \Omega_{\textup{loc}}^1(U)$ can be integrated in $A_{\textup{loc}}(U)$
in many ways, because we can choose a different constant of integration for
each $U_y$. Coleman's theory finds a subclass of forms for which one can
assign canonically an integral in $A_{\textup{loc}}(U)$ defined up to a global
constant. This is done recursively as follows. First one finds integrals to
all forms $\omega\in \Omega^1(U)$. At each stage one integrates all forms
that can be written as $\sum f_i \omega_i$ where $f_i$ are integrals which
have been found in previous stages and $\omega_i\in \Omega^1(U)$. The rules
for finding the integrals are:
\begin{enumerate}
 \item The integral is additive.
 \item When $g\in A(U)$, $\int \dd g=g+C$, for some constant $C$.
 \item We have $\phi^\ast \int \omega = \int \phi^\ast \omega+C$.
\end{enumerate}
The fact that these rules suffice to carry out the integration process
uniquely and
that it is independent of the choice of $\phi$ is the main result of
Coleman (see \cite{Col82}~and~\cite{Col-de88}).
One other result about Coleman integration that will be used is the
following.
\begin{proposition}
  Let $f\in A(U)^*$. Then the Coleman integral of
  the form $\dd f/f$ is $\log(f)$.
\end{proposition}
\begin{proof}
  See \cite[Lemma 2.5.1]{Col-de88}.
\end{proof}

The original reason that Coleman integrals were introduced is probably to
give a $p$--adic analogue of complex iterated integrals. Let $\omega_1$,
$\omega_2$, $\ldots$, $\omega_r$ be forms in $\Omega^1(U)$ and let $x\in U$.
Then we can define an iterated integral
\begin{equation*}
  f_r(z)=\int_x^z \omega_1\circ \omega_2\circ\cdots\circ \omega_r
\end{equation*}
by defining $f_1(z)=\int \omega_1$ normalized so that $f_1(x)=0$ and then
by induction $f_k(z)=\int f_{k-1}\omega_k$ again normalized so that
$f_k(x)=0$.

The definition of $p$--adic polylogarithms is a slight modification of the above.
Here we take $\omega_1= - \dlog(1-z)$ and $\omega_i=\dlog z$ for $i>1$. Notice that
$\dlog z$ has a simple pole at $0$. However, if we normalize $f_i$ at each
step to vanish at $0$ this zero will cancel with the pole and we will
obtain a form which is also integrable at the residue disc of $0$. This
gives the definition of the introduction.

\begin{remark}\label{galequiv}
If $U$ is a wide open defined over a field $L$, containing at least
one $L$-rational point $x$ and suppose we have chosen the branch
$\log_\pi$ such that $\pi\in L$. Then, if one has forms
$\omega_1,\ldots,\omega_r$ which are all defined over $L$ again, then an
iterated Coleman integral $f=\int \omega_1\circ \cdots \circ \omega_r$, where
the constants are fixed so that all the intermediate integrals
$\int \omega_1 \circ \cdots \circ\omega_k$ take an $L$--rational value
at $x$ are
Galois equivariant in the sense that for every automorphism $\sigma$ of
$\C_p$ over $L$ we have for every $z\in U$ that
$f(z^\sigma)=(f(z))^\sigma$. In particular, if $z$ is defined over $L$ then
$f(z)\in L$. For $\Li_n$, since the forms are either $\dlog z$
or $\dlog (1-z)$,
which are all defined over $\Q_p$, this means that if we take a branch
$\log_\pi$, with $\pi\in \Q_p$, then $\Li_n$ is Galois equivariant
over $\Q_p$.
\end{remark}

We now want to collect some facts about the functions $ \Li_n $ and other
things here that we need in the rest of the paper.

We begin with recalling some results from \cite{Col82}.  The
following is contained in Proposition~6.1 and Corollary~7.1a
of loc.~cit.  (Note that Proposition~6.1 of loc. cit. contains
an obvious misprint.)

The functions $ \Li_n(z) $ are defined on $ \C_p \setminus\{1\} $.
If $ L $ is a finitely ramified extension of $ \Q_p $ then the
limit $ \lim_{\overset{z\to1}{z\in L}} \Li_n(z) $ exists for $ n \geq 2 $,
and is independent of $ L $.
Using this limit as the value for $\Li_n$ at 1, $ \Li_n $ extends to
a function on $ \C_p $, which is continuous on finitely ramified
extensions of $ \Q_p $.

If $ m $ and $ n $ are integers at least equal to 2, then on
$ \C_p $
\begin{equation}
\label{poly-distribution}
\Li_n(z^m) = m^{n-1} \sum_{\zeta^m=1}\Li_n(\zeta z)\;.
\end{equation}
Clearly the same formula holds for $ n=1 $ provided $ 1-z^m \neq 0 $.

Let $ \log_a $ and $ \log_b $ be two different branches of the
logarithm.  Denote the corresponding different branches of $ \Li_n $
by $ \Li_{n,a} $
and $ \Li_{n,b} $.  Let $ \b = \log_a p - \log_b p $, and
let $ v $ be the valuation such that $ v(p)=1 $. Note that
\begin{equation}
  \label{eq:diflog}
  \log_a(z) - \log_b(z) = v(z) \b\;.
\end{equation}
\begin{proposition}
\label{branch-prop}
We have
\begin{multline}
\label{branch-dependence}
\Li_{n,a}(z) - \Li_{n,b}(z) =
\\
-\frac{1}{n!} v(1-z) \b
\left(
\log_a^{n-1}z + \log_a^{n-2} z \log_b z + \dots +
\log_a z \log_b^{n-2} z + \log_b^{n-1} z
\right)\;.
\end{multline}
\end{proposition}

\begin{proof}
We first remark that by the construction of Coleman integrals the
polylogarithm depends on the branch of the log chosen only on residue discs
where one of the forms involved in the definition, i.e., $\dd z/z$ and
$\dd z/(z-1)$, has a pole.
This means that $\Li_{n,a}$ and $\Li_{n,b}$ can differ at most
on the residue discs of $0$, $1$ and $\infty$, and in fact only on the
latter two discs because $ \Li_n(z) $ is analytic on $|z|<1$. We note
that apriori it would seem that because the constant of integration is
determined by the value at $0$ the function could depend on the branch
of the log everywhere, but this is not the case exactly because logs
do not appear in $\Li_n$ at the residue disc of $0$. Because
$v(1-z)\ne 0$ only on the residue discs of $1$ and $\infty$ the formula
is proved except in the cases $|z|>1$ and $|z-1|<1$.
Suppose $ |z| > 1 $. Using \eqref{functionaleq} we obtain
\[
\begin{split}
& \,\,\,\,\,\,\,\, \Li_{n,a}(z) - \Li_{n,b}(z)
\\
& =
(-1)^n ( \Li_{n,b}(1/z) - \Li_{n,a}(1/z) )
    -\frac{1}{n!} (\log_a^n z - \log_b^n z)
\\
& =
-\frac{1}{n!} (\log_a^n z - \log_b^n z)
\\
& =
-\frac{1}{n!} (\log_a z - \log_b z)
     ( \log_a^{n-1}z + \log_a^{n-2} z \log_b z + \dots + \log_b^{n-1}(z) )
\\
& =
-\frac{1}{n!} v(1-z) \b ( \log_a^{n-1}z + \log_a^{n-2} z \log_b z + \dots + \log_b^{n-1}(z) )
\end{split}
\]
because $ v(z) = v(1-z) $ for such $ z $.
It remains to consider the case $|z-1|<1$. Note that here $\log(z)$
is independent of the branch so the formula to be proved reads
\[
\Li_{n,a}(z) - \Li_{n,b}(z) =
-\frac{1}{(n-1)!} v(1-z) \b \log^{n-1} z\;.
\]
We prove this by induction on $ n $.
For $ n=1 $ this follows immediately from \eqref{eq:diflog}.
Assume $ n > 1 $.
According to \cite[Proposition~7.1]{Col82},
$ \Li_{n,a}(z) - \Li_{n-1,a}(z) \log(z) $ extends to an analytic
function on $ |1-z| < 1 $.  (Note that $ B(0,1) $ should be replaced
with $ B(1,1) $ everywhere in the formulation and the proof of
loc.\ cit.).
The result will follow from the induction hypothesis if we show
that
\[
 \gamma_n(z):=(\Li_{n,a}(z) - \frac{1}{n-1}\log(z)\Li_{n-1,a})-
 (\Li_{n,b}(z) - \frac{1}{n-1}\log(z)\Li_{n-1,b})=0\;.
\]
When we differentiate $\gamma_n(z)$
we find
\begin{align*}
&\phantom{\mathrel{=}}
\dd (\Li_{n,a}(z) - \frac{1}{n-1}\log(z)\Li_{n-1,a})
-\dd (\Li_{n,b}(z) - \frac{1}{n-1}\log(z)\Li_{n-1,b})
\\
&=
\left(\left(1-\frac{1}{n-1}\right) \Li_{n-1,a}(z)-\frac{1}{n-1}
  \log(z)\Li_{n-2,a}(z)\right) \dlog(z)
\\
&\phantom{\mathrel{=}}
-
\left(\left(1-\frac{1}{n-1}\right) \Li_{n-1,b}(z)-\frac{1}{n-1}
  \log(z)\Li_{n-2,b}(z)\right) \dlog(z)
\\
& =\frac{n-2}{n-1} \gamma_{n-1}(z)\dlog(z)= 0
\end{align*}
by induction. So the $\gamma_n(z)$ is a constant on $|z-1|<1$,
call it $C$, and we must show that $C=0$.
But $\gamma_n(z)$
satisfies the distribution relation corresponding to
\eqref{poly-distribution}. For $|z-1|<1$ this relation now reads
$C=m^{n-1}\cdot m\cdot C$, which shows $C=0$ as required.
\end{proof}

\begin{proposition}
\label{independent}
Assume the Beilinson--\soule\ conjecture holds for fields of
characteristic zero.  Let $ \log_a $ and $ \log_b $ denote two
branches of the logarithm, and denote the corresponding functions
involving $ \Li $'s by a subscript $ a $ or $ b $.  If
the maps
\[
\tM_k(K) \to \C_p
\]
given by mapping $ \symb x k $ to $ \Lmod k {}_a (x) $ (resp.
$ \Lmod k {}_b (x) $)
are well defined for $ k\leq n $, then the induced maps
\[
H^1( \tildescriptM n (K) ) \to \C_p
\]
given by mapping $ \symb x n $ to $ \Lmod n {}_a (x) $ (resp.
$ \Lmod n {}_b (x) $)
are the same.
\end{proposition}

\begin{remark}
$\tM_k(K)$ will be constructed below in Section~\ref{sec:kstuff}, but
for a heuristic approach to working with it see the beginning of
the Introduction.
Our computation of the regulator map in Sections to come will
show that, for a fixed choice of $ \log $, the map
\[
\tM_n (R) \to \C_p
\]
given by mapping $ \symb x n $ to $ \Lmod n (x) $
is well defined, but we have to assume this for $ \tM_n (K) $.
Note also that for the special units, the function $ \Lmod n (x) $
does not depend on the branch of the logarithm by
Proposition~\ref{independent}.
\end{remark}

\begin{proof}[Proof of Proposition~\textup{\ref{independent}}]
Because $ \tM_n(K) $ is by definition~\cite[Corollary~3.22]{Jeu95} the
quotient of $ M_n(K) $
by the subspace spanned by elements of the form $ \symb x n + (-1)^n
\symb 1/x n $,
and $ \Lmod n $ vanishes on those elements by \eqref{zooz}, we can just as well
prove the corresponding result for $ H^1(\scriptM n (K)) $ as
the assumption about the Beilinson--\soule\ conjecture implies
that the map $ H^1(\scriptM n (K) ) \to H^1(\tildescriptM n (K) ) $
is a quasi--isomorphism.
Note that all functions of the form $ \symb x n \mapsto \log^{n-k} (x) \Lmod k (x) $
are well defined on $ M_n(K) $ by assumption as they factor through
the map
\[
M_n(K) \to  M_{n-1}(K) \tensor \KQ \to\dots \to M_k(K) \tensor
 (\KQ)^{\tensor n-k}
.
\]
The function $ \Li_n(z) $ can be expressed as a linear combination of
the functions $\Lmod k (z) \log^{n-k}(z) $ with
$  1 \leq k \leq n $, by the definition of $ \Lmod n $,
and on elements in $ H^1(\scriptM n (K)) $ those functions vanish
for $ k<n $.  Therefore the map from $ M_n(K) $ to $ \C_p $ given
by mapping $ \symb x n $ to $ \Li_n(x) $ is well defined under
our assumptions, and it suffices to show that it induces a map
on $ H^1(\scriptM n (K)) $ which is independent of the branch
of the logarithm.
Proposition~\ref{branch-prop} shows that this is true, as
the right hand side of equation~\ref{branch-dependence}
is zero on elements in $ H^1(\scriptM n (K)) $, as this map factors
through the map
$
M_n(K) \to M_{n-1}(K) \tensor \KQ \to\dots\to (\KQ)^{\tensor n}
$
above, which maps $ \symb x n $ to $  (1-x)\tensor x\tensor\dots\tensor x $.
\end{proof}

We shall also need the distribution relation for elements in
$ M_n(F) $, as given in~\cite[Proposition~6.1]{Jeu95}.

\begin{proposition}
\label{distribution}
If $ F $ is a field of characteristic zero that contains the $ m $--th roots of unity,
then in $ M_n(F) $ (and hence $ \tM_n(F) $) we have
\begin{equation}
\symb x^m n = m^{n-1} \sum_{\alpha^m = 1} \symb {\alpha x} n\;.
\end{equation}
\end{proposition}
\section{Some $ K $--theory}
\label{sec:kstuff}

In this Section we construct the complexes $ \tildescriptM n (\O) $
as quotient complexes of complexes $ \scriptM n (\O) $ for $ n \geq 2 $.
The main idea is the same as in \cite{Jeu95}, but the fact that
we will be working with a discrete valuation ring rather
than a field gives rise to complications.  In order to highlight
the idea we start with a rather gentle exposition.  For the proofs of
the statements that are used in the construction, we refer the
reader to loc.\ cit., especially Sections~2.1 through~2.3, and~3.
In loc.\ cit.\ most of the work was done over $ \Q $,
but in fact the proofs hold over our base $ \O $ , a discrete
valuation ring of characteristic zero, 
without any change.  There is also a very brief introduction to multi
relative $ K $--theory in Appendix~\ref{sec:appendix}.

The idea of the whole construction is the following.
If $ B $ is a regular Noetherian scheme, then the pullback
$ K_*(B) \to K_*(\A_B^1) $ 
is an isomorphism.  We shall be using
an Adams decomposition with respect to weights,
$ K_m(X)_\Q = \oplus_i \K m i X $.  The weight behaves naturally
with respect to pullback, and under suitable hypotheses for a
closed embedding, there
is a pushforward Gysin map with a shift in weights corresponding
to the codimension (see, e.g., \cite[Proposition~2.3]{Jeu95}).

Let $ X_B = \PP_B^1\setminus\{t=1\} $ with $ t $ the standard
affine coordinate on $ \PP^1 $.  Write $ \bbox_B^1 $ for the
closed subset $ \{t=0,\infty\} $ in $ \PP_B^1 $.  Then the relative
exact sequence for the couple $ (X_B; \bbox_B^1) $
gives us
$$
\dots\to
K_{n+1}(X_B) \to K_{n+1}(\bbox_B^1) \to K_{n}(X_B;\bbox_B^1) \to K_{n}(X_B) \to K_{n}(\bbox_B^1)
\to\cdots
$$
for $ n \geq 0 $.  Because the map pullback $ K_{n+1}(B) \to K_{n+1}(X_B) $
is an isomorphism, combining it with the pullback 
$  K_{n+1}(X_B) \to K_{n+1}(\bbox_B^1) = K_{n+1}(B)^2 $
shows that the map $ K_{n+1}(X_B) \to K_{n+1}(\bbox_B^1) $ corresponds
to the diagonal embedding $ K_{n+1}(B) \to K_{n+1}(B)^2 $.  As this holds
for all $ n \geq 0 $, we get that we have an isomorphism
$ K_{n}(X_B;\bbox_B^1) \iso K_{n+1}(B) $ for $ n \geq 0 $.  Note
that we have a choice of sign here in the isomorphism of the
cokernel of $ K_n(B) \to K_n(B)^2 $ with $ K_n(B) $.  This results
in similar choices of signs in the maps $ H^p(\scriptM n (\O)) \to \K 2n-p n {\O} $
(resp. $ H^p(\tildescriptM n (\O)) \to \K 2n-p n {\O} $)
later on in this Section.

We  can iterate this procedure using multi relative $ K $--theory.  
(The construction of this is recalled in Appendix~\ref{sec:appendix}.)
For the sake of exposition we give the argument here for the next level of relativity.
If we let $ \bbox_B^2  = \{t_1=0,\infty\};\{t_2=0,\infty\}$,
then we can get a long exact sequence
\begin{align*}
& \cdots \to
K_{n+1}(X_B^2;\{t_1=0,\infty\}) \to  K_{n+1}(\{t_2=0,\infty\};   \{t_1=0,\infty\} )
\to
\cr
& 
\to K_{n}(X_B^2;\bbox_B^2) 
\to
K_{n}(X_B^2;\{t_1=0,\infty\})
\to K_{n}(\{t_2=0,\infty\} ;  \{t_1=0,\infty\} )
\to\cdots
.
\end{align*}
Using induction on the degree of relativity one sees that the composition
\begin{align*}
&
K_{n+1}(X_B;\{t_1=0,\infty\}) \to K_{n+1}(X_B^2;\{t_1=0,\infty\})
\to 
\cr
& \qquad\qquad \to
K_{n+1}(\{t_2=0,\infty\};   \{t_1=0,\infty\} ) \iso K_{n+1}(X_B;\{t_1=0,\infty\})^2
\end{align*}
(with the first the pullback along the projection $ (t_1,t_2) \mapsto t_2 $)
is the diagonal embedding, hence we obtain an isomorphism
$ K_{n}(X_B^2;\bbox_B^2) \iso K_{n+1}(X_B;\bbox_B^1) $
for $ n \geq 0 $.
Therefore we get
$ K_n(X_B^2;\bbox_B^2) \iso K_{n+1}(X_B;\bbox_B)  \iso K_{n+2}(B)  $ 
for $ n \geq 0 $.  By induction one proves that
\begin{equation}
  \label{eq:downink}
  K_n(X_B^m;\bbox_B^m)  \iso K_{n+m}(B)
\end{equation}
for $ n \geq 0 $ and $ m \geq 1 $, with $ \bbox_B^m $ shorthand for 
$ \{ t_1=0,\infty\};\dots;\{ t_m=0,\infty\} $
($ m $--th order relativity).
One can also do this with weights, and as the weight are compatible
with pullbacks, we get isomorphisms
$ \K n j X_B^m;\bbox_B^m  \iso \K n+m j B  $ for $ n \geq 0 $
and $ m \geq 1 $.

In order to get elements in groups like $ K_{n+m}(X_B^m;\bbox_B^m) $,
we use localization sequences.  We shall explain the idea for
$ m=1 $.  
(For $ m\geq2 $ the localization sequences get replaced
by a spectral sequence, see below.)
If $ u $ is an element in our discrete valuation ring
$ \O $ such that both $ u $ and $ 1-u $
are units, then we get an exact localization sequence
$$
\dots\to
K_{m}(\O) \to
K_m(X_\O;\bbox_{\O}^1)
\to K_m(X_{\O,\loc} ;\bbox_\O^1 ) \to K_{m-1}(\O) 
\to\cdots
$$
where $ X_{\O,\loc}
 = X_\O \setminus \{ t=u \} $ and we identified
$ \{ t=u \} \subset X_\O $ with $ \O $ (or rather $ \Spec(\O) $).
We used here that $ u $ and $ 1-u $ are units in $ \O $ so that
$ \{ t = u \} $ does not meet $ \bbox_\O^1 $ or $ \{ t = 1 \} $,
and that $ \O $ is regular in order to identify $ K_m(\O) $ with
$ K_m'(\O) $.
(If we want to leave out $ \{ t = u \} $ and $ \{ t = v \} $ simultaneously
for two distinct elements $ u $ and $ v $ in $ \O $ such that all of
$ u $, $ v $, $ 1-u $ and $ 1-v $ are units, which we shall do
below, this already becomes far more complicated and one is force
to use a spectral sequence.)
The image of $ K_2(\O) \to K_2(X_{\O};\bbox_{\O}^1) $ can be controlled by
looking at the weights, which for the bit that we are interested in
gives us
$$
\dots\to
\K 2 1 {\O} \to
\K 2 2 {X_\O;\bbox_\O^1}
 \to \K 2 2 {X_{\O,\loc};\bbox_\O^1}  \to \K 1 1 {\O} 
\to\cdots
.
$$
Because of weights in $ K $--theory, one knows that $ \K 2 1 {\O} = 0 $,
so we can analyze $ \K 2 2 {X_\O;\bbox_\O^1} $
 as subgroup of
$ \K 2 2 {X_{\O,\loc};\bbox_\O^1} $.  In \cite[Section~3.2]{Jeu95}
universal elements $ [S]_n $ were constructed, of which we want
to use $ [S]_2 $ here.  It gives rise to an element $ [u]_2 $
in $  \K 2 2 {X_{\O,\loc};\bbox_\O^1} $ with boundary $ (1-u)^{-1} $
in $ \K 1 1 {\O} $.  If we use this for various $ u $ (suitably
modifying the localization sequence above into a spectral sequence)
and consider elements coming from the cup product 
$$
\K 1 1 {X_{\O,\loc};\bbox_\O^1} \times \K 1 1 {\O} \to \K 2 2 {X_{\O,\loc};\bbox_\O^1} 
$$
we can get part of $ \K 2 2 {X_{\O};\bbox_\O^1} \iso \K 3 2 {\O} $
by intersecting the kernel of the map corresponding to
$ \K 2 2 {X_{\O,\loc};\bbox_\O^1} \to \K 1 1 {\O} $
with the space generated by the symbols $ [u]_2 $ and the image
$  \K 1 1 {X_{\O,\loc};\bbox_\O^1} \cup \K 1 1 {\O}  $ of the cup product.

Unfortunately, this gets fairly technical, but after this gentle
introduction we are now ready to begin.  The reader is encouraged
to compare this construction with the simpler construction for
fields, which is carried out in \cite[Section~3]{Jeu95}.

To ease the notation somewhat, we will drop the subscript (indicating
the base scheme) from $ \bbox^n $.

\begin{definition}\label{no-low-weight-def}
A scheme $ B $ has no low weight $ K $--theory if the Beilinson--\soule\ 
conjecture holds for $ B $, i.e., $ \K m j B = 0 $ if $ 2 j \leq m $
and $ m > 0 $. A ring $A$ is said to have no low weight $ K $--theory
if $\Spec(A)$ does.
\end{definition}

We shall use the following notation.  Let $ t $ be the standard
affine coordinate on $ \PP^1 $.  We let $ X= \PP_\Z^1 \setminus \{ t=1 \} $.
If $ B $ is any scheme, we let $ X_B = X \times_\Z B $, and
$ X_B^n = X_B\times_B\dots\times_B X_B $.
If $ U $ is a subset of $ \Gamma(B, \O^*) $ such that if $ b $
is in $ U $, then  $ 1-b $ is also in $ \Gamma(B, \O^*) $, we let
$ X_{B,\loc} = X_B\setminus\{t = b , b \in U \} $, and
$ X_{B,\loc}^n = X_{B,\loc}\times_B\dots\times_B X_{B,\loc} $.
The set $ U $ will normally be clear from the context.  We shall
also abuse notation by writing $ X_{B,\loc}^n $ even after we
took direct limits over finite sets $ U $.
In the multi--relative $ K $--theory below, we shall write
$ \bbox^n $ for $ \{ t_1=0,\infty \};\dots;\{ t_n=0,\infty \} $.

\begin{notation}
For the remainder of the section, $ \O $ will be a discrete valuation
ring with field of fractions $ F $  and residue field $ \kappa $.
\end{notation}

\begin{lemma}\label{no-low-weights}
Assume $ F $ and $ \kappa $ have no low weight $ K $--theory.
Then for $ 2j \leq q+m $ and $ m>q $, all of 
$ \Krloc m j F q $, $ \Krloc m j {\O} q $ and $ \Krloc m j {\kappa} q $
are zero.
\end{lemma}

\begin{proof}
Lemma~3.4 of \cite{Jeu95} shows the statement to be true for $ F $ or
$ \kappa $.
The result for $ {\O} $ follows immediately from the exact localization sequence
\[
\dots \to \Krloc m j-1 {\kappa} q \to \Krloc m j {\O} q \to \Krloc m j F q \to \cdots
\]
\end{proof}

\begin{remark}
\label{tr-deg-one}
$ \kappa $ has no low weights $ K $--theory if $ \kappa $ is algebraic
over $ \F_p $, because all $ K_n(\kappa) $ are torsion for $ n\geq1 $.
It also holds if $ \kappa $  is of transcendence degree 1 over $ \F_p $ by
a result of Harder, see \cite[Korollar~2.3.2]{Har77}.  Because
$K_n(F)$ is torsion for $n \geq 1$ for a finite field $F$, using
localization it is enough to show that $K_n(\O)$ is torsion for a
Dedekind ring in a function field of transcendence degree 1 over a
finite field, which is the result quoted.
\end{remark}

\begin{remark}
$ F $ has no low weights $ K $--theory if $ F $ is a number field,
or more generally a subfield of the algebraic closure of $ \Q $.
As the residue field is an algebraic extension
of $ \F_p $ in this case, the conditions of Lemma~\ref{no-low-weights}
are certainly satisfied, and all constructions in this Section
go through without assumptions about the weights on the $ K $--groups
involved.
\end{remark}

Consider the divisors on $ X_{\O}^n $ defined by putting $ t_i =u_j $
for some $ u_j $ in $ \Osp $.
Put $  W^0 = X_{\O}^n $, and let $ W^1 $ be the union of divisors
$ \{t_i=u_j\} $ for all $ u_j $ in some finite set $ U \subset \Osp $.
Considering the singular locus of $ W^1 $, it is easy to see
that one can extend this to a stratification
$ W^0 \supset W^1\supset \dots\supset W^{n+2}=\emptyset $
on $ X_{\O}^n $, with all $ W^s\setminus W^{s+1} $  for $ s=0,\dots,n+1 $
consisting of a finite union of 
$ \xl {\O} n-s $, $ \xl F n-s $ and $ \xl {\kappa} n-s+1 $'s,
which are regular.  Using the localization sequences
\begin{align*}
&
\dots\to
\K n,W^{s+1} j+s X^n_{\O};\bbox^n \to
\K n,W^{s} j+s X^n_{\O};\bbox^n \to
\\
& \qquad\qquad\qquad\qquad
\K n,{W^s\setminus W^{s+1}} j+s {X^n_{\O}\setminus W^{s+1}};\bbox^n \to
\K n-1,W^{s+1} j+s X^n_{\O};\bbox^n \to
\cdots
\end{align*}
where $ \K n,W^{s+1} j+s X^n_{\O};\bbox^n  $ etc. is K-theory with supports.
we get an exact couple, which gives rise to a spectral sequence
converging to $ \K n j X_{\O}^n;\bbox^n $.
We have isomorphisms
$$
\K n,{W^s\setminus W^{s+1}} j+s {X^n_{\O}\setminus W^{s+1}};\bbox^n \iso
\K n j {W^s\setminus W^{s+1}};\bbox^n
$$
and we can identify the terms in the spectral sequence with terms
of this type.
Note that the components of $( W^s \setminus W^{s+1};\bbox^n )$ are
of the form
$( \xlocb {\O} n-s )$, $( \xlocb F n-s )$ and $( \xlocb {\kappa} n-s+1 )$.
Taking $ j=n+1 $ we get a spectral sequence with $ E_1^{s,t} $ equal to
\[
\Krloc -s-t n-s {\O} n-s-1 \coprod \Krloc -s-t n-s F n-s-1 \coprod \Krloc -s-t n-s {\kappa} n-s
\]
and converging to $ \K -s-t n X_{\O}^{n-1};\bbox^{n-1} $.
If we write $ \Ksnot n j {\O} m $ for $ \Krloc n j {\O} m $ for
typographical reasons,
and similarly for $ F $ and $ \kappa $,
this looks as
\begin{equation}
\label{messy-ss}
\end{equation}
\begin{alignat*}{4}
&\quad\,\,\,\vdots & &\,\,\,\vdots & &\qquad\quad\vdots & &\notag
\\
& \Ks n-1 n   {\O} n-1 \qquad
& \Ks n-2 n-1 {\O} n-2 & \scp \Ks n-2 n-1 F n-2 \qquad
& \Ks n-3 n-2 {\O} n-3 & \scp \Ks n-3 n-2 F n-3 \scp \Ks n-3 n-2 {\kappa} n-2 \qquad
&& \cdots \notag
\\
& \Ks n   n   {\O} n-1
& \Ks n-1 n-1 {\O} n-2 & \scp \Ks n-1 n-1 F n-2
& \Ks n-2 n-2 {\O} n-3 & \scp \Ks n-2 n-2 F n-3 \scp \Ks n-2 n-2 {\kappa} n-2
&& \cdots \notag
\\
& \Ks n+1 n   {\O} n-1
& \Ks n   n-1 {\O} n-2 & \scp \Ks n   n-1 F n-2
& \Ks n-1 n-2 {\O} n-3 & \scp \Ks n-1 n-2 F n-3 \scp \Ks n-1 n-2 {\kappa} n-2
&& \cdots \notag
\\
&\quad\,\,\,\vdots & &\,\,\,\vdots & &\qquad\quad\vdots & & \notag
\\
&\quad & &\quad & & & & 
\end{alignat*}
Observe that, due to the choice of the stratification,
$ \Krloc * n-s {\kappa} n-s $ occurs only when $ s\geq 2 $.  Also,
by Lemma~\ref{no-low-weights}, if both $ F $ and $ \kappa $ have
no low weights $ K $--theory, then in
the spectral sequence \eqref{messy-ss} converging to $ \Krloc * n X_{\O} n-1 $
there are no nonzero terms in the
row below the one beginning with $ \Krloc n n {\O} n-1 $
(i.e., the middle row of \eqref{messy-ss}, 
where $ \Krloc n n {\O} n-1 $ is denoted $ \Ksnot n n   {\O} n-1 $).

\begin{lemma}
\label{injects}
If $ \kappa $ has no low weight $K$-theory,
then the map
\[
\Krloc n n {\O} n-1 \to \Krloc n-1 n-1 F n-1
\]
is injective.
\end{lemma}

\begin{proof}
Immediate from the localization sequence
\[
\dots \to \Krloc n n-1 {\kappa} n-1 \to \Krloc n n {\O} n-1 \to \Krloc n n F n-1 \to \dots
\]
as the first term here is zero by Lemma~\ref{no-low-weights}.
\end{proof}

We now notice that all our localizations are compatible with
localizing in a larger set $ U $, and that we can take direct limits of our
localizations over finite sets $ U $ if we want.  In order not to overburden
the notation we shall suppress $ U $ from the notation.
Notice that this means also that all components in the spectral
sequence \eqref{messy-ss} of codimension at least one with $ \O $ as base
become the corresponding components with $ F $ as base.

\begin{assumption}
We assume for the remainder of the section that the discrete
valuation ring $ \O $ has characteristic zero.
\end{assumption}

We now define symbols in $K$--theory. Let
$G = \Spec(\Z[S,S^{-1},(1-S)^{-1}]) $. In~\cite[Section 3.2]{Jeu95}
universal symbols 
\begin{equation}
  \label{eq:Ssubn}
   \symb S n  \in  \Krloc n n G n-1
\end{equation}
were constructed, with boundary
\[
\sum_{i=1}^{n-1} (-1)^i \symb S n-1|t_i=S
\]
in $ \coprod_i  \Krloc n-1 n-1 G n-2 $ under the boundary in the spectral sequence
corresponding to \eqref{messy-ss} for $ G $.
(Although the proofs in loc. cit. were formulated over
$ \Q $, the constructions hold for a much larger class of base
schemes without any change.)
Recall that we denote by $ \O^* $ the units in $ \O $, and by $ \Osp $
the set
of elements $ u $ in $ \O^* $ such that $ 1-u $ is also in $ \O^* $.
\begin{definition}\label{symbols}
  For $u\in \Osp$ we define the symbol
  \begin{equation*}
     \symb u n  \in \Krloc n n {\O} n-1 
  \end{equation*}
  as the pullback of the universal symbol $ \symb S n $
  under the map $ \Spec(\O) \to G $
  induced by mapping $ S $ to $ u $.
\end{definition}
It was also shown in loc. cit. that the symbol $ \symb 1 n $
exists for $ n\ge2 $, but we shall tacitly ignore this symbol
here, as it can also be defined by the distribution relation
Proposition~\ref{distribution} if there are other roots of unity in $ F $.

We define inductively the symbolic part of the $ K $--theory.
Let
\begin{equation}\label{oneplusistar}
  \I = \Krloc 1 1 {\O} 1 =
  \{
  \prod_j \left(\frac{t-u_j}{t-1}\right)^{n_j} \text{ such that }
  \prod_j u_j^{n_j} =1
  \}\;,
\end{equation}
where the $ u_j $ are in $ \O^* $ and the $ n_j $ are in $ \Z $,
and let $ \Symb k {\O} \subseteq \Krloc k k {\O} k-1 $ be defined
by setting $ \Symb 1 {\O} = \O^* $, and
\[
\Symb k+1 {\O} = \Span {\symb u k+1 , u \in \Osp } + \I\tcup \,\, \Symb k {\O}
\]
for $ k\geq 1 $.
The notation $ \tcup $ means the following.  There are $ k $ projections
of $ X_{\O,\loc}^k $ to $ X_{\O,\loc}^{k-1} $, giving rise to $ k $ cup
products 
\[
\I \times \Krloc k k {\O} k-1 \to \Krloc k+1 k+1 {\O} k 
.
\]
$ \I\tcup \,\, \Symb k {\O} $ is the $ \Q $--subspace spanned by the
image of all those $ k $ cup products.

Because $ \dd \symb x n = \sum_{i=1}^{n-1} (-1)^i \symb x n-1|t_i=x $
where $\dd$ is the differential in the spectral sequence,
we get a complex  $ \cSymb n ({\O}) $
in the row of \eqref{messy-ss} which starts with
$ \Ksnot n n {\O} n-1 
= \Krloc n n {\O} n-1 $, given by
\[
\Symb n {\O} \to \!\coprod \Symb n-1 {\O} \to  \!\coprod \Symb n-2 {\O}
\to \dots \to
\!\coprod \Symb 2 {\O} \to \!\coprod \O^*
\!.
\]
If $\kappa$ has no low weight $K$--theory, then by
Lemma~\ref{injects} we can view those groups as
subgroups of the corresponding spaces for $ F $.  As the components
corresponding to $ \kappa $ in (\ref{messy-ss}) will never
play a role in the boundary for elements in $ \cSymb n (\O) $,
we can view the above complex as a subcomplex of the complex
\[
\Krloc n n F n-1 \to \coprod \Krloc n-1 n-1 F n-2 \to
\coprod \Krloc n-2 n-2 F n-3 \to \cdots
\]
where all coproducts for codimension $ r $ are taken over $ r $--tuples
$ (u_1,\dots,u_r) $ in $ (\Osp)^r $.

\begin{lemma}
\label{localize-more-F}
Assume $ F $ has no low weight $K$-theory.  Then
the map corresponding to further localization
from $ \K n n {\xlocb F n-1 } $ for one set of localizing elements
to $ \K n n {\xlocb F n-1 } $ for a larger one is injective.
\end{lemma}

\begin{proof}
For $ n=1 $ there is nothing to prove.  For $ n\geq2 $, we use
the exact sequence
\[
0 \to \K n n {X_F^{n-1} ; \bbox^{n-1}} \to
\K n n {\xlocb F n-1 } \to \coprod \K n-1 n-1 {\xlocb F n-2 }
\]
for two different set of localizing elements.  As clearly
the right most term injects under localizing more (i.e., make
the coproduct larger as well), we are done by induction.
\end{proof}

By Lemmas~\ref{injects} and~\ref{localize-more-F},
if both $\kappa$ and  $F$ have no low weight
$K$--theory, then we also have an inclusion
$ \Symb n {\O} \subset \Symb n F  $, so that $ \cSymb n ({\O}) $
is a subcomplex of $ \cSymb n (F) $.  We can also forget about
exactly which finite subset  $ U $ of $ \Osp $ or $ F^* \setminus \{1\} $ we use,
and work in the direct limit for such $ U $ from now on.

If both $ F $ and $ \kappa $ have no low weight $K$-theory, then all
this takes place in the lowest non--zero row
of the spectral sequence (\ref{messy-ss}) above, and if we give $ \cSymb n $
a cohomological grading in degrees 1 through $ n $, we get a
commutative diagram of maps
\[
\xymatrix{
H^r(\cSymb n ({\O}))\ar[d]\ar[r]& \Krel n-r+1 n {\O} n-1 \ar[d]\ar[r]^-\iso & \K 2n-r n {\O} \ar[d]
\\
H^r( \cSymb n (F) ) \ar[r]      & \Krel n-r+1 n K n-1 \ar[r]^-\iso       & \K 2n-r n F 
}
\]
because the differentials on $ \Symb n {\O} $ and $ \Symb n F $
are induced from the spectral sequence~\eqref{messy-ss}.

\begin{remark}
\label{injective-by-construction}
Note that if both $\kappa$ and $F$ have no low weight
$K$--theory, then the horizontal maps here are injections by construction
for $ r = 1 $.
\end{remark}

The complex $ \cSymb n $ can be changed into a
tensor complex, using the following subcomplex.
Define $ J_k $ as $ \I \tcup \,\, \Symb k {\O} $.
Let $ \J n $ be the subcomplex of $ \cSymb n $
given by
$$
J_n \to \dd J_n + \coprod J_{n-1} \to \dd(\dots) + \coprod J_{n-2}
\to\dots\to
\dd(\dots) + \coprod J_2 \to \dd(\dots)
.
$$

\begin{prop}\label{J-acyclic}
The subcomplex $ \J n $ is acyclic.
\end{prop}

\begin{proof}
The same as the proof of Lemma~3.7 of \cite{Jeu95}, see also Remark~3.10
in loc. cit.
\end{proof}

Note that the symmetric group $ S_{n-1} $ acts on $ \cSymb n $ and $ \J n $.
Denote the parts of those complexes on which $ S_{n-1} $ acts
alternatingly by the superscript $ \alt $.  Let the complex $ \scriptM n $
be the quotient complex $ \left( \cSymb n \right)^\alt / \left(\J n \right)^\alt $.
It has the form
$$
M_n \to M_{n-1} \tensor \OQ \to M_{n-2} \tensor \W^2 \OQ \to
\dots \to
M_2 \tensor \W^{n-2} \OQ \to \OQ \tensor \W^{n-1} \OQ
$$
with $ M_k({\O}) = \Symb k {\O} ^\alt / J_{(k)}^\alt  $, which is generated
by the classes of the elements $ \symb u k $, where $ u \in \Osp $,
and similar for $ F $.
Denote the class of $ \symb u k $ simply by  $ \symb u k $.
Then the differential is given by
\[
\dd (\symb x k \tensor y_1\w\dots\w y_{n-k}) =
     \symb x k-1 \tensor x\w y_1\w\dots\w y_{n-k}
\]
if $ k \geq 3 $, and 
\[
\dd (\symb x 2 \tensor y_1\w\dots\w y_{n-k}) =
        (1-x)\tensor x\w y_1\w\dots\w y_{n-k}
.
\]

\begin{proposition}
\label{O-F-injective}
If $ \kappa $ and $ F $ have no low weight $ K $--theory,
then the localization map $ \scriptM n ({\O}) \to \scriptM n (F) $ is injective.
\end{proposition}

\begin{proof}
 From the localization sequence
\[
\dots\to
\K 2n-1 n-1 {\kappa} \to
\K 2n-1 n {\O} \to
\K 2n-1 n F \to
\cdots
\]
we get that the map $ \K 2n-1 n {\O} \to \K 2n-1 n F $ is injective.
Consider the commutative diagram
\[
\xymatrix{
0 \ar[r] & H^1(\scriptM k (\O)) \ar[r]\ar[d] & M_{(k)}(\O)\ar[r]\ar[d]
    & M_{(k-1)}(\O)\tensor\OQ\ar[d]
\\
0 \ar[r] & H^1(\scriptM k (F)) \ar[r] & M_{(k)}(F)\ar[r]
    & M_{(k-1)}(F)\tensor\FQ .
}
\]
Because the maps to $ \K 2n-1 n {\O} $ (resp. $ \K 2n-1 n F $)
the $ H^1 $'s are injective by Remark~\ref{injective-by-construction},
we find that the map
$ H^1 (\scriptM k ({\O}) ) \to H^1( \scriptM k (F) ) $
is injective.  It follows by induction
on $ k $ that the maps $ M_k({\O}) \to M_k(F) $ are injective, as
this is clear for $ k = 2 $, where  the vertical map on the right
is the inclusion $ \OQ \tensor \OQ \to \FQ \tensor \FQ $.
\end{proof}

Let $ N_k = N_k(\O)= \Span {\symb u k + (-1)^k \symb u^{-1} k } $
where the $u$ runs through $\Osp$.
This is the analogue of the groups $N_k(F)$
defined in~\cite[Proposition~3.20]{Jeu95}.
We consider the subcomplex $ \N n (\O)$ of $ \scriptM n (\O)$ given by
\[
N_n \to N_{n-1} \tensor \OQ \to\dots\to
N_2 \tensor \W^{n-2} \OQ \to
\dd \left( N_2 \tensor \W^{n-2} \OQ \right)
.
\]

\begin{lemma}
\label{generate}
If $ \kappa $ has more than two elements, or equivalently, if
$ \Osp \neq \emptyset $, then $ \dd N_2({\O}) = \Sym^2 \OQ $, and
similarly for $ F $.
\end{lemma}

\begin{proof}
Because $ \dd \symb u 2 + \dd \symb u^{-1} 2 = u \tensor u $
in $ \Sym^2 \OQ $, $ \dd N_2 = \Span {u \tensor u, u \in \Osp} $.
$ \Sym^2 \OQ $ is spanned by elements of the form $ v \tensor v $
where $ v $ runs through $ \O^* $.  If $ v $ is a special unit,
then it is clear from the formula above that $ v \tensor v $
is in $ \dd N_2 $.  If not, $ v $ reduces to $ 1 $ in $ \kappa $.
Let $ w $ be a non--special unit in $ \O^* $.  then $ wv $, $ w v^{-1} $
and $ w v^2 $ are all non--special, and they give the elements
$ (wv) \tensor (wv) $, $ (w v^{-1}) \tensor (w v^{-1}) $ and
$ (wv^2) \tensor (wv^2) $ in $ \Sym^2 (\OQ) $, and a linear combination
of them gives $ v \tensor v $.
\end{proof}

\begin{prop}\label{n-acyclic}
If the Beilinson--\soule\ conjecture holds for $ \kappa $ and
for fields of characteristic
zero, then the complex $ \N n (\O) $ is acyclic.
\end{prop}

\begin{proof}
If there are no special units, there is nothing to prove as the
complex $ \N n (\O) $ is zero.
If $ \Osp $ is non--empty, we show that $ \N n (\O) $
is quasi--isomorphic to the complex
\[
\Sym^n (\OQ) \to
\Sym^{n-1}\tensor_\Q \OQ \to
\Sym^{n-2}(\OQ) \tensor_\Q \W^2 \OQ \to
\cdots\qquad\qquad\qquad
\]
\[
\qquad\qquad\qquad\dots\to
\Sym^2(\OQ) \tensor_\Q \W^{n-2} \OQ \to
\dd \left( \Sym^2(\OQ) \tensor_\Q \W^{n-2} \OQ \right)
.
\]

It is well known that this last complex is acyclic, with an explicit
homotopy operator given in Corollary~3.22 of \cite{Jeu95}.
It was proved in Proposition~3.20 of \cite{Jeu95} that the map
$ \symb u n + (-1)^n\symb 1/u n \mapsto u \tensor\dots\tensor u $
induces an isomorphism between $ N_n(F) $ and the subspace
of $ \Sym^n(\FQ) $ generated by the elements $ u\tensor\dots\tensor u $
with $ u $ in $ F^* $.  Considering Proposition~\ref{O-F-injective}, our complex
is a subcomplex of the corresponding complex for $ F $.  So it
will suffice to check that the image of $ N_n({\O}) $ is
$ \Sym^n(\OQ) $, which is done as in the proof of Lemma~\ref{generate}
\end{proof}

\begin{remark}\label{n-acyclic-number-field}
In case $ F $ is a number field, we can prove the statement
of Proposition~\ref{n-acyclic} without assuming the Beilinson--\soule\
conjecture.  It is known that $ F $ satisfies the Beilinson--\soule\
conjecture,
and it was shown in \cite[Proposition~5.1]{Jeu95} that the map 
$ \symb u n + (-1)^n\symb 1/u n \mapsto u \tensor\dots\tensor u $
gives an injection from $  N_n(F) $ into $ \Sym^n {\FQ} $.  Because
$ \kappa $ is finite, $ K_m(\kappa) $ is torsion for $  m \geq 1 $,
we have inclusions $ N_n({\O}) \subset N_n(F) $ by Proposition~\ref{O-F-injective}.
The proof that  $ \N n (\O) $ is acyclic then proceeds as in the general case.
\end{remark}

Let $ \tildescriptM n (\O) $ be the quotient complex
$ \scriptM n (\O)/ \N n (\O) $.
It has the form
$$
\tM_n (\O) \to \tM_{n-1} (\O) \tensor \OQ \to \tM_{n-2} (\OQ) \tensor \W^2 \OQ \to
\dots \to
\tM_2 (\O) \tensor \W^{n-2} \OQ \to \W^{n} \OQ
$$
with $ \tM_k (\O) = M_k(\O)/N_k(\O)  $, and is clearly still generated by
the classes of the elements $ \symb u k $, where $ u \in \Osp $.
We have similarly the complex $ \tildescriptM n (F) $ defined
in~\cite[Corollary~3.22]{Jeu95}.  In both
cases, the differential is now given by
\[
\dd (\symb x k \tensor y_1\w\dots\w y_{n-k}) =
     \symb x k-1 \tensor x\w y_1\w\dots\w y_{n-k}
\]
if $ k \geq 3 $, and 
\[
\dd (\symb x 2 \tensor y_1\w\dots\w y_{n-k}) =
        (1-x)\w x\w y_1\w\dots\w y_{n-k}
.
\]

We have now proved the assertions (1) in Theorem~\ref{general-field}
and Theorem~\ref{number-field}.  Namely, either assume the
Beilinson--\soule\ conjecture is true for fields
of characteristic zero as well as for $ \kappa $, or that 
$ F $ is a number field.  We then have a map 
\[
H^r(\tildescriptM n ({\O}) ) \to \K 2n-r n {\O}
,
\]
which is obtained as the composition of the maps
\begin{equation}
  \label{eq:kbasicmap}
  H^r( \tildescriptM n ({\O}) ) \isoot H^r( \scriptM n ({\O}) ) \isoot
  H^r( \cSymb n ({\O}) ) \to \K 2n-r n {\O} \;.
\end{equation}
From left to right, those maps are justified by
Proposition~\ref{n-acyclic} or Remark~\ref{n-acyclic-number-field},
Proposition~\ref{J-acyclic}, and Lemma~\ref{no-low-weights}
as it implies
that we are working in the lowest nonvanishing row of the spectral
sequence~\eqref{messy-ss}.

One has the same maps when replacing $ {\O} $ with $ F $ everywhere,
and the inclusion of $ {\O} $ into $ F $ induces an injection
of complexes $ \scriptM n ({\O}) $ into $ \scriptM n (F) $
by Proposition~\ref{O-F-injective},
which is compatible with all those, maps, i.e., the diagram
\begin{equation}
\label{O-F-diagram}
\xymatrix{
H^r( \tildescriptM n ({\O}) ) \ar[d] & H^r( \scriptM n ({\O}) ) \ar[r]\ar[d]\ar[l]_-{\sim}
                                                                    & \K 2n-r n {\O}  \ar[d]
\\
H^r( \tildescriptM n (F) )        & H^r( \scriptM n (F) ) \ar[r]\ar[l]_-{\sim}
                                                                    & \K 2n-r n F
}
\end{equation}
commutes.  This shows that we can identify  the complex
$ \scriptM n (\O) $ with a subcomplex of $ \scriptM n (F) $.

\begin{remark}
\label{N-injects}
If the Beilinson--\soule\ conjecture holds for fields of characteristic
zero and for $ \kappa $, then one
proves in the same way as in Proposition~\ref{O-F-injective} that
the map $ \tildescriptM n ({\O}) \to \tildescriptM n (F) $ is injective,
so we can identify
$ \tildescriptM n (\O) $ with a subcomplex of $ \tildescriptM n (F) $.
\end{remark}

\begin{remark}
\label{extension-injective}
If $ F'/F $ is an arbitrary field extension, and $ {\O}'\subset F' $
is a discrete valuation ring with $ {\O} \subset {\O}' \bigcap F $, then there
are obvious maps $ \tildescriptM n ({\O}) \to \tildescriptM n ({\O}') $,
and similarly for $ F $ and $ F' $, as well as for the complexes $ \scriptM n $.
Because the map $ \K 2n-r r F \to \K 2n-r r F' $ is always injective,
the corresponding maps for the
complexes of $ F $ to the ones for $ F' $ is injective, provided
either the Beilinson--\soule\ conjecture is true in general for
fields of characteristic zero, or
$ F' $ is algebraic over $ \Q $.
This is proved like in the proof of Proposition~\ref{O-F-injective},
as the assumptions mean that the necessary maps to $ \K 2n-1 n F $
etc.\ exist and are injective.
Similarly, if the map $ \K 2n-1 n {\O} \to \K 2n-1 n {\O}' $ is injective,
we get the corresponding statement for the complexes for $ {\O} $
and $ {\O}' $ provided the Beilinson--\soule\ conjecture holds for
$ F' $ and $ \kappa' $ (and hence for $ F $ and $ \kappa $).
In particular, all those maps are injective either if the Beilinson--\soule\
conjecture is true for fields of characteristic zero, or $ F' $ is algebraic over $ \Q $.
If this is the case, we shall always view all complexes as being
subcomplexes of the corresponding complexes of $ F' $, and view
all $ K $--groups as being contained in the corresponding $ K $--groups
of $ F' $.
\end{remark}

\begin{remark}
\label{without}
We make a few remarks about the above construction without assuming
the Beilinson--\soule\ conjecture.  There are two places where
it plays a role, namely the vanishing of the spectral sequence~\eqref{messy-ss}
below the row starting with $ \Krloc n n {\O} n-1 $, and the acyclicity
of the subcomplex $ \N n $.

As for the first case, it is well known that $ \K n 1 {\O} = 0 $ for
$ n\geq 2 $.
Therefore the last two columns
(which would correspond to degrees $ n $ and $ n+1 $ for our
complexes) of the spectral sequence below our main row are always
zero.  This means that the maps $ H^r( \scriptM n ({\O}) ) \to \K 2n-r n {\O} $
always exists for $ r = n-1 $ and $ r = n $.
As for the second case, it follows from the proof of \cite[Proposition~3.20]{Jeu95}
that the map $ N_2(K) \to \Sym^2(\KQ) $ given by mapping $ \symb x 2 $
to $ x\tensor x $ is an injection as $ \K 2 1 L $ and $ \K 3 1 L $
are zero for any field $ L $.  Therefore the complex $ \N n (F) $
is acyclic in degrees $ n-1 $ and $ n $.
Because we have the exact sequence
$$ \K 3 1 {\kappa} \to \K 3 2 {\O} \to \K 3 2 F $$
and the first term is zero, we see as in Proposition~\ref{O-F-injective}
that $ \Symb 2 {\O} \subseteq \Symb 2 K $,
$ M_2(\O) \subseteq M_2(K) $,
that $ \N n (R) $ is acyclic in degrees $ n-1 $ and $ n $ and
hence
$ \tM_2(\O) \subseteq \tM_2(K) $.
We also get a commutative diagram
\[
\xymatrix{
H^r(\tildescriptM n (\O)) \ar[d] & H^r(\scriptM n (\O))\ar[l]_\sim\ar[d]\ar[r]
          & \K 2n-r n {\O} \ar[d]
\\
H^r(\tildescriptM n (F))         & H^r(\scriptM n (F))\ar[l]_\sim\ar[r]
          & \K 2n-r n F
}
\]
for $ r=n-1 $ and $ r=n $ without any assumptions.
\end{remark}

\section{Syntomic regulators}
\label{sec:syntomic}

In this Section we briefly recall some parts of the theory of rigid syntomic
regulators, originally due in the affine case to Gros~\cite{Gro94},
as described in detail in~\cite{Bes98a}. Our goal is to describe the
theory in the minimal details required to understand constructions to
follow and to develop certain computational tools that are needed in
later Sections.

Recall that $R$ is a complete discrete valuation ring with quotient
field $K$ of characteristic $0$ and residue field $\kappa$ of
characteristic $p$. We will assume that $\kappa$ is algebraic over the
prime field since this is required for some of the versions of
syntomic cohomology we will be using. All schemes will be separated
and of finite type over their respective bases. We describe as little
as we need of the general theory, sending the interested reader
to~\cite{Bes98a}. All versions of syntomic cohomology are defined as
cohomologies of certain huge complexes. These are needed for the
definition of the regulators but are useless when it comes to
calculations. The cohomology can, however, be realized, using some
auxiliary data, as the cohomology of very explicit complexes, and maps
on cohomology can similarly be realized explicitly. The theory
developed in loc.\ cit.\ guarantees that these explicit maps are
indeed the correct maps and we avoid explicit mentioning of that in the
sequel. 

For the purpose of this work, the version
best suited for computations is the Gros style modified syntomic
cohomology denoted by $\htms$ in~\cite{Bes98a}. This is the weakest
version of syntomic cohomology and all other versions, in particular
$\hsyn$ have natural maps to it~\cite[Proposition~9.5]{Bes98a} which,
by definition, are compatible with Chern classes in algebraic
K-theory. Fortunately, according 
to~\cite[Proposition~\ref{msprop}.\ref{ms3}]{Bes98a}, when $X/R$ is proper
and $2n\ne i,i-1,i-2$, the canonical map $\hsyn^i(X,n)\to
\htms^i(X,n)$ is an isomorphism. Therefore, for the purpose of
detecting syntomic regulators in the cohomology of $\Spec(R)$
working with $\htms$ is just as good as working with $\hsyn$. To
further simplify matters, we only give the description of $\htms$
given certain additional data that may not exist in general but do
exist in our situation.

Suppose first that $X$ is a scheme over a field $\kappa$ of
characteristic $p$.
Following Berthelot we define the rigid complex of
$X$ over $K$ as follows: We choose an open immersion $X\xrightarrow{j}
\Xbar$ into a proper $\kappa$-scheme and a closed immersion
$\Xbar\rightarrow \PF$ into a $p$-adic formal $\vv$-scheme which is smooth in a
neighborhood of $X$.
We remark that in general there may be some
difficulty doing this but in the cases we will consider it will be totally
obvious how to do so.

In the above situation we can, following Berthelot, define the complex
\begin{equation*}
  \rgrig(X/K)_{\PF}:=\rg(\tu{\Xbar}{\PF},
  \jdag\Omega_{\tub{\Xbar}}^\bullet)\;.
\end{equation*}
Here,
the notation
$\tu{\Xbar}{\PF}$ stands for the tube of $\Xbar$ in $\PF$, which
roughly means the space of points in the rigid analytic space
associated to $\PF$ that reduce to a point in $\Xbar$.
The functor $\jdag$ of ``sections of overconvergent support'' goes from
the category of abelian sheaves on $\tu{\Xbar}{\PF}$ to itself and is
defined by
\[
  \jdag (F)=\dirlim_{U} (j_U)_\ast (F|_U)\;,
\]
where the direct limit is over all $U$ which are strict neighborhoods of
$\tu{X}{\PF}$ in $\tu{\Xbar}{\PF}$ in the sense of Berthelot
and $j_U$ is the inclusion of $U$ in $\tu{\Xbar}{\PF}$. We recall that
$U$ is a strict neighborhood if $\{U,\tu{\Xbar}{\PF}-\tu{X}{\PF}\}$ is
an admissible cover of $\tu{\Xbar}{\PF}$ in the sense of rigid analysis.
We have indexed
the complex for simplicity by $\PF$ but we should remember the entire
setup leading up to the definition. In any case, Berthelot shows that
in the derived category of $K$-vector spaces this complex is
independent of all choices, so its cohomology, $\hr(X/K)$, is entirely
well defined. To simplify notation we will drop the $\PF$ subscript from the
notation. In the applications it will be clear which additional data is being
used.

We will often need to let a ($\kappa$-linear) Frobenius act on our
complexes. To do that we will consider a morphism
$\varphi:\Xbar\rightarrow \Xbar$ which is a $\kappa$-linear base
change from a model of $\Xbar$ defined over a finite field with
$q=p^r$ elements of the $r$-th
power of the absolute Frobenius. We insist that $\varphi$ 
preserves $X$. Such a $\varphi$ is called a \emph{Frobenius
endomorphism} of $X$. We then assume that there is a lift $\phi$ of
$\varphi$ to $\PF$. We call $q$ the degree of $\varphi$ and $\phi$. It
is then clear that $\phi$ acts on the rigid complex.

Next we describe the construction of the syntomic complexes. Here we
assume that $X$ is a smooth $\vv$-scheme and that we have an open
immersion $X\rightarrow \Xbar$ into a proper $\vv$-scheme and a closed
immersion $\Xbar\rightarrow P$ into a $\vv$-scheme, smooth in a
neighborhood of $X$, and that
there is a $\vv$-morphism $\phi:P\rightarrow P$ inducing on the
special fiber a Frobenius endomorphism. In
this situation we can clearly embed $X_\kappa$ into $\Xbar_\kappa$ and
this last scheme into the $p$-adic completion $\Phat$ of $P$ to get to
the situation we had when we defined the rigid complex and $\phi$ will
induce a lift of a Frobenius endomorphism.

The given data induces a filtration on the complex $\rgrig(X_\kappa/K)$
defines as follows: let $J$ be the sheaf of ideals defining the generic fiber
$\Xbar_K$ inside $\tu{\Xbar_\kappa}{\Phat}$ and consider the filtration of
$\Omega_{\tu{\Xbar_\kappa}{\Phat}}^\bullet$ given by the complexes
\[
  F_J^n \Omega_{\tu{\Xbar_\kappa}{\Phat}}^\bullet := J^n \Omega^0 \to
  J^{n-1} \Omega^1 \to \cdots\;,
\]
where it is understood that $J^r=\O$ for non positive $r$.
This filtration induces a filtration on the rigid complex by
\[
  F^n \rgrig(X_\kappa/K):=
  \rg(\tu{\Xbar_\kappa}{\Phat},\jdag F_J^n
  \Omega_{\tu{\Xbar_\kappa}{\Phat}}^\bullet)\;.
\]
Berthelot shows that these complexes are again independent of the
additional data up to quasi--isomorphism.
We can now define the Gros style modified syntomic complex to be
the complex
\[
  \rgtms(X,n):= \Cone(F^n \rgrig(X_\kappa/K) \xrightarrow{1- \phi^\ast/q^n}
  \rgrig(X_\kappa/K))[-1]\;.
\]
The map in the cone is a shorthand for the composition of the indicated map
$1-\phi^\ast/q^n$ with the natural map of $F^n \rgrig(X_\kappa/K)$ into
$\rgrig(X_\kappa/K)$. To fix notation for Cones we use the following
sign convention here.
If $f:A^\bullet \rightarrow B^\bullet$, then $\Cone(A^\bullet
\rightarrow B^\bullet)[-1]$  is given in degree $i$ by
\begin{equation}
  \label{sign-cone}
  A^i\oplus B^{i-1},\quad \text{with differential}\;
  d(a,b)=(da,f(a)-db)\; .
\end{equation} 
One can show that $\rgtms(X,n)$ is independent of the
additional data except for the choice of the Frobenius endomorphism
$\varphi$. Here, in the general case one takes a direct limit over all
possible Frobenius endomorphisms as described in
\cite[Definition~\ref{mod-sym}]{Bes98a}. When
$X=\vv$ all the connecting homomorphisms of the limit are
quasi--isomorphisms so we may in fact fix a single $\varphi$.

In \cite{Bes98a} syntomic regulators from the K-theory of $X$ into the
various versions of syntomic cohomology were constructed. For the
cohomology theory we are considering these take the form of a map
\[
c_{i,j}: K_j(X)\to \htms^{2i-j}(X,i)\;.
\]
In this work we will need to consider similar maps in the relative and
multi-relative situations. These were not constructed in loc.\ cit.\
but are constructed in Appendix~\ref{sec:appendix}.

We recall the computation of the regulator on
a part of the K-theory of affine $\vv$-schemes. Suppose $X=\Spec(A)$
is such a scheme. We will give an explicit
description of the rigid and syntomic cohomology of $X$. We can choose
an embedding
of $X$ as an open subset in the projective $P=\Xbar$. 
Suppose $X_\kappa$ is defined in $\Xbar_\kappa$ by the nonvanishing of the
reductions of functions
$h_i$.
Then for $\lambda<1$ we define a rigid space $U_\lambda$ by the conditions
$|h_i|>\lambda$.
The $U_\lambda$ are strict neighborhoods of $\tu{X_\kappa}{\Phat}$ in
$\tu{\Xbar_\kappa}{\Phat}$. It follows that there exists a map
\[
  \dirlim_{\lambda<1} \Gamma(U_\lambda,\Omega^\bullet) \to
  \rgrig(X_\kappa/K)\;.
\]
\begin{proposition}\label{func-simp-rep}
  This map is a quasi--isomorphism. In addition, this
  quasi--isomorphism is functorial with respect of maps of pairs
  $(X,\Xbar)$.
\end{proposition}
\begin{proof}
The first statement follows from the proof of Proposition 1.10
in~\cite{Ber97}. The second statement is a consequence of the
construction of the rigid complexes in~\cite{Bes98a}.
\end{proof}

To obtain the modified syntomic complex,
suppose we have a map $\phi:\Xbar \to \Xbar$ whose reduction is a
Frobenius endomorphism fixing $X_\kappa$.
The ideal $J$ considered above is the $0$ ideal in this
case. We thus get a quasi--isomorphism
\[
  \Cone\left(\dirlim_\lambda \Gamma(U_\lambda,\Omega^{\geqslant n})
   \xrightarrow{1-\phi^\ast/q^n} \dirlim_\lambda
  \Gamma(U_\lambda,\Omega^\bullet)\right)[-1]\isom \rgtms(X,n)\;.
\]
We formally write $U$ for the system of spaces $\{U_\lambda\}$ and
define
\begin{equation}
  \label{eq:convention}
  \Omega^i(U):=\dirlim_\lambda \Gamma(U_\lambda,\Omega^i)\;.
\end{equation}

It follows that
\begin{equation}\label{omegeta}
  \htms^i(X,n)=\frac{\{
    (\omega,\varepsilon):\; \omega\in F^n \Omega^i(U),\; \varepsilon\in
 \Omega^{i-1}(U),\; \dd\omega =0,\; \dd\varepsilon = (1-\phi^\ast/q^n)\omega
 \}}{ \{ (\dd\omega, (1-\phi^\ast/q^n)\omega)-\dd\varepsilon ,\; \omega\in
  F^n\Omega^{i-1}(U),\; \varepsilon\in \Omega^{i-2}(U) \}
  }\;,
\end{equation}
with $F^n \Omega^j(U)=0$ if $n<j$ and $\Omega^j(U)$ otherwise.

As mentioned in the introduction, in many cases syntomic cohomology
becomes isomorphic to rigid cohomology. The normalization of this
isomorphism is perhaps not the obvious one and since the computation
of the regulator depends on the particular normalization, we describe
this here at least in a special case (see~\cite{Bes98a} for a fuller
discussion).
Suppose that $X$ has relative dimension $i-1$ over $\vv$. Suppose
in the description above that $(\omega,\varepsilon)\in\htms^i(X,n)$. We see
that
$\omega=0$ so $\dd\varepsilon =0$. Thus $\varepsilon$ defines a class in
$\hr^{i-1}(X_\kappa/K)$ which is easily seen to be well defined up to
an element of $(1-\varphi^\ast/q^n)F^n \hr^{i-1}(X_\kappa/K)$. When $n\ge i >$
relative dimension of $X$, as
will be the case for us, the map
\begin{multline}\label{twist}
(1-\varphi^\ast/q^n): \hr^{i-1}(X_\kappa/K)/F^n \hr^{i-1}(X_\kappa/K)\\ \to
 \hr^{i-1}(X_\kappa/K)/ (1-\varphi^\ast/q^n)F^n \hr^{i-1}(X_\kappa/K)
\end{multline}
is an isomorphism by
\cite[Proposition~\ref{msprop}.\ref{ms3}]{Bes98a}.
\begin{definition}\label{normalization}
  When $n\ge i >$ relative dimension of $X$, we have a canonical
  isomorphism,
\begin{equation*}
  \htms^i(X,n)\isom \hr^{i-1}(X_\kappa)/F^n ,\quad
  (0,\varepsilon)\mapsto 
  (1-\varphi^\ast/q^n)^{-1}(\text{class of } \varepsilon)\;.
\end{equation*}
\end{definition}
The justification for this normalization requires a longer tour into
the general theory of syntomic cohomology than we would like to
present. The reader may refer to~\cite[Proposition~10.1.3]{Bes98a} for
example. In any case, note that this choice is functorial. We will
make this definition in relative situations as well.

We now describe the regulator in this special case. First of all, consider
$f\in A^*$. If $\bar{f}$ is the reduction of $f$ one finds that
$\varphi^\ast \bar{f} = \bar{f}^q$ and therefore that $f_0:= f^q/\phi^\ast
f$ is congruent to $1$ mod the maximal ideal of $\vv$. One can deduce from
that that the function $\log(f_0)$ is analytic on some $U_\lambda$.
\begin{lemma}\label{reg-on-func} (\cite[Proposition~10.3]{Bes98a})
The
syntomic regulator sends the class of $f$ in $K_1(X)$ to the cohomology
class of $(\dlog f,\log (f_0)/q)$ in the
representation~\eqref{omegeta} of $\htms^1(X,1)$.
\end{lemma}
The value of the regulator on a cup product $f_1\cup \cdots \cup f_r$ in $ K_n(X)$
is the cup product of the regulators of the $f_j$'s, so it
is enough to describe the cup product on syntomic cohomology. This is
given, in the notation of \eqref{omegeta}, by any of the formulas, depending on
the parameter $\gamma$,
\begin{equation}
\label{cup-product}
    \begin{split}
      (\omega_1,\varepsilon_1)\cup (\omega_2,\varepsilon_2) = &\Big(\omega_1\wedge
      \omega_2,\\&
      \varepsilon_1 \wedge \left(\gamma +(1-\gamma)
      \left(\frac{\phi^\ast}{q^j}\right)^k\right) \omega_2\\
    +(-1)^{\operatorname{deg}\omega_1}
    &\left(\left((1-\gamma)+ \gamma
      \left(\frac{\phi^\ast}{q^i}\right)^k\right)\omega_1\right)
      \wedge \varepsilon_2\Big)\;.
    \end{split}
\end{equation}

We need to describe the pullback map in syntomic cohomology in certain
special situations. Suppose that $X$ is an affine scheme and $f:Y\to X$ is a
closed embedding on an
affine subscheme and choose the same auxiliary data for $X$ as before. We
may compactify $Y$ by embedding it into its closure $\Ybar$ in $\Xbar$. The
difficulty in describing the pullback map from $X$ to $Y$ is that the
lift of Frobenius morphism $\phi$ will not preserve $\Ybar$ in general.
Note however that we may and do assume that $\varphi$ preserves $Y_\kappa$.
The way to overcome the difficulty is as follows: we use the embedding of
$\Ybar$ into $\Xbar$ to compute the syntomic complex of $Y$.
This gives us the following model for $\rgtms(Y,n)$,
\[
  \rgtms(Y,n)\isom\Cone\left(
  \rg(\tu{\Ybar_\kappa}{\Xbarhat},\jdag F_J^n
  \Omega_{\tu{\Ybar_\kappa}{\Xbarhat}}^\bullet)
  \xrightarrow{1- \phi^\ast/q^n}
  \rg(\tu{\Ybar_\kappa}{\Xbarhat},\jdag
  \Omega_{\tu{\Ybar_\kappa}{\Xbarhat}}^\bullet)\right)
\!\![-1]\, ,
\]
and the pullback map is now simply obtained by restriction to the tube
$\tu{\Ybar_\kappa}{\Xbarhat}$. Here $J$ is the ideal of $\Ybar_K$ in
$\Xbar_K$.

Suppose now that $Y$ is of relative dimension $i-1$ over $\vv$ and that we
are given an element of $\htms^i(X,n)$ represented by the pair of forms
$(\omega,\varepsilon)$ as in \eqref{omegeta}. We would like to study the pullback of this
element to $Y$, identified with an element of $\hr^{i-1}(Y_\kappa/K)/F^n$
as in Definition~\ref{normalization}. Note that this pullback does not
factor through $\hr^{i-1}(X_\kappa/K)/F^n$
because $X$ is of higher dimension than $Y$ in general. Recalling the sets
$U_\lambda$ we see that for each $\lambda$ the set $U_\lambda \cap
\tu{\Ybar_\kappa}{\Xbarhat}$ is a strict neighborhood of
$\tu{Y_\kappa}{\Xbarhat}$ in $\tu{\Ybar_\kappa}{\Xbarhat}$.
It follows that we may factor the map $\rgrig(X_\kappa/K) \to
\rgrig(Y_\kappa/K)$, respectively
$F^n\rgrig(X_\kappa/K) \to F^n \rgrig(Y_\kappa/K)$ as
\[
  \dirlim_\lambda \Gamma(U_\lambda,\Omega^\bullet)\to
  \dirlim_\lambda \Gamma(U_\lambda \cap \tu{\Ybar_\kappa}{\Xbarhat}
  ,\Omega^\bullet)\to
  \rg(\tu{\Ybar_\kappa}{\Xbarhat},\jdag
  \Omega_{\tu{\Ybar_\kappa}{\Xbarhat}}^\bullet)\;,
\]
respectively with $\Omega^\bullet$ replaced by $F_J^n \Omega^\bullet$.
We may therefore factor the map of syntomic complexes $\rgtms(X,n) \to
\rgtms(Y,n)$ via
\begin{equation}
  \label{middlecone}
  \Cone\left(
  \dirlim_\lambda \Gamma(U_\lambda \cap \tu{\Ybar_\kappa}{\Xbarhat}
  ,F_J^n \Omega^\bullet)
  \xrightarrow{1- \phi^\ast/q^n}
  \dirlim_\lambda \Gamma(U_\lambda \cap \tu{\Ybar_\kappa}{\Xbarhat}
  ,\Omega^\bullet)\right)[-1]\;.
\end{equation}
\begin{lemma}
  \label{theta}
  In the situation described above let $\theta\in
  \dirlim_\lambda \Gamma(U_\lambda \cap \tu{\Ybar_\kappa}{\Xbarhat}
  ,F_J^n \Omega^{i-1})$ be such that $\dd\theta =
  \omega|_\tu{\Ybar_\kappa}{\Xbarhat}$.
  Then the image of $f^\ast (\omega,\varepsilon)$ in
  $$
  \hr^{i-1}(Y_\kappa/K)/ (1-\varphi^\ast/q^n)F^n \hr^{i-1}(Y_\kappa/K)
  $$
  is
  the same as the image of $\varepsilon|_\tu{\Ybar_\kappa}{\Xbarhat}-
  (1-\phi^\ast/q^n)\theta\in \dirlim_\lambda \Gamma(U_\lambda
  \cap \tu{\Ybar_\kappa}{\Xbarhat},\Omega^{i-1})$ in the same group.  
\end{lemma}
\begin{proof}
Subtract the boundary $\dd(\theta,0)=(\dd\theta, (1-\phi^\ast/q^n)\theta)$
from $(\omega,\eta)|_\tu{\Ybar_\kappa}{\Xbarhat}$ in
\eqref{middlecone}.
\end{proof}

Finally we specialize even further and show how to compute the
difference of the pullbacks at two nearby points. We assume that
$i=n$. Consider a 
situation where we are given
an affine $X$ together with a smooth affine map $\pi: X\to B$ to another affine
scheme $B$ smooth over $\vv$. Suppose that
$\pi$ extends to $\bar{\pi}: \Xbar\to B$.
Let $z'\in B(\kappa)$ and let
$D$ be the rigid analytic space of all points in
of $B$ reducing to $z'$ (this is the residue disc of $z'$ in the
terminology of Coleman). For any $z\in D(K)$ let $f_z$ be the embedding of
$Y_z:=\pi^{-1}(z)$ in $X$. The $Y_z$ for $z\in D(K)$ have a common
reduction which we denote by $Y_\kappa$, and the $\bar{\pi}^{-1}(z)$
have a common reduction $\Ybar_\kappa$, which is a compactification of
$Y_\kappa$. Finally, the tube $\tu{\Ybar_\kappa}{\Xbarhat}$ is simply
$\bar{\pi}^{-1}(D)$.
\begin{proposition}\label{pullback}
  In the situation described above let $(\omega,\varepsilon)$
  represent a class in $\htms^i(X,i)$. Let $z_1$, $z_2\in D(K)$, let
  $J_j$ be the ideal defining $Y_{z_j}$ and let
 $\theta_j$ in 
 $ \dirlim_\lambda \Gamma(U_\lambda \cap \bar{\pi}^{-1}(D)
  ,F_{J_j}^n \Omega^{i-1})$ be such that $\dd\theta_j =
  \omega|_{\bar{\pi}^{-1}(D)}$.
  Then, the difference of the images in $\hr^{i-1}(Y_\kappa/K)$ of the
  pullbacks $f_{z_j}^\ast (\omega,\varepsilon)$, $j=1,2$, 
  is the image
  of $\theta_2-\theta_1\in \dirlim_\lambda \Gamma(U_\lambda
  \cap \bar{\pi}^{-1}(D),\Omega^{i-1})$.
\end{proposition}
\begin{proof}
As seen in Lemma~\ref{theta}, the image in
$ \!\hr^{i-1}(Y_\kappa/K)/ (1-\varphi^\ast/q^n)F^i \hr^{i-1}(Y_\kappa/K)$
is the image of $(1-\varphi^\ast/q^i)(\theta_2-\theta_1)$ and the
result thus follows from Definition~\ref{normalization} of the
image in $\hr^{i-1}(Y_\kappa/K)$.
\end{proof}

In the situation as above the following immediate
Corollary will also be useful.
For any point $x\in D(K)$ the fiber $\bar{\pi}^{-1}(x)$
is a lift of $\Ybar_\kappa$ and therefore the rigid cohomology of
$Y_\kappa$ can be computed as the cohomology of $\dirlim_\lambda
\Gamma(U_\lambda\cap \bar{\pi}^{-1}(x),\Omega^\bullet)$.
\begin{corollary}\label{pullback1}
  In $\hr^{i-1}(Y_\kappa/K)$ the difference $f_{z_1}^\ast
  (\omega,\varepsilon)-f_{z_2}^\ast(\omega,\varepsilon)$ is the image of
  $(\theta_2-\theta_1)|_{U_\lambda\cap \bar{\pi}^{-1}(x)}\in
  \dirlim_\lambda
  \Gamma(U_\lambda\cap \bar{\pi}^{-1}(x),\Omega^{i-1})$.
\end{corollary}

All of our considerations are also valid for the cohomology of diagrams of
schemes, and in particular for the relative and multi--relative cohomologies
that will be considered in Sections to come.
\section{The integration down process}
\label{sec:down}

A key ingredient in the computation of the regulator will be a
functional on rigid cohomology obtained by repeated integration, which
we now go on to describe.

Let $\kappa$ be a field of characteristic $p$.
Set $X^n=(\PP_\kappa^1-\{t=1\})^n$. Let $B$ be an affine
$\kappa$-variety. Let $Y$ be an open
affine subset of $X^n\times B$.
Let $\bbox^n$ be the subset of
$Y$ where at
least one of the coordinates is either $0$ or $\infty$. We write
cohomology relative to $\bbox^n$ to mean the multi relative cohomology
taken in exactly the same way as was done in Section~\ref{sec:kstuff} for
$K$-theory. We would like to write an explicit complex computing
the multi relative rigid cohomology of $Y$.

We first choose the rigid data $X^n\inject (\PP_\kappa^1)^n\inject
(\PP_{\spf(\vv)}^1)^n$ and $B\inject \Bbar\inject P_B$ where $\Bbar$
and $P_B$ are projective spaces of some degree over $\kappa$ and
$\spf(\vv)$ respectively. We thus obtain a rigid
datum for $Y$ as well. As in Section~\ref{sec:syntomic} we obtain for
$\lambda<1$ a certain inverse system of rigid space
$U=\{U_\lambda^Y\}$. We know that there
exists a canonical quasi-isomorphism
$\Omega^\bullet(U^Y):=\dirlim_{\lambda\rightarrow 1} 
\Omega_{U_\lambda^Y}^\bullet \rightarrow
\rgrig (X^n/K)$. Similarly, the complexes of rigid forms on the
subspace of $U_\lambda$ cut out by
equations of the form $t_i=0$ or $t_i=\infty$ 
are quasi--isomorphic to the  rigid complexes of the various components
in $\bbox^n$.

As discussed (at length) in Appendix~\ref{sec:appendix},
we can now write a complex quasi-isomorphic to the
multi relative  $\rgrig(X^n;\bbox^n/K)$ by taking iterated cones on
the complexes above with respect to the restriction maps to $t\in
\{0,\infty\}$. We want to do the ``battle of signs'' correctly to write this
iterated cone as a simple complex. This can be done as follows:
For $ 0\leq j\leq n
$ consider all strictly increasing functions $ f:[1,\ldots
,j]\rightarrow [1,\ldots ,n] $. To such a function $ f $ we
associate the subspace 
$$ Y_{f}:=\{(x_{1},\ldots ,x_{n})\in Y:\; x_{i}\in \{0,\infty \},i\notin \im f\} .$$
We can similarly define
rigid spaces $U_{f,\lambda}^Y$ forming an inverse system $U_f^Y$ and
like \eqref{eq:convention} we can formally define complexes of
differential forms $\Omega^\bullet(U_f^Y)$.
Let us call $ n-j $
the degree of $ f $ (this includes the empty function $\emptyset$ with
degree $n$), which is the same as the codimension of $Y_f$. The
complex computing our multi relative
cohomology can then be written as $ \oplus _{i+\deg f=k}\Omega
^{i}(U_{f}^Y) $ in degree $ k $ . Let us write an element in the $ f
$ component of this complex as $ (\omega ,f) $. Then the
differential is given by
\begin{equation}\label{multidif}
  \dd(\omega,f) = (\dd\omega,f)- (-1)^{\deg \omega} \sum_g \chi
  (f,g)(\omega |_{Y_{g}},g)
\end{equation}
where
\begin{equation*}
  \chi (f,g)=
  \begin{cases}
    (-1)^{f(r)+r} & \im f=\im g{\dotcup{f(r)}}
\\
    0 & \text{otherwise,}
  \end{cases}
\end{equation*}
and where here $\dotcup$ denotes disjoint union.
\begin{definition}\label{mulcompl}
The complex above is denoted $\Omega^\bullet(U^Y;\bbox^n)$.  We let
$F^j \Omega^\bullet(U^Y;\bbox^n)$ be the subcomplex having
$\oplus _{i+\deg f=k,i\ge j}\Omega^{i}(U_{f}^Y)$  in degree $k$.
\end{definition}
The following Lemma is mostly an exercise in sign fixing.
\begin{lemma}\label{explicitrgrig}
  The complexes $F^j \Omega^\bullet(U^Y;\bbox^n)$ and
  $F^j \rgrig(Y;\bbox^n/K)$ are  quasi-isomorphic.
\end{lemma}
\begin{proof}
We prove this without the filtrations. The result for the filtered
part is then clear. First we note that by
Proposition~\ref{func-simp-rep}
this complex is quasi-isomorphic to the corresponding complex with $\Omega
^{i}(U_{f}^Y)$ replaced by the degree $i$ part of $\rgrig(Y_f/K)$ which
we now denote by $\Gamma^i$ for simplicity. Consider now the double
complex introduced in the appendix (compare
 \eqref{final-complex} and \eqref{eq:explicit-differential})
For $\beta: [1,k]\to [1,n]$  an increasing function define
$Y_\beta=\{(x_{1},\ldots ,x_{n})\in Y:\; x_{\beta(i)}\in \{0,\infty
\}\}$. Then the complex in degree $q$ is $\oplus_{k+|\beta|=q} \Gamma^k
(Y_\beta)$. We can again write elements there as pairs
$(\omega,\beta)$ and the differential is defined by

\begin{equation*}
  \dd(\omega ,\beta)=(\dd\omega ,\beta)+(-1)^q (-1)^{|\beta|} \sum
  _{\beta'}\chi' (\beta,\beta')(\omega |_{Y_{\beta'}},\beta')
\end{equation*}
where
\begin{equation*}
  \chi' (\beta,\beta')=
  \begin{cases}
    (-1)^{r} & \im \beta'=\im \beta\dotcup {\beta'(r)}
\\
    0 & \text{otherwise,}
  \end{cases}
\end{equation*}
Now we want to switch to a dual point of view.
The relation is as follows:  For $f$ a function as
in the Lemma we define $\beta(f)$ to be the
increasing function enumerating $[1,\ldots,n]-\im f$. Then we have
$|\beta(f)|=\deg f$ and $Y_f=Y_{\beta(f)}$. The key thing to check is
the following: If $g$ is obtained from $f$ by deleting
$f(r)$, then $\beta=\beta(f)$ is obtained from $\beta'=\beta(g)$ by
deleting $\beta'(f(r)-r+1)$. This easily gives the result.
\end{proof}

Consider now the case where $B=\Spec(\kappa)$ and $Y=X^n$.
Then the
relative rigid cohomology $\hr^n(X^n;\bbox^n/K)$ is well known to be
isomorphic to $K$. We can explicitly describe this isomorphism.
The basic
idea (compare\cite{Gro94}) is of iterated integration between $ 0 $
and $ \infty $. We can take $U_\lambda^Y$ to be $U_\lambda^n$ where
$U_\lambda$ denotes the space $\PP_K^1 -\{|t-1|<\lambda\}$ and
$U_\lambda^n$ is the $n$th power of $U_\lambda$.
Let $ (\omega ,f) \in \Omega^n(U^n;\bbox^n)$ and suppose that
\begin{equation*}
  \omega =G(t_{f(1)},\ldots ,t_{f(j)})dt_{f(1)}\wedge \cdots \wedge
dt_{f(j)}
\end{equation*}
(here the ordering is critical). Define
\begin{equation*}
  \pi (\omega ,f)=\int _{\infty }^{0}\cdots \int _{\infty
    }^{0}G(t_{f(1)},\ldots ,t_{f(j)})dt_{f(1)}\cdots dt_{f(j)}.
\end{equation*}

Notice that now the order is not critical and we can integrate in
whatever order we want.  Let $\hdr(U^n;\bbox^n)$ be the homology of
$\Omega^\bullet(U^n;\bbox^n)$ We have the following.

\begin{lemma}\label{lem:going-down}
  There is a unique isomorphism $ \hdr ^{n}(U^{n};\bbox^n )\rightarrow
  K $ normalized by the condition that on the class of a closed form
  $ (\omega ,f) $ with $ \deg f=0 $ it is given by $ \pi (\omega
  ,f) $. This functional is given as follows: consider a form $
  (\eta ,g)$ where $g$ has degree $m$ and of the $m$
  coordinates on $U_g^n$ which are fixed $i$ are fixed to be $\infty$.
  Set
\begin{equation*}
  \Pi ((\eta ,g))=(-1)^{i+\sum (g(k)+k)}\pi (\eta ,g).
\end{equation*}
Then the functional is
given by the restriction of the $ K $-linear extension of $ \Pi $
to closed forms.
\end{lemma}
\begin{proof}
  We can find a form $ (\omega ,f) $, with $\deg(f)=0$, whose
  cohomology class is non-trivial.  The required isomorphism is
  determined by its value on such a form and is therefore unique. To
  show that $\Pi$ 
  provided the required map, we only need to show, in view of the fact
  that $ \Pi (\omega ,f)=\pi (\omega ,f) $ when $\deg(f)=0$, that it
  kills exact
  forms. The exact forms are spanned by forms
\begin{align*}
  &\phantom{=}\dd\big(F(t_{g(1)},\ldots ,t_{g(j)})\cdot dt_{g(1)}\wedge
  \cdots \wedge
  \widehat{dt_{g(k)}}\wedge \cdots \wedge dt_{g(j)},g\big)\\
  &=(-1)^{k-1}\left(\frac{\partial }{\partial
      t_{g(k)}}F\cdot dt_{g(1)}
  \wedge \cdots \wedge dt_{g(j)},g\right)\\
  &\phantom{=}+(-1)^{j+g(k)+k}\big(F|_{t_{g(k)}\in \left\{ 0,\infty
  \right\} }
  \cdot dt_{g(1)}\wedge \cdots \wedge \widehat{dt_{g(k)}}\wedge \cdots
  \wedge dt_{g(j)},h\big)
\end{align*}
with $h$ obtained from $g$ by removing the $k$-th value
and $F$ a function on a component of $U_h^n$ with $i$ coordinates forced
to $ \infty $.  Notice that
\begin{align*}
  &\phantom{=}\pi \left(\frac{\partial }{\partial t_{g(k)}}F \cdot
    dt_{g(1)}\wedge
  \cdots \wedge dt_{g(j)},g\right)\\
  &=\pi \big(F|_{t_{g(k)}=0}\cdot dt_{g(1)}\wedge \cdots \wedge
  \widehat{dt_{g(k)}}
  \wedge \cdots \wedge dt_{g(j)},h\big) \\
  &\phantom{=}-\pi \big(F|_{t_{g(k)}=\infty }\cdot dt_{g(1)}\wedge
  \cdots \wedge
  \widehat{dt_{g(k)}}\wedge \cdots \wedge dt_{g(j)},h\big)
\end{align*}
because in the computation of $\pi$ we can begin the integration
on the $g(k)$ coordinate. Now call the two terms on the right hand
side of the last equation $\alpha_0$ and $\alpha_\infty$ respectively,
and write the possible sign in the definition of $\Pi$ as
$\operatorname{Sign}(g,i)$. We therefore obtain
\begin{align*}
&\phantom{=}\Pi \Big(\dd\big(F(t_{g(1)},\ldots ,t_{g(j)})\cdot dt_{g(1)}\wedge
 \cdots \wedge \widehat{dt_{g(k)}}\wedge \cdots \wedge dt_{g(j)},g\big)\Big)\\
&=(-1)^{k-1}\operatorname{Sign}(g,i)(\alpha_0-\alpha_\infty) +
 (-1)^{j+g(k)+k} (\operatorname{Sign}(h,i)\alpha_0+
 \operatorname{Sign}(h,i+1)\alpha_\infty)\;.
\end{align*}
Thus, clearly, to make this cancel, we need to choose
$\operatorname{Sign}(g,i)=(-1)^{\operatorname{sign}(g)+i}$ with
$\operatorname{sign}(g)$ satisfying the relation
\begin{equation*}
  \operatorname{sign}(g)+k-1\equiv \operatorname{sign}(h)+j+g(k)+k+1
  \pmod{2}
\end{equation*}
(the last $1$ is there to make this alternating) when $h$ is obtained
from $g$ by deleting $g(k)$. After cancellations this becomes
\begin{equation*}
  \operatorname{sign}(g)\equiv \operatorname{sign}(h)+g(k)+j
  \pmod{2}\;.
\end{equation*}
It is easily seen that $ \operatorname{sign}(g)=\sum (g(k)+k)$
satisfies this condition, which completes the proof.
\end{proof}

Like in $K$--theory (see the derivation of \eqref{eq:downink}
or~\cite[page 220]{Jeu95}) the isomorphism $\hr^n(X^n;\bbox^n/K)\isom K$ can
be obtained by a repeated application of boundary maps. At each stage
there is a choice of signs.
Here we have taken the
approach of writing down the isomorphism $\Pi$ directly and we would
now like to know how it can be obtained using boundary maps.

We have a short exact sequence
\begin{equation*}
  0\to\hr^{n-1}(X^n;\bbox^{n-1}/K)\to\hr^{n-1}
  (\bbox^n;\bbox^{n-1}/K)\to\hr^n(X^n;\bbox^n/K) \to 0\;, 
\end{equation*}
and an isomorphism
$ \hr^{n-1}(X^n;\bbox^{n-1}/K) \iso \hr^{n-1}(X^{n-1};\bbox^{n-1}/K)$
under pullback. It follows from this that the following two composed
maps are isomorphisms:
\begin{equation*}
  \hr^{n-1}(X^{n-1};\bbox^{n-1}/K) \to
  \hr^{n-1}(\bbox^n;\bbox^{n-1}/K)\to\hr^n(X^n;\bbox^n/K)\;.
\end{equation*}
Here, the two maps come from the two possible choices of the
embeddings of $(X^{n-1};\bbox^{n-1})$ in $(\bbox^n;\bbox^{n-1})$ as
either $x_n=0$ or $x_n=\infty$ and they differ by a minus
sign. Iterating this we get an isomorphism
\begin{equation}
  \label{eq:iteratedisom}
  K= \hr^0(*/K) \xrightarrow{\sim} \hr^1(X^1;\bbox^1/K)
  \xrightarrow{\sim}\cdots
  \xrightarrow{\sim} \hr^n(X^n;\bbox^n/K)\;.
\end{equation}

\begin{proposition}\label{newdownsign}
  The composed map $K\xrightarrow{\eqref{eq:iteratedisom}}
  \hr^n(X^n;\bbox^n/K)\xrightarrow{\Pi} K$
  is the identity provided at each stage we choose the embedding as
  $x_i=0$.
\end{proposition}

\begin{proof}
From the proof of Proposition~\ref{explicit-qi}
it is not difficult to get the following explicit description of the map
$\Omega^\bullet(\bbox^n;\bbox^{n-1}) \to \Omega^{\bullet+1} (X^n;\bbox^n)$
(dual to the map $ \ab C_{\bullet} \to \ab Y_{\bullet} [-1] $ in
the notation of the proof of Proposition~\ref{explicit-qi}):
it is simply given by
$(\omega,f)\mapsto(\omega,f)$ where $f:[1,\ldots,j]\to [1,\ldots,n-1]$
is considered on the right as a function  $f:[1,\ldots,j]\to
[1,\ldots,n]$. It follows that the map \eqref{eq:iteratedisom} with
the choice of signs as in the Proposition corresponds to the map
sending $\alpha\in K$ to $(\alpha,\emptyset)$ on the component with
$x_i=0$ for all $i$. Applying $\Pi$ to this we get $\alpha$.
\end{proof}

Now comes a crucial point. In applications we will want to consider
the cohomology not of $(X^{n};\bbox^n)$
but rather of an open subset $Y$ obtained from $X^n$ by removing
subsets of the form $\{t_j=u\}$ with $u\in \kappa^\ast$ (e.g., sets of
the form $X_{\loc}^n$ as in Section~\ref{sec:kstuff}. That means that
it is no longer
possible to perform the integrals needed to construct $ \Pi $ (and
of course the isomorphism that $ \Pi $ represents does not exist).
It is sometimes possible, however, to replace the integral by a
Coleman integral. We want to show that when this is possible it
corresponds to an operation which can be made sense out of in general.

\begin{lemma}
  There is a short exact sequence
\begin{equation*}
  0\rightarrow \hr ^{n}(X^{n};\bbox^n /K)\rightarrow 
  \hr^{n}(Y;Y\cap \bbox^n
  /K)\rightarrow E\rightarrow 0\;,
\end{equation*}
 where Frobenius acts on $\hr
^{n}(X^{n};\bbox^n)$ trivially and on $ E $ with strictly positive
weights.
\end{lemma}

\begin{proof}
Write $H^i(\ast)$ for $\hr^i(\ast/K)$.
From the diagram of pairs 
\[  
  (X^n-Y;\bbox^n-Y)\to (X^n;\bbox^n)\leftarrow
  (Y;Y\cap \bbox^n)
\]
We get the standard long exact sequence
\begin{equation*}
  \cdots \to
  H^n(X^n;\bbox^n)\to H^n(Y,Y\cap \bbox^n)
  \to H_{(X^n-Y;\bbox^n-Y)}^{n+1}(X^n;\bbox^n)\to \cdots\;.
\end{equation*}
The action of Frobenius on $ \hr ^{n}(X^{n};\bbox^n /K) $ is
trivial because the isomorphism $ \hr ^{n}(X^{n};\bbox^n /K)\isom \hr
^{0}(\operatorname{pt}/K)$ is Frobenius equivariant.
To prove the Lemma we need to show that the first arrow indicated in the
diagram is not $0$ while the last term has strictly positive weights.
The first assertion follows because it is easy to see that the same
integration process described in Lemma~\ref{lem:going-down} also
vanishes on exact relative forms
on the pair $(Y;\bbox^n)$. It remains to show the statement about the
weights. To do that we ``peel off'' the relativity step by step:
We have a long exact sequence
\begin{multline*}
  \cdots \to H_{(\{t_n\in \{0,\infty\}\}-Y;\bbox^{n-1}-Y)}^{n}
  (\{t_n\in \{0,\infty\}\};\bbox^{n-1})
  \to
  \cr
  H_{(X^n-Y;\bbox^n-Y)}^{n+1}(X^n;\bbox^n)
  \to H_{(X^n-Y;\bbox^{n-1}-Y)}^{n+1}(X^n;\bbox^{n-1})\to \cdots
\end{multline*}
and the two terms on the sides fit into similar sequences. The key
observation is that the degree of the cohomology is always one more than
the dimension of the space. The final ``building blocks'' are of the form
$ H_{X^i-Y}^{i+1}(X^i)$. By~\cite{Chi96} such a term has weights between $i+1$
and $2i$ (because $X^i-Y$ is always of codimension $1$ by our assumptions)
except that
the term with $i=0$ clearly vanishes. Thus all terms have positive weights.
\end{proof}
\begin{corollary}\label{cor:pi-col}
  Let $ M\subset \hr ^{n}(Y;\bbox^n /K) $ be any Frobenius invariant
  subspace containing $ \hr ^{n}(X^{n};\bbox^n /K) $. Then there
  exists a unique $ K $-linear functional 
  $$ \widetilde{\Pi} _{M} : M\rightarrow K $$
  which is fixed under Frobenius and
  coincides with the functional induced by $ \Pi $ on $ \hr
  ^{n}(X^{n};\bbox^n /K) $.
\end{corollary}
Of course the conclusion is also true with $ M=\hr ^{n}(Y;\bbox^n /K)
$ in which case we will denote $ \widetilde{\Pi }_{M} $ simply by
$ \widetilde{\Pi } $.  We will need the uniqueness statement,
however, for possibly different subspaces.

The map $\widetilde{\Pi}$ gives a splitting of $V=\hr^n(Y;\bbox^n/K)$ into
a direct sum $V=K\oplus E$ as a $\phi$-module, where $E$ has no
$\phi$-fixed vectors. We will need a certain result about $\phi$-modules with
such a structure.
\begin{lemma}\label{pairi}
  Let $V_i=K\oplus E_i$ for $i=1,2,3$ be three $\phi$-modules such that
  $\phi$ has no invariant vectors on $E_i$ for each $i$ and on $E_1\otimes
  E_2$. Let $\Pi_i:V_i\to
  K$ be the natural projection. Suppose there is a
  $\phi$-equivariant pairing $\pair{~,~}:V_1\otimes V_2 \to V_3$ which
  gives the usual multiplication when restricted to $K\otimes K$.
  Then we have $\Pi_3(\pair{x_1,x_2})=\Pi_1(x_1)\cdot \Pi_2(x_2)$.
\end{lemma}

\begin{proof}
  The conditions of the Lemma imply that the algebraic multiplicity of $1$
  as an eigenvalue of $\phi$ on $V_1\otimes V_2$ is $1$. It follows that
  the space of $\phi$-invariant functionals on $V_1\otimes V_2$ is
  $1$-dimensional. Therefore the statement of the Lemma has to be true up
  to a multiplicative constant. This constant has to be $1$ because the
  statement is true for $x_i=1$.
\end{proof}

By assumption all components of $Y$ are affine. We can therefore
compute relative rigid cohomology using the complex
$\Omega^\bullet(U^Y;\bbox^n)$ of Definition~\ref{mulcompl}.
\begin{definition}
  A relative form in $\Omega^n(U^Y;\bbox^n)$ is called Coleman
  integrable if for each of its component $ (\omega ,f) $ the
  expression defining $ \pi (\omega ,f) $ makes sense when we
  replace ordinary integration with Coleman integration. If $ x $ is
  such a form we denote by $ \Pi _{\Cole }(x) $ the expression
  derived from the $ \pi (\omega ,f) $ as in Lemma~\ref{lem:going-down}.
\end{definition}

\begin{lemma}
  Coleman integrable relative forms form a subspace of
  $\Omega^n(U^Y;\bbox^n)$ which is closed
  under $ \phi $.  Exact forms and forms extending to $U^{n}$
  are Coleman integrable. The functional $ \Pi _{\Cole } $ is
  $\phi$-invariant.
\end{lemma}

\begin{proof}
The only thing which possibly requires proof is the fact that if $
  x $ is a relative form which is Coleman integrable, then so is $
  \phi ^{*}(x) $ and $ \Pi _{\Cole }(\phi ^{*}(x))=\Pi _{\Cole }(x)
  $. This is an easy explicit computation. We may assume that
  \begin{equation*}
    x=(G(t_{f(1)},\ldots ,t_{f(j)})\cdot dt_{f(1)}\wedge \cdots \wedge
    dt_{f(j)},f)\;.
  \end{equation*}
  Then
  \begin{equation*}
    \phi ^{*}(x)=(G(t_{f(1)}^{q},\ldots
    ,t_{f(j)}^{q})\cdot d(t^{q}_{f(1)})\wedge \cdots \wedge
    d(t^{q}_{f(j)})\;.
  \end{equation*}
  The assumption that $ x $ is Coleman integrable means the
  following: There is a function $ F_{1}(t_{f(1)},\ldots ,t_{f(j)})
  $ which is a Coleman function in the first variable and such that
  \begin{equation*}
    \frac{\partial }{\partial t_{f(1)}}F_{1}=G\;.
  \end{equation*}
  Setting $
  G_{1}(t_{f(2)},\ldots ,t_{f(j)})=F_{1}|_{t_{f(1)}=\infty
    }^{t_{f(1)}=0} $ we can find a function $ F_{2}(t_{f(2)},\ldots
  ,t_{f(j)}) $ which is again Coleman in the first variable such that
  \begin{equation*}
    \frac{\partial }{\partial t_{f(2)}}F_{2}=G_{1}
  \end{equation*}
 and we continue like
  this until we reach $ G_{j} $ which is just a number equaling $ \Pi
  _{\Cole }(x) $. Now we start with
  \begin{equation*}
    \widetilde{G}(t_{f(1)},\ldots ,t_{f(j)})=G(t_{f(1)}^{q},\ldots
    ,t_{f(j)}^{q})qt_{f(1)}^{q-1}\cdots qt_{f(j)}^{q-1}\;.
  \end{equation*}
The functoriality of the Coleman integral implies that we may take
  \begin{equation*}
    \widetilde{F}_{1}(t_{f(1)},\ldots ,t_{f(j)})=F_{1}(t_{f(1)}^{q},\ldots
    ,t_{f(j)}^{q})qt_{f(2)}^{q-1}\cdots qt_{f(j)}^{q-1}\;.
  \end{equation*}
Then, as $0^{q}=0$ and $ \infty ^{q}=\infty $ we get
  \begin{equation*}
    \widetilde{G}_{1}(t_{f(2)},\ldots ,t_{f(j)})=G_{1}(t_{f(2)}^{q},\ldots
    ,t_{f(j)}^{q})qt_{f(2)}^{q-1}\cdots qt_{f(j)}
  \end{equation*}
and we can continue this process until we find
$\Pi_{\Cole}(\phi^{*}(x))=\Pi_{\Cole}(x) $.
\end{proof}
We call a cohomology class in $ \hr ^{n}(Y;\bbox^n /K) $ Coleman
integrable if it is represented by a Coleman integrable form. Let $
M_{\Cole } $ denote the space of Coleman integrable cohomology classes.
It is an immediate consequence of the above that $ \Pi _{\Cole } $
induces a functional $ M_{\Cole }\rightarrow K $ which is Frobenius
invariant. By the construction it is also clear that $ \Pi _{\Cole }
$ is just $ \Pi $ on forms that extend to $ U^{n} $. By
corollary~\ref{cor:pi-col} this functional must coincide with the
restriction to $ M_{\Cole }
$ of $ \widetilde{\Pi } $.  We therefore obtain

\begin{proposition}
  For any Coleman integrable form $ x $ representing a cohomology
  class $ [x] $ we have $ \Pi _{\Cole }(x)=\widetilde{\Pi }([x]) $.
\end{proposition}
\section{Regulators for special elements}

Recall that we have universal symbols \eqref{eq:Ssubn}
\begin{equation*}
   \symb S n  \in  \Krloc n n G n-1
\end{equation*}
where $G=\Spec(\Z[S,S^{-1},(1-S)^{-1}])$. Now let $B=\Spec(R[S,S^{-1},
(1-S)^{-1}]) $. Pulling back via the canonical map $B\to G$ we obtain
elements, for which we retain the notation,
\begin{equation*}
   \symb S n  \in  \Krloc n n B n-1 \;.
\end{equation*}
In this Section we obtain some information on the regulators
\begin{equation*}
  \reg(\symb S n )\in  \htms^n(\xlocb B n-1 , n) \;.  
\end{equation*}

We embed $X_{B,\textup{loc}}^{n-1}$ in $P=(\PP_R^1)^{n-1} \times \PP_R^1$
($B$ is mapped
to the last coordinate). Taking the special fiber corresponds to the
compactification discussed at the beginning of
Section~\ref{sec:down}. Therefore, we obtain certain rigid subspaces
$U_\lambda$ of $P_K$. We denote the inverse system of these by
$\Ulocb{n-1}$ and we have complexes
$\Omega^\bullet(\Ulocb{n-1};\bbox^{n-1})$ and 
$F^\bullet \Omega^\bullet(\Ulocb{n-1};\bbox^{n-1})$ as in
Definition~\ref{mulcompl}. A map $\phi$ whose reduction is a Frobenius
endomorphism and which is compatible with all boundaries is given by
raising to $q$'th power for a sufficiently large $q$. One checks
easily that the sign convention 
for cones \eqref{sign-cone} is such that it commutes with taking the
complex computing multi relative cohomology. From
Lemma~\ref{explicitrgrig} we therefore have a
canonical quasi--isomorphism:
\begin{equation*}
  \Cone\left(F^j \Omega^\bullet(\Ulocb{n-1};\bbox^{n-1})\xrightarrow{1-
 \phi^\ast/q^j}\Omega^\bullet(\Ulocb{n-1};\bbox^{n-1})\right)
 \xrightarrow{\sim} \rgtms(\xlocb B n-1 ,j)
\end{equation*}
From degree considerations it is very easy to see that
\[
F^n \Omega^n(\Ulocb{n-1};\bbox^{n-1})=(\Omega^n(\Ulocb{n-1}),[1,\ldots
,n-1]\rightarrow [1,\ldots ,n-1]) \;.
\]
We can identify this space with
$\Omega^n(\Ulocb{n-1})$.
On the other hand, $F^n \Omega^{n-1}(\Ulocb{n-1};\bbox^{n-1})=0$.
Thus we obtain (compare \eqref{omegeta}) the following expression, with the
identification made above.
\begin{multline}\label{relrep}
    \htms^n(\xlocb B n-1 , n)\\ = \! \frac{\{
      (\omega,\varepsilon):\; \omega\in
      \Omega^n(\Ulocb{n-1}),\; \varepsilon\in
      \Omega^{n-1}(\Ulocb{n-1};\bbox^{n-1}),\; \dd\omega =0,\; \dd\varepsilon =
      (1-\phi^\ast/q^n)\omega
      \}}{ \{ (0,\dd\varepsilon ),\;
      \varepsilon\in \Omega^{n-2}(\Ulocb{n-1};\bbox^{n-1}) \} 
  }\,.
\end{multline}
We take this opportunity to consider two other situations that will be
needed later. In these cases we compute the syntomic cohomology of
$(\xloc n )$ and so there is no $B$ present. The corresponding rigid
space were already considered in previous Sections. We denote by
$\Uloc{n}$ the space $U^{X_{\kappa,\textup{loc}}^n}$. We then have
\begin{multline}\label{sixonehalf}
    \htms^n(\xlocb R n , n)\\ =\frac{\{
      (\omega,\varepsilon):\; \omega\in
      \Omega^n(\Uloc{n}),\; \varepsilon\in
      \Omega^{n-1}(\Uloc{n};\bbox^{n}),\; \dd\omega =0,\; \dd\varepsilon =
      (1-\phi^\ast/q^n)\omega
      \}}{ \{ (0,\dd\varepsilon ),\;
      \varepsilon\in \Omega^{n-2}(\Uloc{n};\bbox^{n}) \} 
  }\;.
\end{multline}
Note that this abuses
the notation somewhat since the differential of $\omega$ is its
differential as a relative form. Also, since there are no $n+1$
relative forms on $(\Uloc{n};\bbox^n)$ we have
\begin{equation}\label{sixthreequot}
    \htms^{n+1}(\xlocb R n , n+1) =\frac{\{
      (0,\varepsilon):\; \varepsilon\in
      \Omega^{n}(\Uloc{n};\bbox^{n}),\; \dd\varepsilon = 0
      \}}{ \{ (0,\dd\varepsilon ),\;
      \varepsilon\in \Omega^{n-1}(\Uloc{n};\bbox^{n}) \} 
  }\;.
\end{equation}
This is simply the $n$th rigid cohomology of $(\xlocb {\kappa} n )$
but note the twisted identification that we have by 
Definition~\ref{normalization}. Also note that in some of the
computations we will be using an altogether different model of this
syntomic cohomology group.

Let
\begin{equation}\label{omegan}
  \omega_n:=\dlog(1-S)\w \dlog\frac{t_1-S}{{t_1-1}}
  \w\dots\w\dlog\frac{t_{n-1}-S}{{t_{n-1}-1}}\; .
\end{equation}
Our main result in this
Section gives the following partial data about the regulator of $\symb S n $.
\begin{proposition}\label{regspecial}
  The regulator of $\symb S n $ in $\htms^n(\xlocb B n-1 , n)$ is
  given, in the representation \eqref{relrep}, by $ (\omega_n,
  \varepsilon_n)$, with some $\varepsilon_n\in
  \Omega^{n-1}(\Ulocb{n-1};\bbox^{n-1})$.
\end{proposition}

Forgetting the relativity gives a map
$ \Krloc n n B n-1 \to \K n n {X_{B,\loc}^{n-1} } $. Let us denote
the image of $ \symb S n $ by $ (S)_n $. The corresponding map in
syntomic cohomology, $\htms^n(\xlocb B n-1 , n) \to
\htms^n(X_{B,\loc}^{n-1}, n)$ simply takes the pair
$(\omega,\varepsilon)$ of \eqref{relrep} to
$(\omega,\varepsilon^\prime)$ in the representation \eqref{omegeta}, where
$\varepsilon^\prime$ is the component of $\varepsilon$ corresponding
to the index function $[1,\ldots ,n-1]\rightarrow [1,\ldots
,n-1]$. Thus, our Proposition follows immediately from the following
Proposition.
\begin{proposition}
  The regulator of $(S)_n$ in $\htms^n(X_{B,\loc}^{n-1} , n)$ is
  given, in the representation \eqref{omegeta}, by $ (\omega_n,
  \varepsilon_n^\prime)$, with some $\varepsilon_n^\prime\in
\Omega^{n-1}(\Ulocb{n-1})$.
\end{proposition}
This last Proposition, again follows easily,
using the formulas for the regulator map for functions, and the
cupproduct in syntomic cohomology given by Lemma~\ref{reg-on-func} and
\eqref{cup-product} respectively, from
the following purely K-theoretic result.
\begin{proposition}
  We have
  \[
  (S)_n =
  (1-S) \cup \frac{t_1-S}{t_1-1}\cup\dots
  \cup\frac{t_{n-1}-S}{t_{n-1}-1}\;.
\]
\end{proposition}

\begin{proof}
Forgetting the relativity
is compatible with the construction of the spectral sequence used
in \eqref{messy-ss}, so in the map
$$
\K n n {\xlocb B n-1 } \to \coprod_{i=1,\dots,n-1} \K n-1 n-1 {\xlocb B n-2 } _{|t_i=S}
$$
the element
$ (S)_n $ will be mapped under the differential in the spectral
sequence
to $ \sum_{i=1}^{n-1} (-1)^i (S)_{n-1|t_i=S} $.  From
this we can determine $ (S)_n $  very easily by induction using
Lemma~\ref{injects-without} below.

\begin{lemma}
\label{injects-without}
For $ m>n>0 $, the map
\[
\K m m {\xl B n } \to \coprod_{i=1,\dots,n-1} \K m-1 m-1 {\xl B n-1 } _{|t_i=S}
\]
is injective.
\end{lemma}

\begin{proof}
Induction on $ n $.  For $ n=1 $, this is clear from the localization
sequence
\[
\dots\to
\K m m X_B \to \K m m {\xl B {} } \to \K m-1 m-1 B \to\dots
\]
as $ \K m m X_B \iso \K m m B = 0 $ if $ m >1 $.
For the induction step, consider the commutative diagram
\[
\xymatrix{
\K m m {\xl B n-1 \times_B X_B}  \ar[r]^-{\phi_1}\ar[d]_{\psi_1} &
         \Coprod_{i=1,\dots,n-2} \K m-1 m-1 {\xl B n-2 \times_B X_B } _{|t_i=S} \ar[d]^{\psi_3}
\\
\K m m {\xl B n-1 \times_B \xl B {} } \ar[r]^-{\phi_2}\ar[d]_{\psi_2} &
         \Coprod_{i=1,\dots,n-1} \K m-1 m-1 {\xl B n-2 \times_B \xl B {} } _{|t_i=S}
\\
\K m m {\xl B n-1 }
}
\]
Here the first vertical column is part of an exact localization
sequence.  $ \psi_3 $ is injective because
$ \K m-1 m-1 {\xl B n-2 \times_B X_B } $ is isomorphic to
$ \K m-1 m-1 {\xl B n-2 } $ under pullback from the base, and
we can restrict the image in $ \K m-1 m-1 {\xl B n-2 \times_B \xl B {} } $
to $ t_{n-1}=0 $ in order to find the element back.  $ \phi_1 $
is injective because again using pullback from the base this
reduces to the case $ n-1 $, where it is true by induction.
In particular, if $ \phi_2(\alpha) = 0 $ for some $ \alpha $,
$ \psi_2(\alpha) = 0 $, and $ \alpha =\psi_1(\beta) $ for some
$ \beta $.  Then $ \psi_3(\phi_1(\beta)) = \phi_2(\alpha) =0 $,
which implies $ \beta = 0 $ as both $ \phi_1 $ and $ \psi_3 $
are injective.  Therefore $ \alpha=0 $.
\end{proof}

By Lemma~\ref{injects-without}, $ (S)_n $ is determined by its
image under the boundary.  As $ (S)_1 = (1-S) $
and $ (S)_2 $ has boundary $ - [S]_1 = - (S)_1 = (1-S)^{-1} $,
one obtains immediately
by induction the required formula,
\[
  (S)_n =
  (1-S) \cup \frac{t_1-S}{t_1-1}\cup\dots
  \cup\frac{t_{n-1}-S}{t_{n-1}-1}\;.
\]
(We use normalizations so that the $ K $--theory acts on the
right in localization sequences.)
\end{proof}

To end this Section, we give the following Lemma.

\begin{lemma}
  Let  $ F(t) $ be an element of $ \I(R) = \K 1 1 {\xlocb R 1 }
  $. Then its regulator in $ \htms^1(\xlocb R 1 ,1) $ is given, in the
  representation \eqref{sixonehalf}, by
  \begin{equation}
    \label{I-reg}
    (\dlog F(t) , \log(F_0(t))/q )
  \end{equation}
\end{lemma}

\begin{proof}
 Note that for $n=1$, \eqref{sixonehalf}
reduces to
\begin{equation}\label{newsixonehalf}
\begin{split}
&  \htms^1(\xlocb R 1 ,1)= \\
& \qquad \{
      (\omega,\varepsilon):\; \omega\in
      \Omega^1(\Uloc{1}),\; \varepsilon\in
      \Omega^0(\Uloc{1}),\; \dd\omega =0,\; \dd\varepsilon =
      (1-\phi^\ast/q)\omega \}\;.
\end{split}
\end{equation}
In this way of writing it looks exactly the same as
$\hms^1(X_{R,\loc}, 1)$. As remarked after \eqref{sixonehalf} this is
slightly misleading since the differentials are different and take
relativity into account. Here This means that
the map $\hms^1(X_{R,\loc};\bbox,1) \to
\hms^1(X_{R,\loc}, 1)$ given simply by
$(\omega,\varepsilon)\mapsto (\omega,\varepsilon)$ embeds
$\hms^1(X_{R,\loc};\bbox,1)$ as the subspace of pairs
$(\omega,\varepsilon)$ where $\varepsilon$ vanishes at $0$ and $\infty$.
Thus, our Lemma is an immediate consequence of
Lemma~\ref{reg-on-func}.
\end{proof}
\section{End of the proof}
\label{sec:end}

We denote the composed map
\begin{equation*}
  \Krloc n n R n-1 \xrightarrow{\reg}
  \htms^{n}(\xlocb R {n-1} , n)\xrightarrow{\sim} 
  \hr^{n-1}(\xlocb {\kappa} n-1 /K) 
  \xrightarrow{\widetilde{\Pi}} K
\end{equation*}
by $\rreg$. Here, the isomorphism is normalized according to
Definition~\ref{normalization} and the map $\widetilde{\Pi}$ is
defined immediately following Corollary~\ref{cor:pi-col}.

\begin{proposition}\label{regisrreg}
  We have the following commutative diagram
  \begin{equation*}
    \xymatrix{
      {\Krloc n n R n-1 }\ar[d]_{\rreg} & {\Krel n n R n-1 }\ar[l]
      \ar[r]^{\sim} & {\K 2n-1 n R }\ar[d]^{\reg}\\
      K \ar@{=}[rr] & & K
      }
  \end{equation*}
\end{proposition}

\begin{proof}
The vertical maps factor through the regulator maps. By the
functoriality of the normalization of Definition~\ref{normalization},
the commutativity of the diagram will follows if we show the
commutativity of the diagram
\begin{equation*}
  \xymatrix{
    {\hr^{n-1}(\xlocb {\kappa} n-1 /K)}\ar[d]_{\widetilde{\Pi}} &
    {\hr^{n-1}(\xb {\kappa} n-1 /K)}\ar[l] 
    \ar[r]^{\sim} & \hr^0(\Spec(\kappa)/K) \ar[d]\\
    K \ar@{=}[rr] & & K
    }
\end{equation*}
But as explained in Proposition~\ref{newdownsign} the composed map
$\hr^{n-1}(\xb {\kappa} n-1 /K)\to K$ is simply the map $\Pi$ and
therefore the commutativity follows from Corollary~\ref{cor:pi-col}.
\end{proof}

\begin{proposition}\label{rregfactors}
The composition
\[
\Symb n R \subset \Krloc n n R n-1 \xrightarrow{\rreg} K
\]
factors through the quotient
$  \Symb n R / \I \tcup\,\, \Symb n-1 R = M_n(R) $.
\end{proposition}

\begin{proof}
In fact, we can show that $\rreg$ vanishes on 
$\I \tcup\,\, \Krloc {n-1} {n-1} R n-2 $. This will follow
by symmetry for all possible products involved in $\tcup$ if we show
that the composition
\begin{multline*}
  \htms^1(\xlocb R 1 ;1) \times
  \htms^{n-1}(\xlocb R {n-2} , n-1)\\ \xrightarrow{\cup}
  \htms^{n}(\xlocb R {n-1} , n)\xrightarrow{\sim} 
  \hr^{n-1}(\xlocb {\kappa} n-1 /K) 
  \xrightarrow{\widetilde{\Pi}} K
\end{multline*}
vanishes on pairs where the first coordinate is used for $\I$.
Let $F(t)$ be in  $\I$. By \eqref{I-reg}, its regulator is given, in
the representation
\eqref{sixonehalf} by $(\dlog F(t),?)$, where the first coordinate
belongs to $\Omega^1(\Uloc{1})$ and the precise form of the second coordinate does not
matter as we will
see in a second.  On the other hand, elements of
$\htms^{n-1}(\xlocb R {n-2} , n-1)$ are, by \eqref{sixthreequot},
always of the
form $(0,\delta)$, with $\delta $ in $ \Omega^{n-2}(\Uloc{n-2};\bbox^{n-2})$.
Choosing $\gamma=0$ in \eqref{cup-product}, we see that 
$ (\dlog F(t), ?) \cup (0,\delta) $, will be of the form
$(0,\dlog F(t)\wedge \delta)$ where $\wedge$ here means the product
in complexes of relative differential forms as defined in
Remark~\ref{cupcompatible} in
Appendix~\ref{sec:appendix}. By Definition~\ref{normalization}
the image of $ (\dlog F(t), ?) \cup (0,\delta) $ in 
$\hr^{n-1}(\xlocb{\kappa} n-1 /K)$ is the
inverse of the operator $1-\varphi^\ast/q^{n-1}$ applied to
the cohomology class $[\dlog F(t)\wedge \delta]$. Since the operator
$\widetilde{\Pi}$ is Frobenius equivariant we see that applied to this
image it gives $(1-q^{1-n})  \widetilde{\Pi}([\dlog F(t)\wedge \delta])$ and
so our goal is to show that
\begin{equation*}
  \widetilde{\Pi}([\dlog F(t)\wedge \delta])
  =\widetilde{\Pi}([\dlog F(t)]\cup [\delta])
\end{equation*}
vanishes, where the cup product on the right is a cup product in
multi relative rigid cohomology. By Lemma~\ref{pairi} it equals
$\widetilde{\Pi}([\dlog F(t)])
\cdot \widetilde{\Pi}([\delta])$ and the result follows since
\begin{equation*}
 \widetilde{\Pi}([\dlog F(t)])= \Pi_\Cole(\dlog
 F(t))=\log(F(\infty))-\log(F(0))=0\;.
\end{equation*}
\end{proof}

We continue to denote the induced map by $\rreg$,
\begin{equation}
  \label{eq:composedreg}
  \rreg:  M_n(R) \to K\;.
\end{equation}

Recall that in Definition~\ref{symbols}  we defined for any $z$ in $ R^\flat$
a symbol $\symb z n $ in $ \Krloc n n R n-1 $ by pullback of the
universal symbol $\symb S n  $ in $ \Krloc n n G n-1 $ along the map sending $S$ to $z$,
where $G=\Spec(\Z[S,S^{-1},(1-S)^{-1}])$.
We now compute $\rreg(\symb z n )$. We begin this by
exploring some auxiliary functions.
\begin{definition}\label{def:special-func}
We define a sequence of functions $f_k(z,S)$ inductively as follows:
\begin{equation}
  f_0(z,S) = \frac{S}{1-S},\quad f_{k+1}(z,S)= \int_z^S f_k(z,t)\dlog
  t\;. 
\end{equation}
\end{definition}
At some later point we will want to interpret the integral as a Coleman
integral but at this point it is enough to define it when $z$ and $S$ belong to
the same residue disc, in which case it is well defined without making any
Frobenius equivariance assumptions. It is immediately noticed that $f_k(z,S)$
vanishes to order $k$ at $z=S$.
\begin{lemma}
\label{fnzS}
  We have
  \begin{equation*}
    f_n(z,S)=\Li_n(S) - \sum_{k=0}^{n-1} \frac{1}{k!}
 (\log S-\log z)^k \Li_{n-k}(z)\;.
  \end{equation*}
\end{lemma}

\begin{proof}
The proof is by induction on $n$. For $n=1$ it is immediately verified that
\begin{equation*}
  f_1(z,S) = \log(1-z) - \log(1-S)= \Li_1(S) - \Li_1(z).
\end{equation*}
Suppose that the statement is true for $n$. Then for $n+1$ we get
\begin{align*}
  f_{n+1}(z,S)&= \int_z^S \left(
    \Li_n(t) - \sum_{k=0}^{n-1} \frac{1}{k!}
 (\log t-\log z)^k \Li_{n-k}(z)
  \right) \frac{dt}{t}\\
  &=  \left. \left(
    \Li_{n+1}(t) - \sum_{k=0}^{n-1} \frac{1}{(k+1)!}
    (\log t-\log z)^{k+1} \Li_{n-k}(z)
  \right) \right|_z^S\\
  &= \Li_{n+1}(S) - \Li_{n+1}(z)
  - \sum_{k=0}^{n-1} \frac{1}{(k+1)!}
 (\log S-\log z)^{k+1} \Li_{n-k}(z)\\
    &=\Li_{n+1}(S) - \sum_{k=0}^{n} \frac{1}{k!}
 (\log S-\log z)^k \Li_{n+1-k}(z)\;.
\end{align*}
\end{proof}

\begin{proposition}\label{closelis}
  Let $z_1,z_2\in R^\flat$ belong to the same residue disc. Then we
  have
  \begin{equation*}
    \rreg([z_1]_n)-\rreg([z_2]_n) = (-1)^{n} (n-1)! ( \L_{n}(z_1)-
    \L_{n}(z_2) )\;.
  \end{equation*}
\end{proposition}

\begin{proof}
For any $z\in R^\flat$ we may factor the map $\Spec(R)\to G$ defined
by sending $S$ to $z$ via the map $\Spec(R)\to B=G\otimes_{\Z}R$
defined in the same way. By functoriality of the regulator map it
follows that $\reg(\symb z n )$ equals $i_z^\ast \reg(\symb S n )$,
where $i_z:(\xlocb R n ) \to (\xlocb B n ) $ is the embedding in
$(\xlocb B n )$ of the fiber at $z$. Thus we are in position to apply
Corollary~\ref{pullback1}, but in the
relative case, which, as mention after its statement, also applies.
To carry out the computation we also shift the index from $n$ to $n+1$
as the computation seems to come out a bit cleaner this way. 

We begin with $\htms^{n+1}(\xlocb B n , n+1)$ in the representation
\eqref{relrep} (with $n$ shifted to $n+1$). In there we have the regulator
of $[S]_{n+1}$, given, according to Proposition~\ref{regspecial}, by the
pair $(\omega,\epsilon)$, where 
\begin{equation*}
  \omega=\omega_{n+1} = \dlog (1-S) \wedge \dlog
  \frac{t_{1}-S}{t_{1}-1}\wedge \cdots \wedge \dlog \frac{t_{n}-S}{t_{n}-1}
\end{equation*}
is the form defined in \eqref{omegan}, while $\epsilon$ is unknown.
Note that $\omega$ should really be thought of as the
relative form $\widetilde{\omega }=(\omega ,[1,\ldots ,n]\rightarrow
[1,\ldots ,n])\in \Omega^{n+1}(\Ulocb{n};\bbox^n)$.

We have the projection $\pi: (\xlocb B n )\to B$, which we can
compactify  to $\bar{\pi}:(\PP_B^1)^n\to B$, where the power is taken
over $B$. By assumption,
$z_1$ and $z_2$ belong to the same residue disc, which we call $D$.

The recipe for computing 
\begin{equation*}
  \reg(\symb {z_1} {n+1} ) - \reg(\symb {z_2} {n+1} )=
  i_{z_1}^\ast \reg(\symb S {n+1} ) - i_{z_1}^\ast \reg(\symb S {n+1} )\;,
\end{equation*}
according to Corollary~\ref{pullback1}, calls for computing, for $z=z_1$ and
$z=z_2$, a form $\theta_{z}$ such that
\begin{equation}
  \label{eq:thetacond}
  \theta_{z}\in F_J^{n+1} \Omega^n(\Ulocb{n}\cap
  \bar{\pi}^{-1}(D);\bbox^n)\;\text{and}\; \dd\theta_z
  =\widetilde{\omega}|_{\bar{\pi}^{-1}(D)}\;,
\end{equation}
where $J$ is the ideal defining $\bar{\pi}^{-1}(z)$.
Such a form is given in the following Lemma.

\begin{lemma}
  Let
\begin{equation*}
  \theta _{z}=-\sum _{k=0}^{n}(-1)^k k!\sum _{\deg
    h=k}((-1)^{\sum (h(i)+i)} f_{k+1}(z,S)\bigwedge _{i=1}^{n-k}\dlog
  \frac{t_{h(i)}-S}{t_{h(i)}-1},h).
\end{equation*}
Then $\theta _{z}$ satisfies \eqref{eq:thetacond}.  Here, $f_{k}$ is the
function introduced in Definition~\ref{def:special-func} and the form
indicated only for the
component in which all the constant coordinates are $0$, otherwise
the form is $0$.
\end{lemma}

\begin{proof}
Recall that the condition for being in $F_J^{n+1}$ is that the sum of
the degree of the form and its order of vanishing at $S=z$ is
$n+1$. and this is clear for $\theta_z$. Now we prove  that the
differential is correct. We will show that $\dd
(-\theta_z)=-\widetilde{\omega }$. Consider first the
differential of a single term in $-\theta_z$.
\begin{equation*}
  (-1)^{k+\sum (h(i)+i)} k! \, (f_{k+1}(z,S)\bigwedge _{i=1}^{n-k}\dlog
  \frac{t_{h(i)}-S}{t_{h(i)}-1},h)\;,
\end{equation*}
with $h$ of
degree $k$. Using \eqref{multidif} the differential is
\begin{equation*}
  (-1)^{k+\sum (h(i)+i)} k! \, (f_{k}(z,S)\dlog (S)\wedge \bigwedge
  _{i=1}^{n-k}\dlog \frac{t_{h(i)}-S}{t_{h(i)}-1},h)
\end{equation*}
minus (except when $n=k$) a sum of
terms obtained by restricting one of the coordinates $t_{h(j)}$
to be $0$ (when we restrict any coordinate to $\infty$ we
get $0$ and there is no need to keep track of that). In the
wedge product at $j$ we get $\dlog (S)$ and moving it to the
front gives a sign of $(-1)^{j-1}$.  This form is then
associated with a function $g$ for which $
\chi (h,g)=(-1)^{h(j)+j}$. In addition there is an overall sign of
$(-1)^{n-k}$ on the entire sum. Thus, the sign on the component with
the function $g$ obtained from $h$ by deleting $h(j)$ has a sign of
$(-1)$ to the power
\begin{multline*}
  n-k+k+ \sum (h(i)+i) + h(j) -1\\ \equiv n+1+ \sum g(i)
  +\sum_{i=1}^{n-k} i \equiv k+1+ \sum (g(i)+i) \pmod{2}\;.
\end{multline*}
Thus we find
\begin{align*}
  &\phantom{=}\dd (-1)^{k+\sum (h(i)+i)} k! \, (f_{k+1}(z,S)
  \bigwedge_{i=1}^{n-k}\dlog \frac{t_{h(i)}-S}{t_{h(i)}-1},h) \\
  &= (-1)^{k+\sum (h(i)+i)} k! \, (f_{k}(z,S)\dlog (S)\wedge
  \bigwedge_{i=1}^{n-k}\dlog \frac{t_{h(i)}-S}{t_{h(i)}-1},h)\\
  &\phantom{=}\,\, - \sum _{\chi (h,g)\neq 0} (-1)^{k+1+\sum (g(i)+i)} k!
  \, (f_{k+1}(z,S) 
  \dlog (S)\wedge \bigwedge_{i=1}^{n-k-1}\dlog
  \frac{t_{g(i)}-S}{t_{g(i)}-1},g)\;.
\end{align*}

Now we consider the coefficient in $\dd(-\theta_z)$ in the $g$
component when $\deg g=m$. If $m>0$ then it gets contributions
from both lines in the right hand side of the last equation.
The contributions from the second line correspond to $k=m-1$. There
are exactly $m$ different $h$'s that would give $g$ and the
contributions are identical. Thus it is visibly seen that the
contributions from the second line cancel the ones from the first line.
The only term that survives is the one with $m=0$. Here there is
only a contribution from the first line. So we find
\begin{align*}
  \dd(-\theta_z)&=\,\, (f_{0}(z,S)\dlog (S)\wedge
  \bigwedge_{i=1}^{n}\dlog \frac{t_{i}-S}{t_{i}-1},\id )\\
  &=-(\dlog (1-S)\wedge \bigwedge_{i=1}^{n}\dlog
  \frac{t_{i}-S}{t_{i}-1},\id )\\
  &=-\widetilde{\omega }\;.
\end{align*}
\end{proof}

Now
$ i_{z_1}^\ast \reg(\symb S n ) - i_{z_2}^\ast \reg(\symb S n ) =
\theta_{z_2}-\theta_{z_1}$ 
restricted to the fiber above any $ S $ in
$D(K)$. On this difference we need to apply $\Pi_\Cole$. We will in
fact compute $\Pi_\Cole(\theta_z)$ for any $z$. We integrate with
respect to the
$t$'s keeping $S$ fixed.  We see
that all the terms we need to successively integrate are products of
$\dlog$'s, which the integration process converts into logs.
The extra sign coming from the formula for $\Pi$ is
$(-1)^{\sum (h(i)+i)}$ because we are always in the component where all the
fixed coordinates are $0$. So using Lemma~\ref{fnzS}
and the definition of $ \L_{n+1}(S) $ as in~\eqref{Ldef}, we find  $ \Pi_\Cole(\theta_z) $ equals
\begin{equation*}
\begin{split}
 & \phantom{=} \,\, - \sum _{k=0}^{n}(-1)^{k}k!\cdot
  \binom{n}{k}f_{k+1}(z,S)\log ^{n-k}(S)\\
  & = -n!\sum _{k=0}^{n}(-1)^{k}\frac{1}{(n-k)!}f_{k+1}(z,S)\log ^{n-k}(S)\\
  & =  - n! \, \sum_{k=0}^{n}(-1)^{k} \frac{1}{(n-k)!} \Li_{k+1}(S) \log ^{n-k}(S)\\
  & \qquad + n! \sum_{k=0}^n (-1)^k \frac{1}{(n-k)!}
    \left[ \sum_{l=0}^k \frac{1}{l!} (\log(S)-\log(z))^l\Li_{k+1-l}(z)\right]\log^{n-k}(S)\\
  & =  (-1)^{n-1}n! \, \L_{n+1}(S) \\
   & \qquad + n!\sum_{k=0}^n(-1)^k\frac{1}{(n-k)!}
        \left[\sum_{r=0}^k\frac{1}{(k-r)!}(\log(S)-\log(z))^{k-r}\Li_{r+1}(z)\right]\log^{n-k}(S)\\
  & =  (-1)^{n-1}n! \, \L_{n+1}(S) \\
  &  \qquad + n!\sum_{r=0}^n\sum_{k=r}^n(-1)^n
  \frac{1}{(n-r)!}\binom{n-r}{k-r} (\log(S)-\log(z))^{k-r}(-\log(S))^{n-k}\Li_{r+1}(z)\\
  & =  (-1)^{n-1}n! \, \L_{n+1}(S) \\
  & \qquad + (-1)^n  n!\sum_{r=0}^n
        \frac{1}{(n-r)!}(\log(S)-\log(z)-\log(S))^{n-r}\Li_{r+1}(z)\\
  & = (-1)^{n} n! \left( \L_{n+1}(z) - \L_{n+1}(S) \right)\;.
\end{split}
\end{equation*}
Thus we find
  \begin{align*}
    \rreg([z_1]_{n+1})-\rreg([z_2]_{n+1}) &=
 (-1)^{n} n! \left( (\L_{n+1}(z_2) - \L_{n+1}(S))-(\L_{n+1}(z_1) -
   \L_{n+1}(S)) \right)\\&=
  (-1)^{n+1} n!( \L_{n}(z_1)-  \L_{n}(z_2) )\;.
  \end{align*}
\end{proof}

\begin{proposition}\label{rregisl}
  For $ z $ in $ R^\flat $ we have
  \begin{equation}\label{main2}
    \rreg(\symb z n )=(-1)^n (n-1)! \, \L_n(z)\;.
  \end{equation}
\end{proposition}

\begin{proof}
Let $E_n(z)$ be the difference of the two sides of the equation.
By Proposition~\ref{closelis} $E_n(z)$ is constant on each residue disc.
The function $E_n(z)$ satisfies the distribution relation
\begin{equation}\label{distribute}
  \frac{1}{m} \sum_{\zeta^m=1}E_n(\zeta z) = \frac{E_n(z^m)}{m^n}
\end{equation}
for each positive integer $m$.
The left hand side of \eqref{main2} satisfies the relation because
by~\cite[Proposition 6.1]{Jeu95} we have the relation
\begin{equation*}
    \frac{1}{m} \sum_{\zeta^m=1} \symb {\zeta z} n =
    \frac{\symb {z^m} n }{m^n}
\end{equation*}
in $\Krloc n n R n-1 $ (modulo terms involving $ \I $) and we
then apply $\rreg$.
(Note, in loc.\ cit.\ the relation is stated for elements in a
field containing $\Q(\zeta)$.  But the proof of the
statement shows there is a corresponding universal relation over
$\Z[X,X^{-1}]$ which can be pulled back.  Alternatively, it can
be deduced from the relation in $\Q[X,X^{-1}]$ because
$\F_p[X,X^{-1}]$ has no low weight $K$--theory, and the
localization map (in the limit) corresponding to $\Z \to \Q$ will
induce an injection (up to torsion) on the level of symbols, cf. Proposition~\ref{O-F-injective}.)
The right hand side satisfies the relation
because it is true for $\Li_n$
by~\eqref{poly-distribution} and for the remaining terms by a
straightforward standard computation.
Multiplication by a $p^k$--th root of unity preserves the residue
discs, as does raising to $p^k$ power for sufficiently divisible
$k$ (here we need to extend $R$ to include these roots of unity). Therefore,
it is immediately seen that the function $E_n$ must be $0$.
\end{proof}

\begin{remark}
The following comparison with the work of De Jeu is perhaps interesting. In the
complex case one again relies only on the explicit description of the form
$\omega$ to obtain the corresponding formula for the regulator. A similar
constant must be fixed in that computation as well. There however, one relies
on the fact that the final result should be single valued. One then derives the
distribution relation from the corresponding formula for the complex
polylogarithm and Borel's theorem. Here we have used this distribution relation
so our proof relies on the proof in the complex case. It may be interesting to
mention that in the work of Wojtkowiak on functional equations for
polylogarithms~\cite{Woj91} a similar phenomenon occurs: To obtain a
functional equation
for $p$-adic polylogs one starts with a functional
equation for the complex polylog and uses the multivaluedness to show that a
certain ``motivic'' functional equation is satisfied, which then translates
into a $p$-adic functional equation.
\end{remark}

\begin{proposition}\label{regisln}
  The map
  \begin{equation*}
   H^1( \scriptM n ({\O}) ) \xrightarrow{\Phi} \K 2n-1 n {\O} 
   \xrightarrow{\sigma} \K 2n-1 n R \xrightarrow{\reg} K\;,
  \end{equation*}
  where the map $\Phi$ is part of \eqref{eq:kbasicmap}, is induced by
  the map sending the symbol $\symb x n $ to $(-1)^n (n-1)! \,
  \L_n(\sigma(x))$.
\end{proposition}
\begin{proof}
Suppose  $\sum a_i \symb {x_i} n \in H^1( \scriptM n ({\O}) )$. Let
\begin{equation*}
  \alpha=\Phi(\sum a_i \symb {x_i} n )\in \K 2n-1 n {\O} \isom 
  \Krel n n {\O} n-1 \;.
\end{equation*}
By definition, the image of $\alpha$ in $\Krloc n n {\O} n-1 $ equals
$\sum a_i \symb {x_i} n $ modulo $\I \tcup\,\, \Symb n-1 {\O} $.
By functoriality, the image of $\sigma(\alpha)$ in $\Krloc n n R n-1 $ equals 
$\sum a_i \symb {\sigma(x_i)} n $ modulo $\I \tcup\,\, \Symb n-1 R $.
By Proposition~\ref{regisrreg} and Proposition~\ref{rregfactors} we therefore
have  
\begin{align*}
  \reg(\sigma(\alpha))&=\rreg(\sum a_i \symb {\sigma(x_i)} n )\\
  &=\sum a_i \rreg(\symb {\sigma(x_i)} n )\\
  &=\sum a_i (-1)^n (n-1)! \L_n(\sigma(x_i)) \quad \text{by
    Proposition~\ref{rregisl}}\\
\end{align*}
\end{proof}
\begin{proof}[Proof of Theorems~\ref{general-field}~and~\ref{number-field}]
Part (1) of each Theorem was already proved in Section~\ref{sec:kstuff}.
To prove part (2), note that any of the functions $ \Lmod n (z) $
differs from $\L_n(z)$ by a linear combinations of
the functions $ z\mapsto\log^k(z) L_{n-k}(z) $ for $k\ge 1$. Any
function in this combination, when composed with $\sigma$, factors
through the composed differential
\[
M_n(\O) \to M_{n-1}(F) \tensor \FQ \to\dots\to
M_{n-k}(F)\tensor (\FQ)^{\tensor k}
\]
which maps $ \symb x n $ to $ \symb x n-k \tensor x\tensor\dots\tensor x
$.
Therefore the functions $ L_n(\sigma(x)) $ and $ \Lmod n (\sigma(x)) $
(and in fact, also $ \Li_n(\sigma(x)) $),
coincide on $ H^1(\scriptM n (\O)) $. But any function $\Lmod n (z)$
satisfies $ \Lmod n (z)+(-1)^n\Lmod n (1/z) $. Therefore the map
$\symb x n \mapsto (-1)^n (n-1)! \,\Lmod n (\sigma(x)) $ factors
through the quotient
map $ M_n(\O) \to \tM_n(\O) $.  But the composition
\[
H^1(\scriptM n (\O)) \to H^1(\tildescriptM n (\O))\to K
\]
is still given by $\symb x n  \mapsto (-1)^n (n-1)! \,
\L_n(\sigma(x))$. Thus, the Theorem follows from Proposition~\ref{regisln}.
\end{proof}
\begin{proof}[Proof of Theorem~\ref{cyclotomic}]
Let $ F $ be a number
field.  Note that roots of unity will
not be special units in general, so we have to work in the complex
for $ F $ rather than for $ \O $.
Namely, let $ \zeta $ in $ F $ be a root of unity of order $ m > 1 $.
If $ (m,p) = 1 $ $ \zeta $ is a special unit, and we have the result
already.  If $ m = p^s l $ with $ s>1 $ and $ (p,l) = 1 $, write
$ \zeta = \zeta_1 \zeta_2 $ with
$ \zeta_1 $ of order $ p^s $ and $ \zeta_2 $ of order $ l $.
As the reduction of $ \zeta $ is the same as the reduction of
$ \zeta_2 $, we see that $ \zeta $ is special unless $ m = p^s $.
If $ m = p^s $, let $ r > 1 $ be an integer congruent to $ 1 $ modulo
$ p^s $.  Then $ \zeta^r = \zeta $, and from the distribution
relations Proposition~\ref{distribution} and \eqref{poly-distribution}
we find that
$ \symb {\zeta} n = r^{n-1} \sum_{\alpha^r =1} \symb {\zeta\alpha} n $
in $ M_n(F') $ with $ F' = F(\mu_r) $.
According to~\cite[Lemma~\ref{dont-know}]{Bes98a} the modified syntomic
regulator commutes with finite base change. This means in the case we are
considering that if $R'$ is a finite extension of $R$
there is a commutative diagram
\begin{equation*}
  \xymatrix{
     {\K 2n-1 n R } \ar[r]\ar[d]& K\ar[d]\\
     {\K 2n-1 n R' } \ar[r]& K',
  }
\end{equation*}  
where the map $ K \to K' $ is the natural inclusion.
We can therefore do our computations
for the regulator just as well in $ K'=K(\mu_r) $.
As all $ \zeta\alpha $
in the sum are special units except when $ \alpha =1 $, we can
solve for $ \reg(\symb {\zeta} n ) $ as $ r^{n-1} \neq 1 $.
As $ \reg(\symb x n ) = \Lmod n (x) $ if $ x $ is a special unit,
and $ \Lmod n $ satisfies the corresponding distribution relation
by \cite[Proposition 6.1]{Col82}, our result follows for $ m>1 $.
For $ \zeta =1 $ one uses the distribution relation similarly.
\end{proof}
\begin{remark}
Although not needed for the purposes of this paper, we would like to
sketch a somewhat less explicit method of doing the computations of
this Section. This was in fact our original method.
The idea is quite easy to explain: Out of our regulator
computations, we obtained the fact that for certain constants
$\alpha_k$ the function
\begin{equation*}
  P(z,S)=\sum _{k=0}^{n}\alpha _{k}f_{k}(z,S)\log
  ^{n-k}S
\end{equation*}
is a sum of a function of $z$ and a function of $S$ (first line of
final computation). In other words,
its mixed derivatives vanish. We get this relation initially only for
$z$ and $S$ which belong to the same residue disc, even though we know
that using Coleman integration the functions $f_n$ themselves extend
to all $z$ and $S$. What we did was to
write explicitly $f_n(z,S)$ in terms of logarithms and polylogarithms
in $z$ and $S$ separately and then show that $P(z,S)$ can indeed be
written as a sum of a function of $S$ and a function of $z$, the
latter being our sought after regulator function. Note that the separation of
variables now holds for any $z$ and $S$.

The alternative approach is to deduce the ``global'' separation of
variables from the same result for $z$ and $S$ in the same residue
disc using Coleman theory in $2$-variables. The theory developed
in~\cite{Bes99} defines a
notion of a Coleman function in several variables. One then shows that
iterated integrals of the kind used to define $f_n(z,S)$ make it a
Coleman function in both variables (in particular for fixed $S$ it is
a Coleman function in $z$, but note that the notion of being Coleman
in two variables is stronger than the notion of begin Coleman in each
variable separately). Coleman functions form a ring, which shows that
$P(z,S)$ is
also a Coleman function, and so are its mixed derivatives. The theory
then shows that the fact that the mixed derivatives vanish on some
residue disc imply that they vanish identically which in turn implies
a global separation of variables.

Knowing a global separation of variables is very convenient, for if
$P(z,S)=P(z) + g(S)$, and all we need to know is $P(z)$ up to a
constant, then this is supplied by $P(z,S_0)$ for any $S_0$ we take.
In our particular situation, if one substitutes $S_0=-1$, whose $\log$ is
$0$, in the first line of the final computation we find that up to a constant
$P(z)=-(-1)^n n! \cdot f_{n+1}(z,-1)$. So one is then left with
computing this function.
\end{remark}
\appendix
\section{Chern classes in relative cohomology}
\label{sec:appendix}

In this appendix we give the necessary constructions for the
main paper as far as relative K--theory and Chern classes in syntomic
cohomology are concerned.  Most of this material is rather standard
and has to be modified in a rather minor way in order to fit
the current context, so we are sketchy in places.  One thing
that we work out in glorious detail is the description of a complex
that computes the multi relative syntomic cohomology encountered
in Section~\ref{sec:down}.

In~\cite{Bes98a} the first named author described a theory of rigid
syntomic regulators. This is not sufficient in all applications, for
example those  described in the present work, since one often needs to
extend the regulators (= Chern classes) to the relative situation. Our
goal in this appendix is to explain how the construction extends to
more general ``spaces'', yielding in particular Chern classes from
relative $K$-theory to relative syntomic cohomology. The construction
is completely formal, relying primarily on~\cite{Gil-Sou99} with some
input from~\cite{Jeu95}.

We consider the category $\bcat$ of Noetherian finite dimensional
schemes over $\vv$
and the topos $T$ of sheaves on it with respect to the Zariski
topology. Following~\cite{Gil-Sou99} we call a simplicial object of
$T$ a space.
If $ X\in \bcat $, $ Y \mapsto \hom(Y,X) $ gives us an element of $ T $ ,
which we still denote by $ X $.  With $ X_+ $ we will denote
the space consisting of the disjoint union of the constant sheaf
$ X $ and the constant simplicial basepoint $ * $.  

Let $\bcat'\subset \bcat$ be the subcategory of schemes which are in
addition smooth and separated over $\vv$, again equipped with the
Zariski topology.
Following~\cite[Section~2.2]{Jeu95} we make the following definition.
\begin{definition}\label{pss-def}
A space $X_\bullet$ is called a smooth separated pointed simplicial
scheme if it is represented in each degree by a scheme in $\bcat'$
(interpreted as sheaf) together with
a disjoint base point and if it is furthermore degenerate above a finite
simplicial degree. In the sequel we shall refer to these spaces as
pointed simplicial schemes for short, unless stated otherwise.
\end{definition}
Clearly, if $ X\in \bcat' $, then $ X_+ $ is a pointed simplicial
scheme.  
According to~\cite[Lemma~2.1]{Jeu95} a pointed simplicial scheme is $K$-coherent in the
sense of~\cite[Section 3.1, Definition 1]{Gil-Sou99} (this notion is
also defined in~\cite[top of page 201]{Jeu95}).

When $X_\bullet$ and $Y_\bullet$ are spaces, $[X_\bullet,Y_\bullet]$ denotes the set of morphisms
between $X_\bullet$ and $Y_\bullet$ in the homotopy category of spaces constructed
in~\cite{Gil-Sou99}.
The point is that there is an element $ K $ in $ T $ such that
if $ X $ is a smooth finite type scheme over $ \vv $, then
$ [S^m \w X_+, K]  \iso K_m(X) $ (loc.\ cit.~Proposition~5 of
Section~3.2).
Gillet--Soul\'e therefore make the following Definition.

\begin{definition}\label{kdef}
If $ X_\bullet $ is a space, we define 
$$
K_m(X_\bullet) = H^{-m}(X_\bullet, K) = [S^m \w X_\bullet, K]
$$
for $ m \geq 0  $.
\end{definition}

In the body of the paper we have systematically used the $ K $--theoretic
notation for this, but in this Appendix we shall also use the
notation $ H^{-m}(X_\bullet, K) $.

In \cite{Gil-Sou99}, it is then shown that if $ X_\bullet $ is
$ K $--coherent, then $ H^{-m}(X_\bullet) $ has a special $ \lambda $--ring
structure.  In particular, 
when $ X_\bullet $ is $ K $--coherent of cohomological dimension
at most $ d $, then according to loc.\ cit.\
Proposition~8 of Section~4.4, 
$ K_m(X_\bullet) \tensor_\Z \Q = \bigoplus_{i=\a}^{m+d} K_m^{(i)}(X_\bullet) $,
with $ K_m^{(i)}(X_\bullet) = H^{-m}(X_\bullet, K)^{(i)} $
 the
$\Q$ sub vector space of $ K_m( X_\bullet)_\Q = H^{-m}(X_\bullet,K)_\Q$
of elements $x$
such that $ \psi^k(x) = k^i x $ for all $ k\geq 2 $, and $ \a = 2 $
if $ m \geq 2 $, $ \a = 1 $ if $ m=1 $ and $ \a=0 $ if $ m=0 $.
This will certainly apply to a pointed simplicial space
$ X_\bullet $ which is degenerate above simplicial degree $N$ and
where the maximum relative dimension of scheme components over
$\vv$ is $M$, with $d=M+N+1$, cf.\ loc.\ cit.\ Lemma~1.2.2.2 or~3.2.4.

For two spaces $ X_\bullet $ and $ Y_\bullet $ we get a map 
$ K_m(X_\bullet) \times K_n(Y_\bullet) \to K_{m+n}(X_\bullet \w Y_\bullet)$
from the composition of 
$$
S^{m+n} \w X_\bullet \w Y_\bullet \to S^m \w X_\bullet \w S^n \w Y_\bullet
\to K \w K \to K
$$
because $ S^m \w S^n = S^{m+n} $ and $ K $ comes equipped with
a map $ K \w K \to K $.
Under this product, we get
$ K_m^{(i)}(X_\bullet) \times K_n^{(j)}(X_\bullet) \to K_{m+n}^{(i+j)}(X_\bullet \w Y_\bullet) $
(loc.\ cit.~top of page~48).

If $ X_\bullet $ and $ Y_\bullet $ are pointed spaces with a base
pointed preserving map $ Y_\bullet \to X_\bullet $, then the reduced
mapping cone $ C_\bullet = C(Y_\bullet, X_\bullet) $,
whose definition will be recalled later, see \eqref{mapcone},
is also a
pointed space, and one gets an exact sequence
$$
\dots \to
K_{m+1}(X_\bullet) \to K_{m+1}(Y_\bullet) \to K_m(C_\bullet) \to K_m(X_\bullet)
\to \cdots
\,.
$$

The most important applications of these are if $ Y_\bullet $
is a pointed closed simplicial subscheme of $ X_\bullet $
(i.e., the map corresponds to a closed embedding of schemes on
all the scheme components of $ X_\bullet $), in which case one
gets the $ K $--theory of $ X_\bullet $ relative to $ Y_\bullet $:
$ K_m(C_\bullet) = K_m(X_\bullet; Y_\bullet) $. Iterating this
one gets multi relative $ K $--groups as in the body of the paper.
E.g., if $ Y_{1,\bullet} $ and $ Y_{2\bullet} $ are closed simplicial subschemes of
$ X_\bullet $, and $ Y_{12\bullet} = Y_{1\bullet} \bigcap Y_{2\bullet} $,
with $ C_{1\bullet} = C(Y_{1\bullet} , X_\bullet) $,
$ C_{2\bullet} = C(Y_{12\bullet} ,  Y_{2\bullet}) $ and
$ C_{3\bullet} = C(C_{2\bullet} , C_{1\bullet}) $ 
we get an exact sequence
$$
\cdots\!\to\!
K_{m+1}(X_\bullet ; Y_{1\bullet}) \to K_{m+1}(Y_{2\bullet}; Y_{12\bullet})
\to K_m(X; Y_{1\bullet}; Y_{2\bullet}) \to K_m(X_\bullet ; Y_{1\bullet}) 
\to\!\cdots
$$
where we write $ K_m(X; Y_{1\bullet}; Y_{2\bullet}) $ for
$ K_m(C_{3\bullet}) $.

The other application is $ K $--theory with support, in which
case $ Y_\bullet $ is an open pointed subscheme of $ X_\bullet $.
Let $ Z_\bullet $ be the closed pointed simplicial scheme complement of
$ Y_\bullet $ in $X_\bullet $ (i.e., the closed
complement in every scheme component, together with the base
point $ * $), and  assume all scheme components of $ Z_\bullet $
are regular.
If also conditions (TC1) and (TC2) of \cite[page~202]{Jeu95}
hold for the embeddings in $Z_\bullet \to X_\bullet $,
then $ K_m(Z_\bullet) \iso K_m(C(Y_\bullet \to X_\bullet)) $.
The sequence then becomes a localization sequence
$$
\dots\to
K_{m+1}(X_\bullet) \to K_{m+1}(Y_\bullet) \to K_m(Z_\bullet) \to K_m(X_\bullet)
\to\cdots
\,.
$$

Under very restrictive hypothesis where $ Z_\bullet $ is of codimension
$ d $ in all scheme components
(see \cite[Proposition~2.3]{Jeu95})
one can prove a Gysin exact sequence
$$
\dots\to
K_{m+1}^{(i+d)}(X_\bullet) \to K_{m+1}^{(i+d)}(Y_\bullet) \to K_m^{(i)}(Z_\bullet) \to K_m^{(i+d)}(X_\bullet)
\to\cdots
\,.
$$

In order to be able to define regulators with values in the cohomology
of a complex of Abelian groups, we briefly review how this is
put into the context of spaces.

When $ A_\bullet $ is a homological chain complex of Abelian
objects in $T$, $X_\bullet$ is a space, and $ n \geq 0$,
we write $H^{-n}(X_\bullet,A_\bullet):=[S^n \w X_\bullet,K(A_\bullet)]$, where $K$ is the
Dold-Puppe functor, see \cite[II~4.11]{Qui67}.

In~\cite{Bes98a} the different
versions of syntomic cohomology are constructed as cohomologies of
bounded below complexes of presheaves $\Gamma_?^\bullet(i)$ on $\bcat'$,
where $?$ could stand for any of the versions of syntomic cohomology
considered. By~\cite[Proposition~\ref{syn-cech}]{Bes98a} these presheaves are
pseudo-flasque in the sense that there is a Mayer--Vietoris exact
sequence involving $ U $, $ V $, $X= U \bigcap V $ and
$ U \bigcup V $ for two open subsets $ U $ and $ V $ of $X$. We
associate to these presheaves spaces as follows: Consider first the
sheaf $ \tilde\Gamma_?^\bullet(i)' $ on $\bcat'$ associated to the
presheaf $ \Gamma_?^\bullet(i) $. Next consider the map of sites
$\smap: \bcat \to \bcat'$ corresponding to the inclusion
$\bcat'\subset \bcat$ and set $ \tilde\Gamma_?^\bullet(i)= \smap^\ast
 \tilde\Gamma_?^\bullet(i)'$, a complex of sheaves on $\bcat$.
As those complexes are complexes of sheaves with cohomological
numbering, we put, for any space $X_\bullet$,
$$
H^{2i-n}(X_\bullet, \Gamma_?^\bullet(i)) = [S^n \w X_\bullet, K(2i,
\tilde\Gamma_?^\bullet(i))]
$$
where,  for any cohomological complex $ A^\bullet $ in nonnegative degree,
$ K(2i, A^\bullet ) $ the Dold--Puppe construction applied to
the homological complex $ A^0 \to A^1 \to \dots \to A^{2i-1} \to
\ker(A^{2i} \to A^{2i+1})  $ 
in degrees $ 2i $ through 0.  
(We shall see shortly that, as the presheaves are pseudo flasque,
using the associated sheaves does not change much.)

The remainder of this Appendix mainly consists of getting an
explicit complex that computes $ H^{-m}(X_\bullet, \Gamma_?^\bullet(i)) $,
first for pointed simplicial schemes and then more specifically for
those pointed simplicial schemes underlying the
multi-relative $ K $--theory.  Together with the construction
of Chern classes in Proposition~\ref{chern-classes} below, this
provides the reference for the regulator and the complexes used
in the explicit calculations in the body of the paper.

So let $X_\bullet$ be a (smooth separated) pointed simplicial
scheme. If we write $ D_T $ for the derived category
of Abelian chain complexes in $ T $, 
$$
H^{2i-m}(X_\bullet, \Gamma_?^\bullet(i)) = [S^m \w X_\bullet, K(2i, \tilde\Gamma_?^\bullet(i)) ] 
$$
is also isomorphic to 
$ [ N_*(S^m \w X_\bullet), \tilde\Gamma_?^\bullet(i)]_{D_T} $,
cf. \cite[(24) on page 213]{Jeu95}, where $ N_*(\cdot) $ denotes
the reduced chain complex associated to the pointed simplicial
objects involved.
As in loc.\ cit., the Alexander--Whitney map,
in degree $n$ given by 
\begin{equation}\label{alexander-whitney}
X_n\wedge Y_n \mapsto
\sum_{i=0}^n d_{i+1}\cdots d_{n-1}d_n X_n\otimes d_0^{i}Y_n
\end{equation}
induces a quasi
isomorphism of $ N_*(X_\bullet \times Y_\bullet) $ with
$ N_*(X_\bullet) \tensor N_*(Y_\bullet) $. As $ N_*(S^n) $ is
quasi isomorphic to $ N_*(S^1) \tensor \dots \tensor N_*(S^1) $
($ n $ times) and $ N_*(S^1)  = \Z[-1] $ (a copy of $ \Z $ in
homological degree 1),
we find that we have to compute $ [ N_*(X_\bullet)[-n],\tilde\Gamma_?^\bullet(i)]_{D_T} $.
If $ \tilde\Gamma_?^\bullet(i) \to I^\bullet $ is an injective resolution,
this equals
$ [ N_*(X_\bullet)[-n],I^\bullet] $ (maps up to chain homotopy),
because the complex $ I^\bullet $ is bounded below (\cite[p.67]{Har66}).

Using the Yoneda lemma,  we can compute this as in~\cite[pages~214---216]{Jeu95} as the
$ 2i-n $--th cohomology of the
complex $ \calC^\bullet (X_\bullet, I^\bullet)$ given by 
\begin{equation}\label{first-complex}
\calC^q(X_\bullet, I^\bullet) = \bigoplus_{t+s = q}\hom(X_s,I^t) = \bigoplus_{s+t=q} \Gamma(X_s,I^t)
\end{equation}
with
$\dd^{(s,t)} = (-1)^q (\dd_s^{N_*(X_\bullet)})^* + \dd_t^{I^\bullet} $,
where of course we ignore the degenerate part of the 
scheme component of $X_s$ as well as the basepoint $*$
as we are working with the complex $N_*(X_\bullet)$.

We would like to be able to replace $ I^\bullet $ with $ \Gamma_?^\bullet(i) $
in this complex.  By using a filtration in the simplicial
index and the associated spectral sequence, it follows immediately
from the following two Lemmas that we can do this without changing
the cohomology of the complex.

\begin{lemma}\label{use-pseudoflasque}
Let $ P^\bullet $ be a complex of pseudo flasque presheaves
of Abelian groups on $\bcat'$
and let $ I^\bullet $ be an injective resolution
of the associated complex of sheaves. If $X\in \bcat'$, then the
natural map on global sections
$$
\Gamma (X, P^\bullet ) \to \Gamma (X, I^\bullet ) 
$$
is a quasi isomorphism of complexes of Abelian groups.
\end{lemma}
\begin{proof}
The proof follows the
proofs of Theorem~4 and Theorem~1' of \cite{BrGe73} extremely
closely, but it is easier as it is in the context of complexes of
Abelian groups rather than simplicial sets.
Namely, for every open set $ U $ of $ X $, let
$ F^\bullet(U) = \Cone(P^\bullet(U) \to I^\bullet(U)) $
Then it follows
immediately from the five lemma that the cohomology of $ F(U) $ satisfies a
Mayer--Vietoris exact sequence associated to two open subsets $ U $ and $ V $,
hence is pseudo flasque.   If we let $ T^q(U) $ be the $ q $--th
cohomology group of the complex $ F^\bullet(U) $, then the proof
of Theorem~1' applies verbatim if we replace $ * $ with $ 0 $
everywhere, and take into account that our indexing is cohomological
rather than homological.
\end{proof}
\begin{lemma}\label{use-smoothsep}
  Let $F^\bullet$ be a complex of sheaves on $\bcat'$ and let $X\in
  \bcat'$. Then the canonical map $H_{\bcat'}^i(X,F^\bullet)\to
  H_{\bcat}^i(X,\smap^\ast F^\bullet)$ is an isomorphism.
\end{lemma}
\begin{proof}
This is essentially part of~\cite[Proposition~III.3.1]{Mil80}, where
it is stated for a single sheaf and for the big site of schemes over
$X$, but the proof is the same.
\end{proof}

Because the complexes $ \Gamma_?^\bullet(i) $ are pseudo flasque,
by Lemma~\ref{use-pseudoflasque} and Lemma~\ref{use-smoothsep}
we could replace the complex~\eqref{first-complex} with the complex
$ \calC^\bullet(X_\bullet, \Gamma_?^\bullet(i)) $ given by
\begin{equation}\label{computable-complex}
\calC^q(X_\bullet, \Gamma_?^\bullet(i))
 = \bigoplus_{s+t=q} \Gamma(X_s,\Gamma_?^t(i))
\end{equation}
with
$\dd^{(s,t)} = (-1)^q (\dd_s^{N_*(X_\bullet)})^* + \dd_t^{\Gamma_?^\bullet(i)} $.
\begin{remark}\label{finalremark}
  Note that we are not using the fact that $X_\bullet$ is degenerate
  above a certain degree. Therefore, the description of cohomology as
  the cohomology of~\eqref{computable-complex} is valid for such
  spaces as well, provided the components belong to $\bcat'$. In
  particular, it is valid for the classifying
  spaces $BGL_n$ over $\vv$.
\end{remark}
In the paper, we have to use this complex on an iterated simplicial
reduced mapping cone $ C_\bullet $, which we recall is defined
as follows.

For $ f : X_\bullet \to Y_\bullet $ a map of pointed simplicial
schemes, define the reduced mapping cone of $ f $ by
\begin{equation}\label{mapcone}
C(X_\bullet,Y_\bullet) = Y_\bullet\coprod X_\bullet\times I/\sim,
\end{equation}
where $I$ is the simplicial version of the unit interval, given in degree
$s$ by all sequences $\{0,\dots,0,1\dots,1\}$ of length
$s+1,$ and pointed by $\{1,\dots,1\}$, and $\sim$ are the usual identifications
to obtain the reduced mapping cone.

Let $X$ be a scheme, and let $ Y_1,\dots,Y_n $, be
subschemes.  Denote by $ X_+ $ the pointed simplicial scheme
consisting of $ X $ in every degree, together with a disjoint
basepoint $ * $.
Consider the iterated mapping cone 
$ \Cr X Y 1 s $ inductively defined by
\begin{align*}
C(X,\{Y_1\}) & = C(Y_{1+},X_+)
\cr
\Cr X Y 1 s+1 & = C(\Cr Y_{s+1} Y 1,s+1 s,s+1 ,\Cr X Y 1 s ) 
\end{align*}
where $Y_{i,j}=Y_i\bigcap Y_j.$
Using induction, one sees easily that
the space $C_\bullet$ one finds for $X, Y_1,\dots,Y_n$ is as
follows. 
\begin{equation}\label{explicit-relativity}
C_m=
*{\scriptstyle\coprod}
\coprod_{\alpha_1,\dots,\alpha_n} Y_{\alpha_1,\dots,\alpha_n}
\end{equation}
with
$\alpha_i\in\{ \{0,\dots,0\} , \{0,\dots,0,1\} ,\dots,
\overbrace{ \{0,1,\dots,1\}}^{m+1}\}, $
$
Y_{\alpha_1,\dots,\alpha_s}=\bigcap_{\alpha_i\ne \{ 0,\dots,0\} } Y_i
$
and $\cap_\emptyset Y_i=X$.
The boundary and degeneracy maps are the natural maps coming from the
inclusions and the identity, which we get by deleting or doubling the
$i$-th place in the zeroes and ones, with the
convention that we identify $Y_{\alpha_1,\dots,\alpha_s}$ with $*$ if at
least one of the $\alpha$'s consists of only 1's.
Clearly,
$ C_\bullet $ is a pointed simplicial scheme, smooth if
$ X $, all $Y_j$ and all of their possible intersections are smooth.
Due to our definition of the mapping cone, the reduced chain
complex $ N_*(C_\bullet) $ no longer looks like an iterated
mapping cone of reduced chain complexes as there are too many
nondegenerate copies of intersections for $ n \geq 2 $, and neither
does the complex in (\ref{computable-complex}).

So we also define $ \ab C_{\bullet} $ to be the homological chain complex
given in degree $k$ by $ \coprod _{|\b|=k} \Z[ Y_\b ] $,
with $ \b $ a subset of $ \{ 1,\dots,n \} $,
$ Y_\b = \bigcap_{i \in \b} Y_i $, and $Y_\emptyset = X$.
The boundary is given on generators of $ \ab C_{\bullet} $
by
\begin{equation}\label{eq:explicit-differential}
d (Y_\b) = (-1)^{k-1} \sum_{j=1}^k (-1)^j Y_{\b\setminus \{\b_j\}}
\end{equation}
if $ \b = \{ \b_1,\dots,\b_k \} $ with $ \b_1<\b_2<\dots<\b_k $.
(Just as in the complexes $N_*(\cdot)$, the maps here are the
ones induced from the maps in the pointed spaces, which means the
correspond to the maps of sheaves that the scheme component
represent in our topos $T$.)

\begin{proposition}\label{explicit-qi}
$N_*(C_\bullet) $ and $ \ab {C_\bullet} $ are quasi isomorphic.
\end{proposition}

\begin{proof}
Define a map
$$
\xymatrix{
\Psi : \ab C_\bullet \ar[r] & N_*(C_\bullet) \\
}
$$
via
$$
Y_\b \mapsto \sum_{\s \in S_k} (-1)^\s Y_{(\b,\s)}
$$
in degree $ k $, where $ (-1)^\s $ is the sign of $ \s $, and
$ (\b,\s) = \a_1,\dots,\a_n  $ is an index defined as follows.
We make $\a_j = \{0,\dots,0\} $ unless $ j $ is an element of $\b$.
The remaining $k$ $\a_j$ are indexed by $\b$.
We consider the $ k $ standard $ k+1 $ tuples
$ \{ 0,\dots,0,1 \} $, 
$ \dots $,
$ \{ 0,1\dots,1 \} $,
and put $\a_{\b_{\s(j)}} $ equal to the $ j $-th $k+1$--tuple in
this list.

We have to check that $ \Psi $ defines a map of complexes.  This
is clear if $ k = 0 $.
For $ k \geq 1 $, 
$ \Psi_{k-1} \circ d $ is given by mapping $ Y_\b $ (with $ |\b|=k $)
to
$$
\Psi_{k-1}( (-1)^{k-1} \sum_{j=1}^k (-1)^j Y_{\b\setminus\{\b_j\}}) =
\sum_{j=1}^k (-1)^{k+j-1} \sum_{\tau \in S_{k-1}} (-1)^\tau Y_{(\b\setminus\{\b_j\},\tau)}
.
$$
On the other hand, $ \Psi_k $ maps $ Y_\b $ to 
$ \sum_{\s \in S_k} (-1)^\s Y_{(\b,\s)} $, which $ d $ maps to
$$
\sum_{\s \in S_k} (-1)^\s \sum_{j=0}^k (-1)^j Y_{d_j(\b,\s)}
$$
where $ d_j $ is the $ j $--th simplicial face.
Now notice that the $ j=0 $ term here is zero, as it introduces
$ \{1,\dots,1\} $ among the indices so this corresponds to $ * $
in $ C_\bullet $, which maps to zero in $ N_*(C_\bullet) $.  Also, for $ j=1,\dots,k-1 $
the $ j $--th and $ j+1 $--st indices become the same after applying $ d_j $, so that
$ \s $ and $ \s \circ (j \, j+1) $ give the same contribution,
which cancels due to $ (-1)^\s $.  Therefore only one term survives,
corresponding to $ j=k $, i.e., $ d_k $, which eliminates the last element of
all the indices.  So we are left with 
\begin{equation*}
(-1)^k \sum_{\s \in S_k} (-1)^\s  Y_{d_k(\b,\s)} =
(-1)^k \sum _{j=1}^k \sum_{\tau \in S_{k-1}} (-1)^{j-1} (-1)^{\tau} Y_{( \b-\{\b_j\}, \tau) }
\end{equation*}
because if $ \s(1) = j $, the $ j $--th index in $ d_k\b $ consists
only of zeroes as well, and we can write the corresponding term
as a contribution in the sum on the right with $ k =j  $
and $  \tau = (k \dots j) \circ \s \circ (1 \dots k) $ due to
the renumbering involved.

In order to check that $ \Psi $ defines a quasi isomorphism,
we proceed by induction on the degree of relativity $ n $.
We investigate how $ \Psi $ behaves with respect to taking cones.

So let $ Y_\bullet \to X_\bullet $
correspond to taking the last
($ n $--th) relativity into account, with $ C_\bullet $ the corresponding
reduced mapping cone.
Let us first note that $ \ab C_{\bullet} $ is the cone
of the map $ \ab Y_{\bullet} \to \ab X_{\bullet} $.
Namely, any component $ C_{\b} $ of $ \ab C_{\bullet} $ comes from
$ \ab X_{\bullet} $ if and only if $ \b $ does not contain $ n $,
and that the components containing $ n $ correspond to 
$ Y_{\b\setminus\{n\}} $'s, i.e., to $ \ab Y_{\bullet} [-1] $.
So $ \ab C_{\bullet} $ is the
mapping cone of $ \ab Y_{\bullet} \to \ab X_{\bullet}  $ provided
the differential is the one on the cone.  As $ \ab X_{\bullet} $ is a subcomplex,
we only need to check what the differential does on $  C_\b $
with $ n $ in $ \b $.
Let $ k= |\b| $.
Applying $ d $ we get
$$ (-1)^{k-1} \sum_{j=1}^k (-1)^j C_{\b\setminus\{\b_j\}}  = 
-C_{\b\setminus\{n\}} - (-1)^{k-2} \sum_{j=1}^{k-1} (-1)^j C_{\b\setminus\{\b_j\}} ,$$
which is exactly what we want.

We shall verify that we have a map of triangles
\begin{equation*}
\xymatrix{
\ar[r] & \ab Y_{\bullet} \ar[r]\ar[d]_\Psi & \ab X_{\bullet} \ar[r]\ar[d]_\Psi & \ab C_{\bullet} \ar[d]_\Psi \ar[r] & \ab Y_{\bullet} [-1] \ar[r]\ar[d]_{\Psi[-1]}&
\\
\ar[r] & N_*(Y_{\bullet}) \ar[r]      & N_*(X_{\bullet}) \ar[r]      & N_*(C_\bullet)
 \ar[r] & N_*(Y_{\bullet})[-1] \ar[r] &
}
\end{equation*}
with the maps as follows.
The map $ \ab Y_{\bullet} \to \ab X_{\bullet} $ comes directly
from $ Y_\bullet \to X_\bullet $, and similarly for $ N_*(Y_\bullet) \to N_*(X_\bullet) $.
The map $ \ab X_{\bullet} \to \ab C_{\bullet} $ 
just views $ \b $, a subset of $ \{1,\dots,n-1\} $, as a subset of $ \{1,\dots,n\} $.  
The map $ N_*(X_\bullet) \to N_*(C_\bullet) $ is the natural
map from the map $ X_\bullet \to C_\bullet $.
The map $ \ab C_{\bullet} \to \ab Y_{\bullet} [-1] $ maps $ Y_\b $
to 0 if $ n $ is not in $ \b $, and to $ Y_{\b\setminus\{n\}} $
if $ n $ is in $ \b $.
The map $ N_*(C_{\bullet}) \to N_*(Y_{\bullet}) [-1] $ is
the composition of the natural map $ N_*(C_{\bullet}) \to N_*(S^1 \w Y_{\bullet}) $
corresponding to contracting $ X_{\bullet} $ to $ * $,
and the Alexander--Whitney map $ N_*(S^1 \w Y_{\bullet}) \to N_*(Y_{\bullet}) [-1] $,
a quasi isomorphism (see (\ref{alexander-whitney})) as $ N_*(S^1) $
is isomorphic to $ \Z[-1] $ via the projection to the $ \{ 0,1\} $--component
in degree 1.

We shall check below that those maps give a map of triangles.
It is well known (and using Mayer--Vietoris for 
an open cover $ U \bigcup V $ for nonreduced mapping cones of the realization
of simplicial sets easily seen) that the bottom triangle  gives
rise to a long exact sequence in homology.  As the top triangle
also gives a long exact sequence, we know by induction that
$ \Psi : \ab C_{\bullet} \to N_*(C_\bullet) $ is a quasi isomorphism,
as $ \Psi $ is clearly an isomorphism if $ n=0 $.

In the diagram, the first square commutes because of the naturality
of $ \Psi $.  For the second square, we note that
applying $ \Psi_{C_\bullet} $ on the image
of $ \ab X_{\bullet} $ inside  $ \ab C_{\bullet} $ is the
same as applying $ \Psi_{X_\bullet} $ and tagging on
an index $ \{0,\dots,0\} $ to the indices already used in 
$ N_*(X_\bullet) $.  (The extra $ \{0,\dots,0\} $
corresponds to $ n \in \{1,\dots,n\} $.)
This is exactly the result as going around the second square counterclockwise,
as $ X_\bullet $ is the 
simplicial subscheme of $ C_{\bullet} $, given by the components
of $ C_{\bullet} $ that acquire a copy of $ \{0,\dots,0\} $  from
the simplicial interval involved in constructing the mapping
cone.

For the third square, the map $ \ab C_{\bullet} \to \ab Y_{\bullet} [-1] $
corresponds to mapping to zero any $ Y_\b $ with
$ \b \subseteq \{1,\dots,n-1\}\subset\{1,\dots,n\} $,
and to $ Y_{\b \setminus\{n\}} $ if $ \b \not \subseteq \{1,\dots,n-1\}  $.
If $ n \not\in \b $, $ \Psi(Y_\b) $ will always have the last index
$ \a_n $ in $ Y_{\a_1,\dots,\a_n} $ equal to $ \{0,\dots,0\} $,
which already goes to zero in $ N_*(S^1\w Y_\bullet) $.
If we have a term $ Y_\b $ with $ n \in \b $,  let $ k = |\b| \geq 1 $.
Going clockwise, we get
$ \Psi_{k-1}[-1](Y_{\b\setminus\{n\}}) = \sum_{\s \in S_{k-1} } (-1)^\s Y_{( \b\setminus\{n\} , \s )} $
in $ N_*(Y_\bullet)[-1] $.  
Going in the other direction, we get
\begin{equation*}
Y_{\b}  \mapsto \sum_{\s \in S_k} (-1)^\s Y_{(\b,\s)}
\end{equation*}
in $ N_*(C_\bullet) $.
The Alexander--Whitney map (\ref{alexander-whitney}) maps this to
$$
\sum_{j=0}^k \sum_{\s \in S_k} 
    (-1)^\s d_{j+1}\circ d_{j+2} \circ \dots \circ d_k \left( (\b,\s)_n \right)
       \tensor d_0^j Y_{  (\b,\s)\setminus\{(\b_n,\s)_n\} } 
.
$$
After the projection to $ \{0,1\} $ in $ N_*(S^1) $, 
we only get a nonzero contribution if $ (\b, \s)_n= \{0,1,\dots,1\} $
and $ j=1 $.  Considering the definition of $ (\b,\s) $,
this means that $ \s(n) = n $, so that $ \s \in S_{k-1} \subset S_k $.
So we find
$$
\sum_{\s \in S_{k-1}} (-1)^\s \{0,1\} \tensor d_0 Y_{(\gamma,\s)} 
,
$$
with $ \gamma = (\b,\s) \setminus \{ (\b,\s)_n \} $.
Considering that $ d_0 $  removes the first coordinate in all
the first $ n-1 $ tuples, $  d_0 Y_{(\gamma,\s)}  = Y_{(\b\setminus\{n\}, \s)}$.
Projecting to the $ \{0,1\} $--component in $ N_*(S^1) = \Z[-1] $,
we therefore find
$$
 \sum_{\s \in S_{k-1}} (-1)^\s Y_{(\b \setminus\{n\},\s)}
$$
as before.

\end{proof}

We now return to our original problem of computing the cohomology
groups using explicit complexes.  For $ C_\bullet =  \Cr X Y 1 n $ as above,
if all scheme components are smooth of finite type over the base
ring $ R $, we can replace $ N_*(C_\bullet) $ by the quasi isomorphic
complex $ \ab {C_\bullet} $ by Proposition~\ref{explicit-qi}
from the very beginning, so instead of~\eqref{first-complex} or
\eqref{computable-complex}, we can
also use the complex $ \calfinalC^\bullet(C_\bullet, \Gamma_?^\bullet(i)) $
with
\begin{equation}\label{final-complex}
\calfinalC^q(C_\bullet, \Gamma_?^\bullet)
= \bigoplus_{t+s = q} \oplus_{|\b|= s }\Gamma(X_{\b},\Gamma_?^t(i))
\end{equation}
and $\dd^{(s,t)} = (-1)^q (\dd_s^{ \ab C_{\bullet} } )^*  + \dd_t^{\Gamma_?^\bullet(i)} $.

\begin{remark}\label{cupcompatible}
In order to get products in $ K $--theory taking the relativity
into account, we define maps 
$$
\Cr X Y 1 s+t  \to \Cr X Y 1 s \w \Cr X Y s+1 s+t 
$$
by the diagonal embedding
$ Y_{\a_1,\dots,\a_{r+s}} \to Y_{\a_1,\dots,\a_s} \times Y_{\a_{s+1},\dots,\a_{s+t}} $
and identifying  anything of the form $ \dots \times *  $
or $  * \times \dots $ with $ * $ in the right hand side.
Taking reduced chain complexes, and using the Alexander--Whitney
map gives us a map 
$$
N_*( \Cr X Y 1 s+t ) \to N_*( \Cr X Y 1 s ) \tensor N_*( \Cr X Y s+1 t )  
$$
which we want to compare with a similar map using the $ \ab {\cdot} $
complexes.
Namely, let $ Y_\b $ be a component in $ \ab {\Cr X Y 1 s+t } $, and 
let $ \b_1 = \{1,\dots,s\} \bigcap \b $, $ \b_2 = \{s+1,\dots,s+t\} \bigcap \b $.
Then we define the map 
$$
\ab {\Cr X Y 1 s+t } \to  \ab {\Cr X Y 1 s }  \tensor \ab {\Cr X Y s+1 s+t } 
$$
via the map $  Y_\b \mapsto (-1)^ {|\b_1|\cdot |\b_2|} Y_{\b_1} \tensor Y_{\b_2} $.
(This is a map of formal sums of sheaves represented by schemes,
and the maps are induced from the scheme embeddings $ Y_\b \to Y_{\b_1} $
and $ Y_\b \to Y_{\b_2} $.) 
We claim that with this definition we have a commutative diagram
(where we suppress $X$ and $Y$ from the notation for typographical reasons)
\begin{equation*}
\xymatrix{
\ab {\Crshort X Y 1 s+t } \ar[r]\ar[d]^{\Psi} & \ab {\Crshort X Y 1 s }  \tensor \ab {\Crshort X Y s+1 s+t }\ar[d]^{\Psi\tensor\Psi}
\\
N_*( {\Crshort X Y 1 s+t } ) \ar[r] & N_*( {\Crshort X Y 1 s ) }  \tensor N_*( {\Crshort X Y s+1 s+t ) }
}
\end{equation*}
Namely, write $ k = |\b| $, $ k_1 = |\b_1| $ and $ k_2 = |\b_2| $.
Starting in the top left corner of the diagram, going the bottom
left corner, and then to the bottom right corner (using also
the Alexander--Whitney map), we get on
$ Y_\b $:
\begin{align*}
Y_\b
& \mapsto
\sum_{\s\in S_k} (-1)^\s Y_{(\b,\s)} 
\\
& \mapsto
\sum_{\s \in S_k} (-1)^\s
                Y_{(\b,\s)_1,\dots,(\b,\s)_s } \tensor Y_{(\b,\s)_{s+1},\dots,(\b,\s)_{s+t} }
\\
& \mapsto
\sum_{j=0}^k \sum_{\s \in S_k} (-1)^\s d_{j+1} d_{j+2} \dots d_k
                Y_{(\b,\s)_1,\dots,(\b,\s)_s } \tensor d_0^j Y_{(\b,\s)_{s+1},\dots,(\b,\s)_{s+t} }
.
\end{align*}
Now observe that the nonzero indices involved are $ k $ in total,
of length $ k+1 $, i.e.,
\begin{align*}
& \{0,0,\dots,0,0,0,\dots,0,1\} \\
& \qquad\qquad\quad\, \vdots \\
& \{0,0,\dots,0,0,1,\dots,1,1\} \\
& \{0,0,\dots,0,1,1,\dots,1,1\} \\
& \{0,0,\dots,1,1,1,\dots,1,1\} \\
& \qquad\qquad\quad\, \vdots  \\
& \{ \underbrace{0,1,\dots,1}_{0 \dots j-1},\underbrace{1,1,\dots,1,1}_{j \dots k}\} \smash{.}
\end{align*}
$ d_{j+1} d_{j+2} \dots d_k$ deletes the last $ k-j $
columns, $ d_0^j $ deletes the first $ j $ columns of all tuples
involved.
For fixed $ j $, if any of the last $ j $ tuples ends up among
the last $ t $ tuples of $ (\b,\s) $,
then one of the
tuples becomes $ \{1,\dots,1\} $ under $ d_0^j $ and the corresponding
component is mapped to zero in $ N_* $.
So for a nonzero contribution, the last $ j $ tuples must end up among
$ (\b,\s)_1,\dots,(\b,\s)_s $.
If any more tuples end up in $ (\b,\s)_1,\dots,(\b,\s)_s $,
then after applying $ d_{j+1} d_{j+2} \dots d_k $, we get at
least two tuples consisting entirely of zeroes.
But as these contributions are summed alternatingly over $ S_k $,
and the original indices can be swapped by a transposition which
yields a minus sign, those contributions cancel.   As $ k_1 $
tuples must end up among $ (\b,\s)_1,\dots,(\b,\s)_s $, this
shows that for a nonzero contribution we must have
$ j = k_1 $, the nonzero tuples among $ (\b,\s)_1,\dots,(\b,\s)_s $
are the last $ k_1 $ rows above, and the nonzero tuples among
$ (\b,\s)_{s+1},\dots,(\b,\s)_{s+t} $ are the first $ k - k_1 = k_2 $
rows above.
The sum then simplifies to 
\begin{align*}
& \phantom{=}\,
\sum_{\s \in S_k} (-1)^\s d_{k_1+1} d_{k_1+2} \dots d_k Y_{ (\b,\s)_1,\dots,(\b,\s)_s}
                                        \tensor d_0^{k_1} Y_{(\b, \s)_{s+1},\dots,(\b,\s)_{s+t}}
\\
& =
(-1)^{k_1 k_2} \left( \sum_{\tau_1 \in S_{k_1}} (-1)^{\tau_1} Y_{(\b_1,\tau_1)} \right) \tensor 
                     \left( \sum_{\tau_2 \in S_{k_2}} (-1)^{\tau_2} Y_{(\b_2,\tau_2)} \right)
\end{align*}
because the permutation $ \s $ must map $ \{1,\dots,k-k_1\} $
to $ \{ k_1+1,\dots,k\} $ as well as $ \{k-k_1+1,\dots,k\} $ to $ \{1,\dots,k_1\} $,
so we must have 
$$
(1\dots k)^{-k_1} \s =  \tau_2 (1 \, k_2+1)\dots (k_1 \, k) \tau_1 (1 \, k_2+1)\dots (k_1 \, k)
$$
for some $ \tau_1 $ in $ S_{k_1} $ and $ \tau_2 $ in $ S_{k_2} $.
As this equals $ (-1)^{|\b_1|\cdot |\b_2|} \Psi(Y_{\b_1}) \tensor \Psi(Y_{\b_2}) $
the diagram commutes as required.

Now suppose that $ A $ and $ B $ are homological chain complexes,
with a cup product $  A \tensor B \to C $.  There is a map
$ \phi : K(A) \w K(B) \to K(A \tensor B) $
(with $ K $ the Dold--Puppe construction as before), which 
gives rise to a map 
\begin{align*}
[S^n \w X_\bullet, K(A) ] \times [S^m \w Y_\bullet, K(B)]
& \to
[S^{m+n} \w X_\bullet \w Y_\bullet, K(A \tensor B) ]
\cr
& \to [S^{m+n}\w
X_\bullet \w Y_\bullet \to K(C)]
.
\end{align*}
It is shown on \cite[page~215]{Jeu95} that under our identifications
$ [S^n \w X_\bullet, K(A) ]  $ with $ [N_*(X_\bullet)[-n], A]_{D_T} $
etc.,
this corresponds to the composition
$$
N_*( X_\bullet \w Y_\bullet )[-n-m] \to
N_*(X_\bullet)[-n] \tensor N_*(Y_\bullet)[-m] \to A \tensor B \to C
$$
with the first map the Alexander--Whitney map and the last map
the given product. 
In the cases we are interested in this becomes a cup product
of sections in (pre)sheaves, and it follows from these formulae
that the product on components corresponds to cup products
$ \Gamma(X_s,A) \times \Gamma(Y_t,B) \to \Gamma(X_s \times Y_t,C) $
up to signs.  In particular, for the explicit map at the very
beginning of this Remark, the diagram tells us that the product
in the complex (\ref{final-complex})
is up to a sign given by the map
$$
\Gamma(Y_{\b_1},A) \times \Gamma(Y_{\b_2},B) \to \Gamma(Y_{\b1\b2},C)
,
$$
which is the composition of
$$
\Gamma(Y_{\b_1},A) \times \Gamma(Y_{\b_2},B) \to \Gamma(Y_{\b1}\times
Y_{\b2},C) \to \Gamma(Y_{\b1\b2},C)
,
$$
the last map being the pullback corresponding to the \lq\lq diagonal\rq\rq\
$  Y_{\b_1\b_2} \to Y_{\b_2} \times Y_{\b_2} $.
\end{remark}

The previous constructions are applied in the body of the paper with the 
schemes 
$$
X = (\PP_B^1 \setminus \{t=1\})^n \text{ and } Y_i = \{t_i=0,\infty\} 
,
$$
with $ t $ the standard affine coordinate on $ \PP^1 $, and $t_i$
the $i$--th coordinate in the $n$--fold product, or localizations
of those schemes.

After this rather explicit exercise, we now turn our attention to the theory of Chern classes.
The theory of syntomic Chern classes of~\cite{Bes98a} can be
extended 
from schemes to arbitrary spaces as follows. In loc.\ cit.\ before
Theorem~\ref{7.5}
universal Chern classes
\begin{equation*}
  c_n\in H^{2i}(BGL,\Gamma^\bullet_?(i))
\end{equation*}
were constructed. Again this was
explicitly done in some of the theories but it can easily be done
in all the others. Further, this was done with the cohomology defined
as the cohomology of the complex \eqref{computable-complex}, but since
the components of $BGL_n$ belong to $\bcat'$ it follows from
Remark~\ref{finalremark} that this is
the same as the definition we have been using here.
Now a standard procedure~\cite[6.1]{Gil-Sou99}
produces, for each $\alpha\in H^{-m}(X_\bullet,K)$ a Chern class
\begin{equation*}
  c_i(\alpha)\in H^{2i-m}(X_\bullet,\Gamma^\bullet_?(i))\; .
\end{equation*}
More precisely, if $ K = \Z \times \Z_\infty BGL $ is the sheaf used to define algebraic
$ K $-theory of spaces as $ K_n(X_\bullet) = [S^n \w X_\bullet, K] $,
then each $ c_i $ defines a map
$ K \to K(2i, \tilde\Gamma_?^\bullet(i)) $.
If $ \a $ is an element in $ K_m(X_\bullet) $, then by composition we get
the element $ c_i(\a) $ in
$ [S^m \w X_\bullet, K(2i, \tilde\Gamma_?^\bullet(i)) ] =
H^{2i-m}(X_\bullet, \Gamma_?^\bullet(i))$.

For a $K$-coherent space $X_\bullet$, both
\begin{equation*}
  H^\sim(X_\bullet,\Gamma^\bullet_?) := H^0(X_\bullet,\Z)\times (\{1\}\times \left(\oplus_{i>1}
  H(X_\bullet,\Gamma^\bullet_?(i))\right))
\end{equation*}
and $\oplus_{m\ge 0} H^{-m}(X_\bullet,K)$ have $\lambda$-ring with involution
structures
described in loc.\ cit.\ 6.1. and there is a total Chern class
\begin{equation*}
  c:\oplus_{m\ge 0} H^{-m}(X_\bullet,K) \to  H^\sim(X_\bullet,\Gamma^\bullet_?)\; .
\end{equation*}

\begin{proposition}\label{chern-classes}
  When $X_\bullet$ is $K$-coherent the total Chern class in a morphism of
  $\lambda$-rings with involutions.
\end{proposition}

\begin{proof}
Everything is reduced to the properties of the universal Chern classes
(see for example the proof of~\cite[Theorem 5]{Gil-Sou99} for the
$\lambda$-structure). These properties are deduced in the following
way. There is a map of complexes of sheaves (in the derived category)
$\Gamma^\bullet_?(n) 
\to \Gamma^\bullet_{\dr}$, where the latter complex is the complex of
differential forms on the generic fiber. We get an induced map of
cohomology theories which is compatible with cup products and
therefore also with $\lambda$-operations. This map
gives an injection
\begin{equation}\label{inji}
  \oplus_i H^{2i}(BGL_N,\Gamma^\bullet_?(i))\inject \oplus_i\hdr^{2i}(BGL_N/K)
\end{equation}
on the part of the 
cohomology of $BGL_N$ containing the Chern classes for any $N$.
The syntomic universal Chern classes are defined to map to the
corresponding de Rham Chern classes. Since both sides of \eqref{inji}
are closed under products, all required properties of
syntomic universal Chern  classes follow from the corresponding results
for the universal de Rham classes.
\end{proof}

As all the cohomology groups are $ \Q $--vector spaces, 
one gets a Chern character from this in the usual way
(cf. \cite[\S4]{Schn88} or \cite[Definition~2.34]{Gil81}),
which gives a ring homomorphism
$$
\reg : K_*(X_\bullet) = H^{-*}(X_\bullet,K) \to  H^*(X_\bullet,\Gamma^\bullet_?(*))\; .
$$
with the property that
$
\reg (K_m^{(j)}(X_\bullet,K)) \subseteq 
    H^{2j-m}(X_\bullet,\Gamma^\bullet_?(j))
,
$ 
cf. \cite[Corollary on page 28]{Schn88}.

\end{document}